%% file: Klein-11-Second-Version.tex
\documentclass{article}

\advance\textheight -\headheight
\advance\textheight -\headsep
\oddsidemargin\paperwidth
\advance\oddsidemargin -\textwidth
\divide\oddsidemargin2 
\ifdim\oddsidemargin<.5truein \oddsidemargin.5truein \fi
\advance\oddsidemargin -1truein
\evensidemargin\oddsidemargin
\topmargin\paperheight \advance\topmargin -\textheight
\advance\topmargin -\headheight \advance\topmargin -\headsep
\divide\topmargin2
\ifdim\topmargin<.5truein \topmargin.5truein \fi
\advance\topmargin -1.4truein\relax

\usepackage{graphicx,amssymb,latexsym,amsfonts,txfonts}
\usepackage{pdfsync,color,tabularx,rotating}
\usepackage[all,cmtip]{xy}
\usepackage{url}
\usepackage{tikz}
\usepackage{amssymb}
\usepackage{comment}
\usepackage{epsfig}

\newcommand{\M}{\mathcal{M}}

\newcommand{\PP}{\mathbb{P}}
\newcommand{\F}{\mathbb{F}}

\newtheorem{theo}{Theorem}[section]

\newtheorem{coro}[theo]{Corollary}
\newtheorem{lemm}[theo]{Lemma}
\newtheorem{conj}[theo]{Conjecture}

\newtheorem{rema}[theo]{\bf Remark}
\newtheorem{exam}[theo]{\bf Example}

\newtheorem{conv}[theo]{\bf Convention}

\begin{document}

\title{Klein's ten planar dessins of degree $11$, \\ and beyond}

\author{Gareth A. Jones and Alexander K. Zvonkin}





\maketitle

\begin{abstract}
We reinterpret ideas in Klein's paper on transformations of degree~$11$ from the modern point of view of dessins d'enfants, and extend his results by considering dessins of type $(3,2,p)$ and degree $p$ or $p+1$, where $p$ is prime. In many cases we determine the passports and monodromy groups of these dessins, and in a few small cases we give drawings which are topologically (or, in certain examples, even geometrically) correct. We use the Bateman--Horn Conjecture and extensive computer searches to support a conjecture that there are infinitely many primes of the form $p=(q^n-1)/(q-1)$ for some prime power $q$, in which case infinitely many groups ${\rm PSL}_n(q)$ arise as permutation groups and monodromy groups of degree $p$ (an open problem in group theory).
\end{abstract}


\medskip


\section{Introduction}
\label{sec:Intro}

In 1878--9 Klein published two papers~\cite{Kle78, Kle79} with almost identical titles. The first concerned equations of degree~$7$, while the second concerned those of degree~$11$, with the properties and actions of the groups ${\rm PSL}_2(p)$, for $p=7$ and $11$ respectively, playing a major role in them. The paper~\cite{Kle78} is deservedly famous, not least for introducing the {\em Hauptfigur}, the iconic diagram of a tessellation of a hyperbolic $14$-gon which, after identification of sides, yields the quartic curve of genus~$3$ now named after Klein: see~\cite{Lev}, a book entirely devoted to this curve.  In contrast with the largely geometric approach taken in~\cite{Kle78}, Klein's emphasis in~\cite{Kle79} was more algebraic. However, that paper contained a diagram\footnote{Regrettably this diagram, like that of the {\em Hauptfigur} in~\cite{Kle78}, is omitted from the on-line version of the paper, both having been provided on separate inserted sheets in the original journal. However, both diagrams can be found in Klein's collected works~\cite[vol.~3, pp.~126, 143]{KleGMA}, integrated into updated versions of the two papers; see also~\cite[pp.~115, 320]{Lev} for the {\em Hauptfigur}}. showing ten plane trees, which anticipated Grothendieck's concept of dessins d'enfants~\cite{Gro} by over a century. Our aim here is to reinterpret this diagram and the mathematics underlying it from a modern point of view, and to explore some generalisations which Klein might have been able to achieve if he had had access to such tools as character theory, computer algebra and the classification of finite simple groups.

In particular, we study dessins of type $(3,2,p)$ and of degree $p$ or $p+1$ for arbitrary primes $p$, extending Klein's work in~\cite{Kle79} on the case $p=11$. These dessins are best studied through their monodromy groups, which are permutation groups of these degrees, and in principle such groups are all known, as a consequence of the classification of finite simple groups. However, because of certain open problems in number theory it is unknown whether two families of these groups, consisting of projective and affine groups of degrees $p$ and $p+1$, are finite or infinite. Using the Bateman--Horn Conjecture and extensive computer searches, we present what we believe is strong evidence that the first of these two families is infinite. (The second family is based on the Mersenne primes, and we have nothing new to add about its cardinality.)

In Section~\ref{sec:klein-trees} we describe Klein's ten plane trees, firstly as he presented them in~\cite{Kle79}, as pictorial representions of the possible 11-sheeted coverings of the Riemann sphere  $\Sigma={\mathbb P}^1({\mathbb C})$ with branching patterns $3^31^2$, $2^41^3$ and $11^1$ over $0$, $1$ and $\infty$. We then use more modern terminology and graphic conventions to describe them as dessins d'enfants, specifically as bipartite maps (Figure~\ref{L2(11)dessins}) and as maps with free edges (Figure~\ref{Kleinmaps}). Klein asserted without proof that his list was complete, so we outline a method he might have used to see this.

An important invariant of any dessin is its monodromy group (equivalently the monodromy group of the covering of the sphere which it represents), discussed in Section~\ref{sec:monodromy}. This is a transitive permutation group on the sheets of the covering (more precisely, on the fibre over a base-point), generated by the monodromy permutations describing the branching over the critical values. Klein showed that just one chiral pair of his trees, corresponding to our maps ${\mathcal M}_1$ and $\overline{\mathcal M}_1$, have monodromy groups isomorphic to 
${\rm PSL}_2(11)$, acting with degree $11$; we show that the other eight have the alternating group ${\rm A}_{11}$ as their monodromy group.

A powerful technique for enumerating dessins with a given monodromy group, but unavailable to Klein when he wrote his papers, is the character-theoretic triple-counting formula introduced by Frobenius~\cite{Fro} in 1896. We use this in Section~\ref{triples}, with the aid of character tables provided by GAP~\cite{GAP} and by the ATLAS~\cite{ATLAS}, to confirm Klein's enumeration of the trees and our description of their monodromy groups; this includes the elimination of the Mathieu group ${\rm M}_{11}$, another possible candidate.

Any dessin $\mathcal D$ with monodromy group $G$ has a unique minimal regular cover, a regular dessin $\mathcal R$ with ${\rm Aut}\,{\mathcal R}\cong G$. Much of Klein's paper~\cite{Kle79} is devoted to the common regular cover ${\mathcal R}_1$ of ${\mathcal M}_1$ and $\overline{\mathcal M}_1$, a map on the modular curve $X(11)$ of genus $26$. We show in Section~\ref{regcov} that the other eight trees, with $G\cong{\rm A}_{11}$, have mutually non-isomorphic regular covers of genus $756\,001$.

Each of Klein's ten trees corresponds to an $11$-sheeted covering $\Sigma\to\Sigma$ realised by a Shabat polynomial $P$ of degree $11$. In Section~\ref{sec:CompGal} we determine these polynomials for ${\mathcal M}_1$ and $\overline{\mathcal M}_1$: they are defined over ${\mathbb Q}(\sqrt{-11})$, and are transposed by the Galois group of that field, while those for the other eight trees form a single Galois orbit. Whereas the illustrations of dessins given so far are only topologically correct, we compute and show in Figure~\ref{Klein-tree-new} a geometrically correct drawing of ${\mathcal M}_1$ as the inverse image under $P$ of the unit interval.

In Section~\ref{sec:(3,2,r)} we extend Klein's classification by determining all the dessins of degree $11$ and type $(3,2,r)$ for $r\ne 11$: there are sixteen of them (see Figure~\ref{othermaps}), all of genus $0$ and with the symmetric group ${\rm S}_{11}$ as their monodromy group. 

In Sections~\ref{sec:mono(3,2,p)}, \ref{sec:BHC} and \ref{sec:dessins(3,2,p)} we generalise Klein's work further by considering dessins $\mathcal D$ of degree $p$ and type $(3,2,p)$ for arbitrary primes $p$. The investigation is quite difficult, because we encounter some important open problems in group theory and number theory. The monodromy group $G$ of $\mathcal D$ must be a transitive permutation group of degree $p$, and if $p>3$ (as we shall assume) then as a quotient of the triangle group $\Delta(3,2,p)$ with mutually coprime periods, $G$ must be perfect. In Section~\ref{sec:mono(3,2,p)} we use the classification of finite simple groups, and hence of permutation groups of prime degree, to restrict the possibilities for $G$ to the alternating groups ${\rm A}_p$, ${\rm PSL}_2(11)$ and the Mathieu groups ${\rm M}_{11}$ and ${\rm M}_{23}$ for $p=11$, $11$ and $23$, and ${\rm PSL}_n(q)$ in cases where its natural degree $(q^n-1)/(q-1)$ is a prime $p$. We show that ${\rm A}_p$ arises for $p=5$ and each $p\ge11$, with exponentially many dessins (far too many to classify) as $p\to\infty$. The group ${\rm PSL}_2(11)$, acting with degree~$11$, was dealt with by Klein, as we have already seen, while the Frobenius triple-counting formula eliminates the Mathieu groups as quotients of $\Delta(3,2,p)$. 

This leaves the groups ${\rm PSL}_n(q)$, acting naturally with degree $(q^n-1)/(q-1)$ for some prime power $q$. In Section~\ref{sec:BHC} we consider the open problem of whether this degree is prime for finitely or infinitely many pairs $(n,q)$; such {\sl projective primes}, as we call them, include the Fermat and Mersenne primes, for $n=2$ and $q=2$ respectively. For each $n\ge 3$ (necessarily prime) the Bateman--Horn Conjecture (a powerful but apparently little-known open problem in number theory) provides a heuristic estimate for the number of primes 
$t\le x$ such that $q=t^e$ and
$(q^n-1)/(q-1)$ is prime; since these estimates agree closely with results obtained from computer searches revealing large numbers of examples, we conjecture that for each prime $n\ge 3$ there are infinitely many projective primes $p$. This case is explored in more detail in~\cite{JZ20}.

In Section~\ref{sec:dessins(3,2,p)} we return to the dessins of type $(3,2,p)$ and degree~$p$ for projective primes $p=(q^n-1)/(q-1)$. Using character theory we show that in the cases $n=2$ and $3$ there are such dessins with monodromy group ${\rm PSL}_n(q)$ for each projective prime $p$; in Theorems~\ref{Fermat} and \ref{th:PSL3q} we enumerate them and give their genus. Theorem~\ref{th:Mersenne} is a similar but weaker result for $q=2$, proving the existence of such dessins for each Mersenne prime $p=2^n-1$, and we conjecture that it extends to all projective primes.

In Sections~\ref{sec:modular} and \ref{sec:dessinsp+1} we consider dessins of type $(3,2,p)$ and degree $p+1$ for arbitrary primes~$p$. The parameters $3$ and $2$ in the types of the dessins we have considered imply that all the monodromy groups and automorphism groups associated with them are quotients of the modular group ${\rm PSL}_2({\mathbb Z})\cong{\rm C}_3*{\rm C}_2$. In most cases the corresponding kernels are non-congruence subgroups, but in Section~\ref{sec:modular} we consider the modular dessins ${\mathcal D}_0(p)$ arising from reduction of coefficients mod~$(p)$, corresponding to the action of ${\rm PSL}_2(p)$ on the projective line ${\mathbb P}^1({\mathbb F}_p)$. We do this firstly in the planar cases $p\le 7$ and $p=13$, and then in the more difficult case $p=11$, where we outline a calculation by John Voight which determines both the underlying elliptic curve $E$ and the Bely\u\i\/ function $E\to\Sigma$. For $p=11$ and $13$ we compute and display geometrically correct drawings of these dessins.

We extend this further in Section~\ref{sec:dessinsp+1} by considering arbitrary dessins of type $(3,2,p)$ and degree $p+1$ for primes $p>3$. We use the classification by M\"uller~\cite{Mul} of primitive permutation groups containing a cycle with one fixed point to restrict the possible monodromy groups to ${\rm A}_{p+1}$, ${\rm PSL}_2(p)$, ${\rm AGL}_n(2)$ for $p=2^n-1$, ${\rm M}_{11}$ and ${\rm M}_{12}$ for $p=11$, and ${\rm M}_{24}$ for $p=23$. We show that ${\rm A}_{p+1}$ arises for all $p>7$, again with too many dessins to classify for all except small $p$, whereas ${\rm PSL}_2(p)$ arises only for ${\mathcal D}_0(p)$. The Mathieu groups ${\rm M}_{12}$ and ${\rm M}_{24}$ each yield two chiral pairs of dessins, whereas as before ${\rm M}_{11}$ yields none. Finally we show that the affine group ${\rm AGL}_n(2)$ yields at least one dessin for each Mersenne prime $p=2^n-1>7$; however, this fails for $p=7$ since ${\rm AGL}_3(2)$ is not a Hurwitz group.

In Section~\ref{sec:3dim} we briefly outline another realisation of ${\rm PSL}_2(11)$, this time as the isometry group of the hendecachoron, a tessellation of a non-orientable $3$-orbifold by $11$ hemi-icosahedra, discovered independently by Gr\"unbaum~\cite{Gru} and Coxeter~\cite{Cox84}.

Section \ref{sec:App} contains monodromy permutations and diagrams for some 
of the dessins which appear in this paper.


\section{Klein's ten plane trees}
\label{sec:klein-trees}

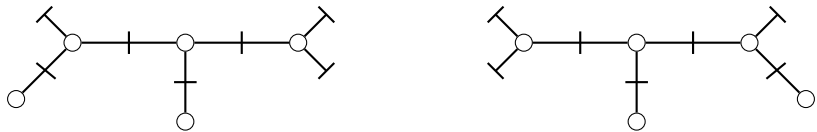
\begin{figure}[!ht]

\begin{center}
 \begin{tikzpicture}[scale=0.15, inner sep=0.8mm]
 
 \node (A) at (-30,0) [shape=circle, draw] {};
 \node (B) at (-20,0) [shape=circle, draw] {};
 \node (C) at (-10,0) [shape=circle, draw] {};
 \node (D) at (-20,-7) [shape=circle, draw] {};
 \node (E) at (-35,-5) [shape=circle, draw] {};

 \draw [thick] (A) to (B) to (C); 
 \draw [thick] (E) to (A) to (-32.5,2.5);
 \draw [thick] (-7.5,2.5) to (C) to (-7.5,-2.5);
 \draw [thick] (B) to (D);
 
 \draw [thick] (-33.2,-1.8) to (-31.5,-3.2);
 \draw [thick] (-25,1) to (-25,-1);
 \draw [thick] (-15,1) to (-15,-1);
 \draw [thick] (-21,-3.5) to (-19,-3.5);
 \draw [thick] (-8.2,3.2) to  (-6.8,1.8);
 \draw [thick] (-8.2,-3.2) to  (-6.8,-1.8);
 \draw [thick] (-31.8,3.2) to  (-33.2,1.8);
 
 
  \node (a) at (30,0) [shape=circle, draw] {};
 \node (b) at (20,0) [shape=circle, draw] {};
 \node (c) at (10,0) [shape=circle, draw] {};
 \node (d) at (20,-7) [shape=circle, draw] {};
 \node (e) at (35,-5) [shape=circle, draw] {};

 \draw [thick] (a) to (b) to (c); 
 \draw [thick] (e) to (a) to (32.5,2.5);
 \draw [thick] (7.5,2.5) to (c) to (7.5,-2.5);
 \draw [thick] (b) to (d);
 
 \draw [thick] (33.2,-1.8) to (31.5,-3.2);
 \draw [thick] (25,1) to (25,-1);
 \draw [thick] (15,1) to (15,-1);
 \draw [thick] (21,-3.5) to (19,-3.5);
 \draw [thick] (8.2,3.2) to  (6.8,1.8);
 \draw [thick] (8.2,-3.2) to  (6.8,-1.8);
 \draw [thick] (31.8,3.2) to  (33.2,1.8);
 
 \end{tikzpicture}

\end{center}
\caption{The chiral pair I of Klein's plane trees.} 
\label{Klein'sdessins}
\end{figure}

In~\cite{Kle79} Klein drew a diagram showing ten plane trees, each with eleven edges and twelve vertices. In each tree there is a bipartite partition of the vertices, with five vertices (coloured white) of valencies dividing $3$, and seven vertices (indicated by short cross-bars, perpendicular to their incident edges) of valencies dividing $2$; in each case the unique face is an $11$-gon. Among these plane trees there are four chiral (mirror-image) pairs, numbered I to IV, and two others, numbered V and VI, which exhibit bilateral symmetry. The chiral pair~I, which play a major role in Klein's paper, are shown (slightly distorted for simplicity, but combinatorially correct) in Figure~\ref{Klein'sdessins}.

\subsection{Klein's trees as dessins d'enfants}

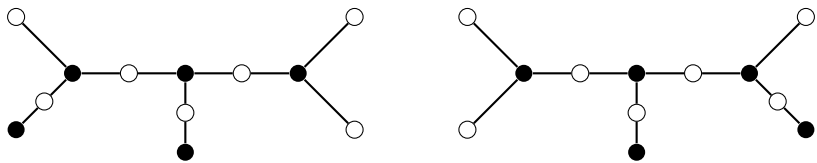
\begin{figure}[!ht]

\begin{center}
 \begin{tikzpicture}[scale=0.15, inner sep=0.8mm]
 
 \node (A) at (-30,0) [shape=circle, fill=black] {};
 \node (B) at (-20,0) [shape=circle, fill=black] {};
 \node (C) at (-10,0) [shape=circle, fill=black] {};
 \node (D) at (-20,-7) [shape=circle, fill=black] {};
 \node (E) at (-35,-5) [shape=circle, fill=black] {};
\node (F) at (-35,5) [shape=circle, draw] {};
\node (G) at (-5,5) [shape=circle, draw] {};
\node (H) at (-5,-5) [shape=circle, draw] {};
\node (I) at (-25,0) [shape=circle, draw] {};
\node (J) at (-15,0) [shape=circle, draw] {};
 \node (K) at (-32.5,-2.5) [shape=circle, draw] {};
 \node (L) at (-20,-3.5) [shape=circle, draw] {};

 \draw [thick] (A) to (I) to (B) to (J) to (C); 
 \draw [thick] (E) to (K) to (A) to (F);
 \draw [thick] (G) to (C) to (H);
 \draw [thick] (B) to (L) to (D);
 
 
 \node (a) at (30,0) [shape=circle, fill=black] {};
 \node (b) at (20,0) [shape=circle, fill=black] {};
 \node (c) at (10,0) [shape=circle, fill=black] {};
 \node (d) at (20,-7) [shape=circle, fill=black] {};
 \node (e) at (35,-5) [shape=circle, fill=black] {};
 \node (f) at (35,5) [shape=circle, draw] {};
 \node (g) at (5,5) [shape=circle, draw] {};
\node (h) at (5,-5) [shape=circle, draw] {};
\node (i) at (25,0) [shape=circle, draw] {};
\node (j) at (15,0) [shape=circle, draw] {};
 \node (k) at (32.5,-2.5) [shape=circle, draw] {};
 \node (l) at (20,-3.5) [shape=circle, draw] {};

 \draw [thick] (a) to (i) to (b) to (j) to (c); 
 \draw [thick] (e) to (k) to (a) to (f);
 \draw [thick] (g) to (c) to (h);
 \draw [thick] (b) to (l) to (d);

\end{tikzpicture}

\end{center}
\caption{The dessins ${\mathcal M}_1$ and $\overline{\mathcal M}_1$ corresponding to Klein's plane trees I.} 
\label{L2(11)dessins}
\end{figure}

If we regard the underlying surface as the complex plane~${\mathbb C}$, and compactify it to give the Riemann sphere or complex projective line $\Sigma={\mathbb P}^1({\mathbb C})={\mathbb C}\cup\{\infty\}$, we can recognise these plane trees as dessins d'enfants in Grothendieck's sense~\cite{Gro}, namely bipartite graphs embedded in compact Riemann surfaces $X$; as such they represent finite covers $X\to\Sigma$ of the sphere, unbranched outside $\{0, 1, \infty\}$, or equivalently, by Bely\u\i's Theorem~\cite{Bel}, projective algebraic curves $X$ defined over algebraic number fields. (See~\cite{GG,JW,LZ} for background on dessins d'enfants.) All ten of these dessins have genus $0$, so that $X=\Sigma$. They have degree~$11$, meaning that they represent $11$-sheeted coverings of $\Sigma$; they have type~$(3,2,11)$, meaning that there are branch-points of orders dividing $3, 2$ and $11$ over $0, 1$ and $\infty$, represented in Klein's diagrams by the white vertices, cross-bars and the point at infinity. A modern convention is to use black and white vertices for points over $0$ and $1$, with the edges corresponding to the points over the unit interval $[0,1]\subset\Sigma$, and to regard those over $\infty$ as the face-centres. We will denote the dessins corresponding to Klein's chiral pairs I to IV by ${\mathcal M}_i$ and $\overline{\mathcal M}_i$ for $i=1,\ldots, 4$, with ${\mathcal M}_i$ on the left in Klein's diagrams and ours; similarly ${\mathcal M}_i\;(i=5,6)$ is the dessin corresponding to the tree V or VI in~\cite{Kle79}\footnote{The temptation to follow Klein by using Roman numerals as subscripts here was almost (but not totally) irresistible.}. The dessins ${\mathcal M}_1$ and $\overline{\mathcal M}_1$ corresponding to Klein's plane trees I are shown in Figure~\ref{L2(11)dessins}.

\begin{figure}[!ht]

\begin{center}
 \begin{tikzpicture}[scale=0.14, inner sep=0.8mm]
 
 \node (A1) at (-30,30) [shape=circle, fill=black] {};
 \node (B1) at (-20,30) [shape=circle, fill=black] {};
 \node (C1) at (-10,30) [shape=circle, fill=black] {};
 \node (D1) at (-20,23) [shape=circle, fill=black] {};
 \node (E1) at (-35,25) [shape=circle, fill=black] {};

 \draw [thick] (A1) to (B1) to (C1); 
 \draw [thick] (E1) to (A1) to (-35,35);
 \draw [thick] (-5,25) to (C1) to (-5,35);
 \draw [thick] (B1) to (D1);
 
 \node at (-40,30) {${\mathcal M}_1$};
 
 
 \node (a1) at (30,30) [shape=circle, fill=black] {};
 \node (b1) at (20,30) [shape=circle, fill=black] {};
 \node (c1) at (10,30) [shape=circle, fill=black] {};
 \node (d1) at (20,23) [shape=circle, fill=black] {};
 \node (e1) at (35,25) [shape=circle, fill=black] {};

 \draw [thick] (a1) to (b1) to (c1); 
 \draw [thick] (e1) to (a1) to (35,35);
 \draw [thick] (5,25) to (c1) to (5,35);
 \draw [thick] (b1) to (d1);
 
  \node at (40,30) {$\overline{\mathcal M}_1$};
 
 
 \node (A2) at (-30,15) [shape=circle, fill=black] {};
 \node (B2) at (-20,15) [shape=circle, fill=black] {};
 \node (C2) at (-10,15) [shape=circle, fill=black] {};
 \node (D2) at (-20,8) [shape=circle, fill=black] {};
 \node (E2) at (-35,20) [shape=circle, fill=black] {};

 \draw [thick] (A2) to (B2) to (C2); 
 \draw [thick] (E2) to (A2) to (-35,10);
 \draw [thick] (-5,10) to (C2) to (-5,20);
 \draw [thick] (B2) to (D2);
 
 \node at (-40,15) {${\mathcal M}_2$};
 
 
 \node (a2) at (30,15) [shape=circle, fill=black] {};
 \node (b2) at (20,15) [shape=circle, fill=black] {};
 \node (c2) at (10,15) [shape=circle, fill=black] {};
 \node (d2) at (20,8) [shape=circle, fill=black] {};
 \node (e2) at (35,20) [shape=circle, fill=black] {};

 \draw [thick] (a2) to (b2) to (c2); 
 \draw [thick] (e2) to (a2) to (35,10);
 \draw [thick] (5,10) to (c2) to (5,20);
 \draw [thick] (b2) to (d2);
 
 \node at (40,15) {$\overline{\mathcal M}_2$};
 
 
\node (A3) at (10,0) [shape=circle, fill=black] {};
 \node (B3) at (20,0) [shape=circle, fill=black] {};
 \node (C3) at (30,0) [shape=circle, fill=black] {};
 \node (D3) at (5,5) [shape=circle, fill=black] {};
 \node (E3) at (5,-5) [shape=circle, fill=black] {};

 \draw [thick] (A3) to (B3) to (C3); 
 \draw [thick] (E3) to (A3) to (5,5);
 \draw [thick] (35,-5) to (C3) to (35,5);
 \draw [thick] (B3) to (20,-7);
 
  \node at (-40,0) {${\mathcal M}_3$};
 

 \node (a3) at (-10,0) [shape=circle, fill=black] {};
 \node (b3) at (-20,0) [shape=circle, fill=black] {};
 \node (c3) at (-30,0) [shape=circle, fill=black] {};
 \node (d3) at (-5,5) [shape=circle, fill=black] {};
 \node (e3) at (-5,-5) [shape=circle, fill=black] {};

 \draw [thick] (a3) to (b3) to (c3); 
 \draw [thick] (e3) to (a3) to (-5,5);
 \draw [thick] (-35,-5) to (c3) to (-35,5);
 \draw [thick] (b3) to (-20,-7);
 
 \node at (40,0) {$\overline{\mathcal M}_3$};
 

  \node (A4) at (-30,-15) [shape=circle, fill=black] {};
 \node (B4) at (-20,-15) [shape=circle, fill=black] {};
 \node (C4) at (-10,-15) [shape=circle, fill=black] {};
 \node (D4) at (-5,-10) [shape=circle, fill=black] {};
 \node (E4) at (-35,-20) [shape=circle, fill=black] {};

 \draw [thick] (A4) to (B4) to (C4); 
 \draw [thick] (E4) to (A4) to (-35,-10);
 \draw [thick] (-5,-20) to (C4) to (-5,-10);
 \draw [thick] (B4) to (-20,-22);
 
 \node at (-40,-15) {${\mathcal M}_4$};
 
 
 \node (a4) at (30,-15) [shape=circle, fill=black] {};
 \node (b4) at (20,-15) [shape=circle, fill=black] {};
 \node (c4) at (10,-15) [shape=circle, fill=black] {};
 \node (d4) at (5,-10) [shape=circle, fill=black] {};
 \node (e4) at (35,-20) [shape=circle, fill=black] {};

 \draw [thick] (a4) to (b4) to (c4); 
 \draw [thick] (e4) to (a4) to (35,-10);
 \draw [thick] (5,-20) to (c4) to (5,-10);
 \draw [thick] (b4) to (20,-22);
  
 \node at (40,-15) {$\overline{\mathcal M}_4$};
 

 \node (A5) at (-10,-30) [shape=circle, fill=black] {};
 \node (B5) at (0,-30) [shape=circle, fill=black] {};
 \node (C5) at (10,-30) [shape=circle, fill=black] {};
 \node (D5) at (-15,-35) [shape=circle, fill=black] {};
 \node (E5) at (15,-35) [shape=circle, fill=black] {};

 \draw [thick] (A5) to (B5) to (C5); 
 \draw [thick] (D5) to (A5) to (-15,-25);
 \draw [thick] (15,-35) to (C5) to (15,-25);
 \draw [thick] (B5) to (0,-37);
 
  \node at (-20,-30) {${\mathcal M}_5$};
 
 
 \node (a6) at (10,-45) [shape=circle, fill=black] {};
 \node (b6) at (0,-45) [shape=circle, fill=black] {};
 \node (c6) at (-10,-45) [shape=circle, fill=black] {};
  \node (d6) at (15,-40) [shape=circle, fill=black] {};
 \node (e6) at (-15,-40) [shape=circle, fill=black] {};

 \draw [thick] (a6) to (b6) to (c6); 
 \draw [thick] (d6) to (a6) to (15,-50);
 \draw [thick] (-15,-50) to (c6) to (-15,-40);
 \draw [thick] (b6) to (0,-52);
  
 \node at (-20,-45) {${\mathcal M}_6$};

\end{tikzpicture}

\end{center}
\caption{The maps ${\mathcal M}_i$ and $\overline{\mathcal M}_i$ corresponding to Klein's plane trees I -- VI.} 
\label{Kleinmaps}
\end{figure}
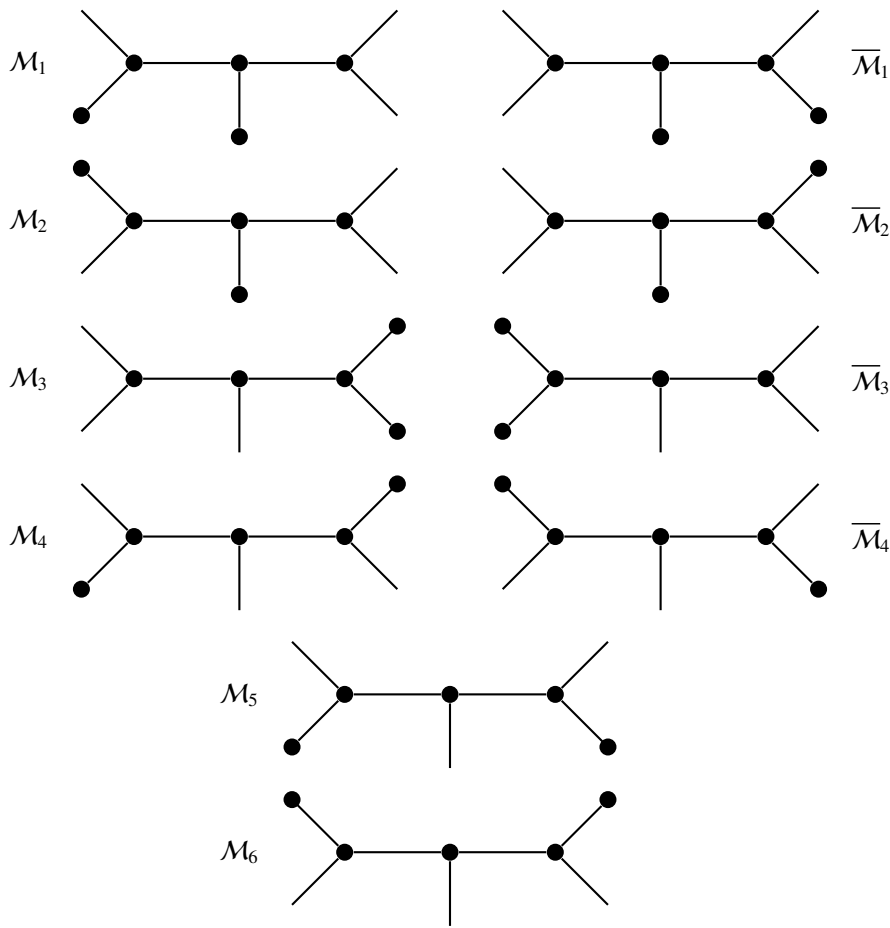

When the order of branching over $1$ divides $2$, as in these dessins, so that white vertices all have valencies $1$ or $2$, it is sometimes convenient to omit these vertices, leaving a map with only black vertices and possibly some free edges where white vertices of valency $1$ have been removed. Of course, this operation is reversible, so no information is lost.
With this convention, the dessins corresponding to all ten of Klein's plane trees are shown in Figure~\ref{Kleinmaps}, with the same layout as in his original diagram. 

In fact these ten dessins form a complete list (up to isomorphism) of the dessins of type $(3,2,11)$ and degree $11$, which is the least possible degree for dessins of this type. Klein presented his diagram as representing a classification of the $11$-sheeted coverings of the sphere $\Sigma$ with branching patterns $3^31^2$ and $2^41^3$ over $0$ and $1$. He explained this in terms of the modular function $J:{\mathbb H}\to\Sigma$, which parametrises isomorphism classes of elliptic curves, an approach which leads naturally to the modular group $\Gamma={\rm PSL}_2({\mathbb Z})$, the subgroup leaving $J$ invariant in the automorphism group ${\rm PSL}_2({\mathbb R})$ of the Riemann surface $\mathbb H$. This part of Klein's work will not be discussed in our paper.


\subsection{The classification problem}

Klein presented his ten planar trees in~\cite[\S1]{Kle79} in answer to a question he had posed in~\cite[\S5]{Kle}: in modern terminology, he had asked for the number of dessins $\mathcal D$ of degree $11$ and passport $(3^31^2; 2^41^3; -)$, that is, with three black vertices of valency $3$ and two of valency $1$, together with four white vertices of valency $2$ and three of valency $1$, and with no restrictions on the faces.
 
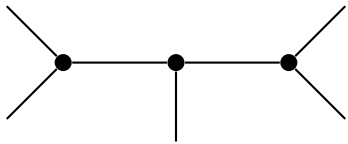
\begin{figure}[!ht]

\begin{center}
 \begin{tikzpicture}[scale=0.15, inner sep=0.8mm]
 
 \node (A) at (-10,0) [shape=circle, fill=black] {};
 \node (B) at (0,0) [shape=circle, fill=black] {};
 \node (C) at (10,0) [shape=circle, fill=black] {};

 \draw [thick] (A) to (B) to (C); 
 \draw [thick] (-15,5) to (A) to (-15,-5);
 \draw [thick] (15,5) to (C) to (15,-5);
 \draw [thick] (B) to (0,-7);

\end{tikzpicture}

\end{center}
\caption{The planar map $\mathcal M_0$.}
\label{9dartmap}
\end{figure}

In his diagram Klein drew ten planar trees, equivalent to these dessins, but did not explain why this list was complete, writing simply `Dass es auch nicht mehr giebt, ist ebenso evident'\footnote{`It is also clear that there are no others.'}. The argument, which we will now give in the language of dessins, is indeed straightforward. Since each dessin $\mathcal D$ has eleven edges and twelve vertices the embedded graph, being connected, must be a tree, so the underlying surface is the sphere and there is a single face, of valency $11$. Since the white vertices have valency $1$ or $2$ we can omit them and for simplicity represent $\mathcal D$ as a map $\mathcal M$, which must have three vertices of valency $3$ and two of valency $1$, and seven edges, three of them free. By connectedness, one of the three vertices of valency $3$ must be adjacent to the other two, so $\mathcal M$ must be formed from the basic map $\mathcal M_0$ in Figure~\ref{9dartmap} by choosing two of its five free edges and adding a vertex to each. This gives us 
$\displaystyle \left(\!\!\begin{array}{c} 5 \\ 2 \end{array}\!\!\right)=10$ 
non-isomorphic maps $\mathcal M$; they are shown in Figure~\ref{Kleinmaps}, with the same layout as the corresponding plane trees in Klein's paper~\cite{Kle79}.
 
An alternative approach to this and similar problems, based on the character theory and
the Frobenius formula, will be developed in Section~\ref{triples}.


\section{Monodromy groups}\label{sec:monodromy}

Any dessin $\mathcal D$ of degree $n$ can be represented as an ordered pair $x,y$ of permutations of its $n$ edges, obtained by using the chosen orientation of the underlying surface to follow the rotation of edges around their incident black and white vertices. The monodromy group $G=\langle x, y\rangle$ of $\mathcal D$, through its action permuting the $n$ edges, can be identified with a subgroup of the symmetric group ${\rm S}_n$. By the connectedness of the embedded graph it is transitive. It is convenient
to introduce a third permutation $z:=(xy)^{-1}$, so that $xyz=1$. It is easy to verify
that, while the cycles of $x$ and $y$ correspond respectively to black and white vertices,
the cycles of $z$ correspond to faces.

\begin{figure}[htbp]
\begin{center}

\includegraphics[scale=0.3]{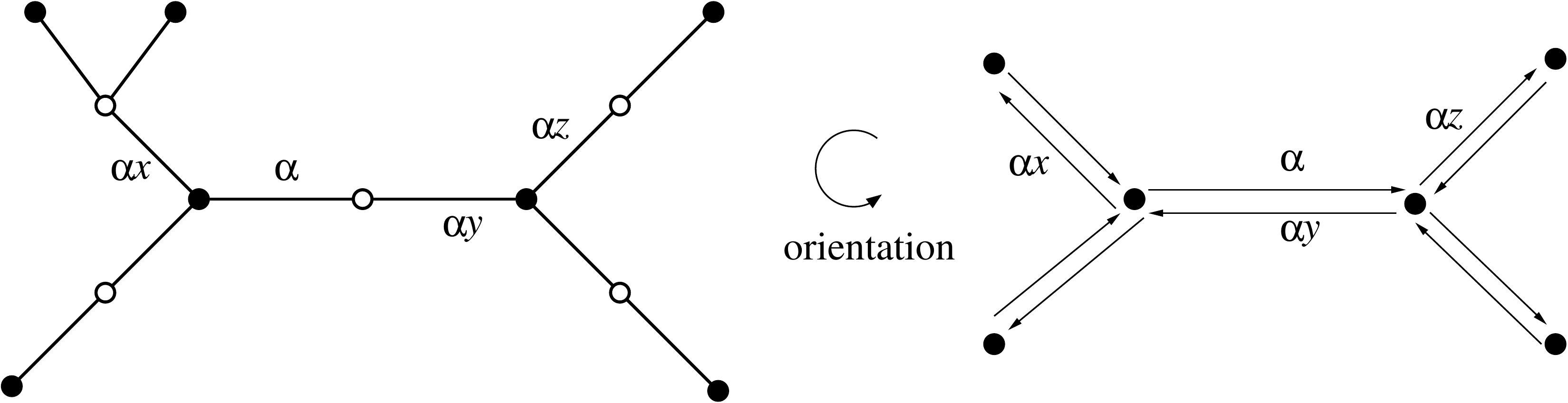}
\end{center}
\caption{Monodromy permutations: here $\alpha x$, $\alpha y$ and $\alpha z$ are the 
images of $\alpha$ under the permutations $x$, $y$ and $z=(xy)^{-1}$.}
\label{monodromy}
\end{figure}

\begin{conv}[Where to put labels]\label{conv:labels}\rm
We put a label of an edge on its left side while moving from its black end to 
its white one. In this way the labels corresponding to a cycle of $z$ or, equivalently,
to a face, will be situated inside this face and will be rotated by $z$ around the center
of the face in the direction corresponding to the orientation of the surface.
This convention is illustrated in Figure~\ref{monodromy}.

(When a dessin is drawn on the plane, its outer face gives an impression that its
labels are rotated in the direction opposite to the orientation of the plane. 
But this is an illusion. In fact, the dessin should be considered not on the plane 
but on a sphere, and the center of the outer face is situated on the ``opposite side'' 
of this sphere. Looking from this center, the rotation of the labels of the outer face
does correspond to the orientation of the sphere.)

When $y^2=1$ and the white vertices are omitted, so that $\mathcal D$ is a map, one 
can regard $G$ as permuting its half-edges, or else directed edges. Thus, an edge 
of a map bears two labels while a free edge bears only one. A free edge may be considered
as an outgoing directed edge.
\end{conv}

\begin{figure}[htbp]
\begin{center}

\includegraphics[scale=0.3]{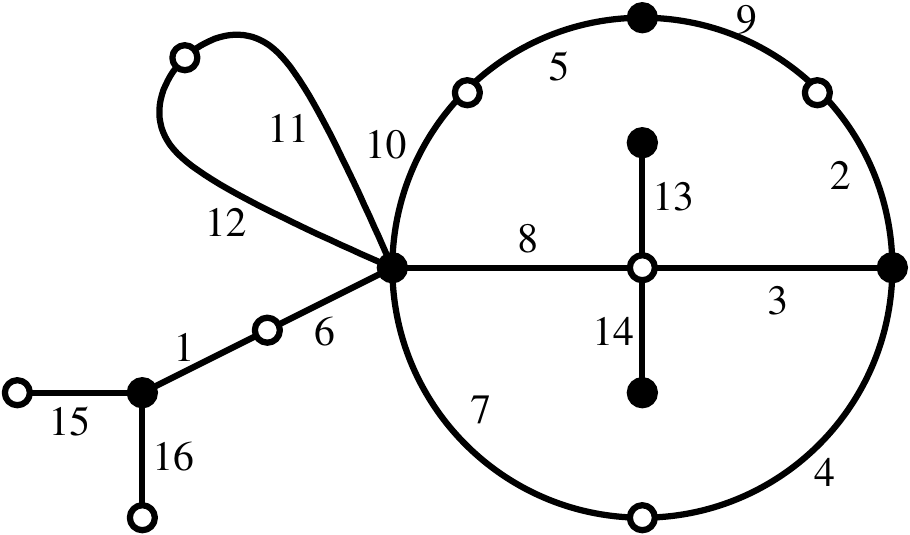}
\end{center}
\caption{A map with labelled edges.}
\label{fig:labels}
\end{figure}

\begin{exam}\rm
The permutations corresponding to the map in Figure \ref{fig:labels} are as follows:
\begin{eqnarray*}
x & = & (1,15,16)(2,3,4)(5,9)(6,7,8,10,11,12), \\
y & = & (1,6)(2,9)(3,13,8,14)(4,7)(5,10)(11,12), \\
z & = & (1,12,10,9,4,6,16,15)(2,5,8,13)(3,14,7).
\end{eqnarray*}
The monodromy group of this map is the alternating group ${\rm A}_{16}$. See also 
Figures \ref{deg12map} and \ref{fig:psl-5-2} for other examples of maps with labelled edges.
\end{exam}

The ten dessins $\mathcal D$ considered here have degree $11$, so their monodromy groups $G$ are transitive subgroups of ${\rm S}_{11}$. In each case $x^3=y^2=1$, so $G$ is a quotient of the modular group
\[\Gamma={\rm PSL}_2({\mathbb Z})=\langle X, Y, Z\mid X^3=Y^2=XYZ=1\rangle\cong C_3*C_2,\]
where the generators
\[X=\pm\left(\,\begin{array}{cc}
0  & 1 \\
-1 & -1 
\end{array}
\right),
\quad
Y=\pm\left(\,\begin{array}{cc}
0  & 1 \\
-1 & 0 
\end{array}
\right),
\quad
Z=\pm\left(\,\begin{array}{cc}
1 & 0 \\
1 & 1 \\
\end{array}
\right)
\]
correspond to the M\"obius transformations
\[t\mapsto\frac{-1}{t-1},\quad t\mapsto \frac{-1}{t}\quad\hbox{and}\quad t\mapsto t+1.\]

(Here we have associated M\"obius transformations $\displaystyle t\mapsto\frac{at+b}{ct+d}$ with pairs of matrices 
$\displaystyle \pm\,\left( \begin{array}{cc} a & c \\ b & d \end{array} \right)$, 
the transposes of those normally used in this context. This is because we will later need to regard M\"obius transformations, especially those over finite fields, as permutations; it is algebraically and computationally convenient to compose these from left to right, instead of using the analytic right to left convention, so in order to have a homomorphism from matrices to M\"obius transformations we will ignore the deceptively suggestive notation for the latter and use transposes.) 

Klein was interested in which of these monodromy groups $G$ are congruence quotients, that is, which are quotients by congruence subgroups of $\Gamma$. The only possible level for such a congruence subgroup is $11$ (since this is the order of $z:=(xy)^{-1}$), and indeed the resulting quotient group $L:={\rm PSL}_2(11)$, a simple group of order $660$, does act as a transitive group of degree $11$, namely on the cosets of a subgroup $H\cong {\rm A}_5$, as we will now show. This is one of the three cases, all known to Galois (see \cite{last-letter}, pp.~411--412), in which a simple group ${\rm PSL}_2(p)$, where $p$ is prime, has a proper subgroup of index less than $p+1$; the others are ${\rm A}_4<{\rm PSL}_2(5)$ and ${\rm S}_4<{\rm PSL}_2(7)$.

As generators of $L$ one can take the images
\[x=\pm\left(\,\begin{array}{cc}
0 & 1 \\
-1 & -1 \\
\end{array}
\right),
\quad
y=\pm\left(\,\begin{array}{cc}
0 & 1 \\
-1 & 0 \\
\end{array}
\right),
\quad
z=\pm\left(\,\begin{array}{cc}
1 & 0 \\
1 & 1 \\
\end{array}
\right)
\]
of $X, Y$ and $Z$, satisfying the relations
\begin{equation}\label{2311eqns}
x^3=y^2=z^{11}=xyz=1.
\end{equation}
The elements
\[u=\pm\left(\,\begin{array}{cc}
4 & -1 \\
2 & -3 \\
\end{array}
\right),
\quad
v=y=\pm\left(\,\begin{array}{cc}
0 & 1 \\
-1 & 0 \\
\end{array}
\right),
\quad
w=\pm\left(\,\begin{array}{cc}
2 & -4 \\
-3 & 1 \\
\end{array}
\right)
\]
of $L$ (note that $\det(u)=\det(w)=-10\equiv 1$ mod~$(11)$) satisfy
\[u^3=v^2=w^5=uvw=1,\]
so they generate a quotient $H\le L$ of the triangle group $\Delta(3,2,5)\cong {\rm A}_5$. By the simplicity of ${\rm A}_5$ we have $H\cong {\rm A}_5$, so $|L:H|=11$. Conjugation by the element
\[\pm\left(\,\begin{array}{cc}
-1 & 0 \\
0 & 1 \\
\end{array}\right)\in {\rm PGL}_2(11)\cong{\rm Aut}\,L\]
induces an outer automorphism of $L$, transposing $u, v$ and $w$ with elements
\[\overline u=\pm\left(\,\begin{array}{cc}
4 & 1 \\
-2 & -3 \\
\end{array}
\right),
\quad
\overline v=v=y=\pm\left(\,\begin{array}{cc}
0 & -1 \\
1 & 0 \\
\end{array}
\right),
\quad
\overline w=\pm\left(\,\begin{array}{cc}
2 & 4 \\
3 & 1 \\
\end{array}
\right)
\]
which satisfy the same relations and therefore generate a subgroup $\overline H\cong {\rm A}_5$. This subgroup is not conjugate in $L$ to $H$: if it were, some non-identity element of 
${\rm PGL}_2(11)$ would commute with $u, v$ and $w$, whereas it is straightforward to check that only the identity element does this. It follows that $L$ has two inequivalent transitive representations of degree $11$, on the cosets of the subgroups $H$ and $\overline H$, or equivalently (since these subgroups are equal to their normalisers in $L$), by conjugation on the two conjugacy classes of subgroups they represent.

Nowadays, by using a program such as GAP~\cite{GAP} one can easily determine the monodromy groups of the ten dessins, showing that those with monodromy group $L={\rm PSL}_2(11)$ are the chiral pair $\mathcal M_1$ and $\overline{\mathcal M}_1$. However, one can also distinguish this pair by hand: for $i=1,\ldots, 6$ the commutator $[x,y]=x^{-1}y^{-1}xy$ has cycle structure $5^21^1, 7^13^11^1, 8^12^11^1, 5^14^12^1, 7^12^2$ and $4^23^1$ in ${\mathcal M}_i$ (and hence also in $\overline{\mathcal M}_i$), so only when $i=1$ does it have one of the possible orders $1, 2, 3, 5, 6$ or $11$ of an element of $L$. Conversely, it is easy to see that when $L$ is represented on the cosets of $H$ or $\overline H$ the generators $x$ and $y$ have cycle structures $3^31^2$ and $2^41^3$, so $L$ is the monodromy group of at least one of the ten dessins, and this must be the pair $\mathcal M_1$ and $\overline{\mathcal M}_1$, with the two dessins corresponding to the two inequivalent representations of degree $11$.

Klein~\cite[\S2]{Kle79} used a similar argument, based on a word equivalent to our $yz^3$, to show that the plane trees labelled II--VI in his diagram do not have monodromy group $L$; he gave details only for tree~V, corresponding to the dessin ${\mathcal M}_5$. His paper is otherwise entirely devoted to the two trees labelled~I (our Figure~\ref{Klein'sdessins}) with monodromy group $L$, together with their minimal regular cover, which correspond to congruence subgroups of the modular group; the other eight trees, labelled~II -- VI and corresponding to noncongruence subgroups, are ignored, as are their monodromy groups.

In fact, the cycle structures of $[x,y]$ given above show that for $i=2, 4,5$ and $6$ the monodromy group $G$ is ${\rm A}_{11}$. Being transitive of prime degree, $G$ must be primitive; now a theorem of Jordan (see~\cite[Theorem~13.9]{Wie}, for example) shows that any finite primitive permutation group containing a cycle of prime length with at least three fixed points must contain the alternating group; in each of these four cases a suitable power of $[x,y]$ has this property, so $G\ge {\rm A}_{11}$. Since the generators $x$ and $y$ are even permutations, it follows that $G={\rm A}_{11}$. However, a separate argument is required for $i=3$, for instance using the fact that in this case $xz^3xz^5$ has cycle structure $5^13^11^3$. (The theorem of Jordan used here follows immediately from results in his papers~\cite{Jor71, Jor73} published in 1871 and 1873, so Klein could have known and used this argument; for further historical comments on Jordan's Theorem, and for a modern extension of it, see~\cite{Jon14}.)

The fact that these eight dessins have monodromy groups $G\cong {\rm A}_{11}$ shows that, as claimed above, they correspond to noncongruence subgroups of $\Gamma$: the only nonabelian composition factors of any congruence quotient of $\Gamma$ are isomorphic to ${\rm PSL}_2(p)$ for primes $p\ge 5$ dividing the level~\cite{McQ}, and by comparing orders it is easily seen that the simple group ${\rm A}_{11}$ does not have this form.


\section{Counting triples}\label{triples}

We have already given a very simple argument showing why there are exactly ten dessins with the given characteristics. However, this simple example allows us an opportunity to illustrate another important enumeration technique based on character theory. Quite often, in more complicated cases where the degree or the genus is greater, this method is the only one available. It is based on Frobenius's formula~\cite{Fro}
\begin{eqnarray}\label{Frobformula}
\frac{|{\mathcal X}|\cdot|{\mathcal Y}|\cdot|{\mathcal Z}|}{|G|}\sum_{\chi}
\frac{\chi(x)\chi(y)\chi(z)}{\chi(1)}
\end{eqnarray}
for the number of triples $(x,y,z)\in\mathcal X\times\mathcal Y\times\mathcal Z$ with $xyz=1$, where $\mathcal X$, $\mathcal Y$ and $\mathcal Z$ are conjugacy classes in a finite group $G$, and the sum is over the irreducible complex characters $\chi$ of $G$. (In specific applications, one often finds that many characters take the value $0$ on $\mathcal X$, $\mathcal Y$ or $\mathcal Z$, so they can be omitted from the summation.) Since~\cite{Fro} was published 17 years after~\cite{Kle79}, this powerful technique was not available to Klein. Here we will show how it works in this fairly simple situation, ignoring results explained earlier; as a bonus we obtain the monodromy groups of all ten dessins.

Before doing this, we will give an alternative form of~(\ref{Frobformula}) which is sometimes more convenient, for instance when using the ATLAS~\cite{ATLAS} where orders of centralisers, rather than conjugacy classes, are given. We have $|\mathcal X|=|G|/|C(x)|$ where $C(x)=C_G(x)$ denotes the centraliser of $x$ in $G$, with similar equations for $|\mathcal Y|$ and $|\mathcal Z|$, so we can rewrite ~(\ref{Frobformula}) as
\begin{eqnarray}\label{Frobformula2}
\frac{|G|^2}{|C(x)|\cdot|C(y)|\cdot|C(z)|}\sum_{\chi}
\frac{\chi(x)\chi(y)\chi(z)}{\chi(1)}.
\end{eqnarray}

If a dessin $\mathcal D$ has type $(3,2,11)$ and degree $11$, then its monodromy permutations $x, y, z\in {\rm S}_{11}$ satisfy the relations~(\ref{2311eqns}); having odd order, $x$ and $z$ are even, and hence so is $y$, so the monodromy group $G=\langle x, y, z\rangle$ is a subgroup of $A:={\rm A}_{11}$. Our aim is to count the triples $x$, $y$, $z$ in $A$ satisfying~(\ref{2311eqns}),
and to determine the corresponding dessins $\mathcal D$ and their monodromy groups $G$. In more general applications of~(\ref{Frobformula}) one has to exclude triples generating intransitive subgroups $G$, since these cannot be monodromy groups of connected covers. However, here the fact that $z$ must be an $11$-cycle guarantees transitivity, so this problem does not arise. Since the orders $3, 2$ and $11$ of $x, y$ and $z$ are mutually coprime, $G$ must be perfect (have trivial abelianisation), and therefore cannot be solvable. If $G\ne A$ then $G$ must be contained in a nonsolvable maximal subgroup $M$ of $A$, of order divisible by $11$. According to the list of maximal subgroups of ${\rm A}_{11}$ in~\cite{ATLAS}, the only such maximal subgroups $M$ are isomorphic to the Mathieu group ${\rm M}_{11}$, forming two conjugacy classes of subgroups, of index 2520, transposed by ${\rm Out}\,A$. Similarly the list of maximal subgroups of ${\rm M}_{11}$ in~\cite{ATLAS} tells us that if $G\ne M$ then $G$ is isomorphic to $L={\rm PSL}_2(11)$, each subgroup $M\cong {\rm M}_{11}$ containing one conjugacy class of such subgroups, of index $12$ in $M$. Thus $G$ is equal to $A$, or is isomorphic to ${\rm M}_{11}$ or ${\rm PSL}_2(11)$; we will count the triples generating these three groups in turn, starting with the smallest.

\begin{exam}\label{ex:PSL_2(11)}\rm
We will show that the group $L={\rm PSL}_2(11)$ has, up to automorphisms (equivalently, up to conjugacy in ${\rm PGL}_2(11)\cong {\rm Aut}\,L$), a single generating triple $x, y, z$ satisfying~(\ref{2311eqns}). The character table of $L$, produced by GAP~\cite{GAP} and given here, has columns corresponding to the conjugacy classes, with the headings $n$a, $n$b, etc. ($n$A, $n$B, etc in ATLAS notation~\cite{ATLAS}) indicating successive classes of elements of order $n$; the rows, labelled X.j, correspond to irreducible characters $\chi_j$ for $j=1, 2,\ldots$. The table includes entries A, $^*$A\,$:=(-1\mp\sqrt 5)/2$ and ${\rm B}, /{\rm B}:=(-1\pm i\sqrt{11})/2$, while a full stop denotes an entry $0$. 
\begin{quote}
\begin{verbatim}
       1a 3a 2a 5a 5b 6a 11a 11b

X.1     1  1  1  1  1  1   1   1
X.2     5 -1  1  .  .  1   B  /B
X.3     5 -1  1  .  .  1  /B   B
X.4    10  1 -2  .  .  1  -1  -1
X.5    10  1  2  .  . -1  -1  -1
X.6    11 -1 -1  1  1 -1   .   .
X.7    12  .  .  A *A  .   1   1
X.8    12  .  . *A  A  .   1   1

A = E(5)^2+E(5)^3
  = (-1-Sqrt(5))/2 = -1-b5
B = E(11)+E(11)^3+E(11)^4+E(11)^5+E(11)^9
  = (-1+Sqrt(-11))/2 = b11
\end{verbatim}
\end{quote}

\medskip

In this example we take $\mathcal X$ and $\mathcal Y$ to be the classes 3a and 
2a, and $\mathcal Z$ to be either of the classes 11a or 11b.
For each choice of $\mathcal Z$ the character sum in formula~(\ref{Frobformula2}) has the value
\[\Sigma_L=1+\frac{(-1)\cdot 1\cdot {\rm B}}{5}+\frac{(-1)\cdot 1\cdot /{\rm B}}{5}+\frac{1\cdot (-2)\cdot (-1)}{10}+\frac{1\cdot 2\cdot (-1)}{10}
=1-\frac{{\rm B}+/{\rm B}}{5}=\frac{6}{5}.\]
Now the centralisers $C_L(x)$, $C_L(y)$ and $C_L(z)$ of typical elements $x, y$ and $z$ of $\mathcal X$, $\mathcal Y$ and $\mathcal Z$  have orders $6$, $12$ and $11$,
so that
\[\frac{|L|^2}{|C(x)|\cdot|C(y)|\cdot|C(z)|}=\frac{2^4\cdot 3^2\cdot 5^2\cdot 11^2}{6\cdot 12 \cdot 11}=2\cdot 5^2\cdot11.\]
Multiplying this by $\Sigma_L$ gives $2^2\cdot 3\cdot 5\cdot 11$ triples $(x, y, z)\in{\mathcal X}\times{\mathcal Y}\times{\mathcal Z}$ with $xyz=1$ for each choice of $\mathcal Z$, and hence $2^3\cdot 3\cdot 5\cdot 11$ in total. Each triple must generate $L$ since this group has no proper subgroups of order divisible by $66$ (in fact, none of index less that $11$). Since ${\rm Aut}\,L$ acts semi-regularly on generating triples, dividing by its order $|{\rm PGL}_2(11)|=2^3\cdot 3\cdot 5\cdot 11$ shows that ${\rm Aut}\,L$ has a single orbit on such triples, as claimed.

Since the normal subgroups of any group $\Delta$ with quotient isomorphic to a group $G$ are in bijective correspondence with the orbits of ${\rm Aut}\,G$, acting by composition, on epimorphisms $\Delta\to G$, it follows that the triangle group
\[\Delta:=\Delta(3,2,11)=\langle X, Y, Z\mid X^3=Y^2=Z^{11}=XYZ=1\rangle\]
 has a single normal subgroup with quotient isomorphic to $L$. Now $L$ has two conjugacy classes of subgroups of index $11$, all isomorphic to ${\rm A}_5$, and these classes are transposed by its automorphism group ${\rm PGL}_2(11)$ (see~\cite{ATLAS} or \cite[\S259]{Dic}); they lift to two conjugacy classes of subgroups of index $11$ in $\Delta$, transposed by conjugation in the extended triangle group
\[\Delta[3,2,11]=\langle R_0, R_1, R_2\mid R_i^2=(R_1R_2)^3=(R_2R_0)^2=(R_0R_1)^{11}=1\rangle,\]
which contains $\Delta$ with index $2$, and they correspond to a chiral pair of dessins $\mathcal D$ of type $(3,2,11)$ degree $11$ with monodromy group isomorphic to $L$. This calculation shows that they are the only such dessins, so they are the dessins ${\mathcal M}_1$ and $\overline{\mathcal M}_1$ discussed earlier.

\end{exam}

\begin{exam}\label{ex:M_11}\rm
Next we count triples satisfying~(\ref{2311eqns}) in the simple group $M={\rm M}_{11}$ of order $7920=2^4\cdot 3^2\cdot 5\cdot 11$. The character table of $M$, given by GAP, is shown here:

\medskip

\begin{quote}
\begin{verbatim}
        1a 11a 11b 2a 4a 8a 8b 3a 6a 5a

X.1      1   1   1  1  1  1  1  1  1  1
X.2     10  -1  -1  2  2  .  .  1 -1  .
X.3     10  -1  -1 -2  .  B -B  1  1  .
X.4     10  -1  -1 -2  . -B  B  1  1  .
X.5     11   .   .  3 -1 -1 -1  2  .  1
X.6     16   A  /A  .  .  .  . -2  .  1
X.7     16  /A   A  .  .  .  . -2  .  1
X.8     44   .   .  4  .  .  . -1  1 -1
X.9     45   1   1 -3  1 -1 -1  .  .  .
X.10    55   .   . -1 -1  1  1  1 -1  .

A = E(11)^2+E(11)^6+E(11)^7+E(11)^8+E(11)^10
  = (-1-Sqrt(-11))/2 = -1-b11
B = -E(8)-E(8)^3
  = -Sqrt(-2) = -i2
\end{verbatim}
\end{quote}

\medskip

We must take $\mathcal X$ and $\mathcal Y$ to be the unique conjugacy classes 3a and 2a 
of elements of order $3$ and~$2$, and $\mathcal Z$ to be either of the classes 11a or 11b 
of elements of order $11$. We see that the character sum for $M$ is
\[\Sigma_M=1+\frac{1\cdot 2\cdot (-1)}{10}+\frac{1\cdot (-2)\cdot (-1)}{10}+\frac{1\cdot (-2)\cdot (-1)}{10}=\frac{6}{5}.\]
By~\cite{ATLAS} the centralisers in $M$ of $x, y$ and $z$ have orders $18$, $48$ and $11$, so
\[\frac{|M|^2}{|C(x)|\cdot|C(y)|\cdot|C(z)|}=\frac{2^8\cdot3^4\cdot 5^2\cdot 11^2}{18\cdot 48\cdot 11}=2^3\cdot 3\cdot 5^2\cdot 11.\]
Multiplying this by $2\Sigma_M$ gives $2^5\cdot 3^2\cdot 5\cdot 11$ triples in $M$. Now $M$ has a single conjugacy class of $12$ subgroups $L\cong {\rm PSL}_2(11)$, and we saw in 
Example~\ref{ex:PSL_2(11)} that each of these is generated by  $2^3\cdot 3\cdot 5\cdot 11$ triples, so $2^5\cdot 3^2\cdot 5\cdot 11$ triples in $M$ generate proper subgroups. Thus none of these triples can generate $M$. (See~\cite{Con91, Wol90} for confirmation that ${\rm M}_{11}$ is not a quotient of $\Delta$.)

The groups ${\rm M}_{11}$ and ${\rm M}_{12}$ were introduced by Mathieu~\cite{Mat61} in 1861, so they would have been known to Klein; unlike ${\rm M}_{12}$ (see~\cite[\S2.1.6]{JW}, for example), ${\rm M}_{11}$ is not the monodromy group of any orientable cubic map, since an extension of the argument used here shows that it is not a quotient of $\Gamma$.
\end{exam}

\begin{exam}\label{ex:A_11}\rm
Now we count triples satisfying~(\ref{2311eqns}) in $A={\rm A}_{11}$. In this case there are three conjugacy classes of elements of order $3$, and two classes of elements of order~$2$. A simple argument shows that if elements $x, y\in A$ of orders $3$ and $2$ generate a transitive group than they must have cycle structures  $3^31^2$ and $2^41^3$, so they belong to the conjugacy classes 
$\mathcal X={\rm 3c}$ and $\mathcal Y={\rm 2b}$ of $A$, while $z$ must be in one of the two classes 11a and 11b of $11$-cycles in $A$.

The character table of $A$ is much too large to give here. However, the only irreducible characters $\chi$ with $\chi(x)\chi(y)\chi(z)\ne 0$ are the principal character $\chi_1$ and the characters $\chi_2$, $\chi_6$ and $\chi_{11}$ of degrees $10$, $120$ and $210$, so we can give here the relevant portion of the character table, concerning these classes and characters:

\medskip

\begin{quote}
\begin{verbatim}
          1a  2b  3c 11a 11b

X.1        1   1   1   1   1
X.2       10   2   1  -1  -1
X.6      120  -8   3  -1  -1
X.11     210   2   3   1   1
\end{verbatim}
\end{quote}

\medskip

 The character sum
\[\Sigma_A:=\sum_{\chi}\frac{\chi(x)\chi(y)\chi(z)}{\chi(1)}\]
has the value
\[1+\frac{1\cdot 2\cdot(-1)}{10}+\frac{3\cdot(-8)\cdot (-1)}{120}+\frac{3\cdot 2\cdot 1}{210}=1-\frac{1}{5}+\frac{1}{5}+\frac{1}{35}=\frac{36}{35}=\frac{2^2\cdot 3^2}{5\cdot 7}.\]
Now $|A|=11!/2=2^7\cdot 3^4\cdot 5^2\cdot 7\cdot11$, and by~\cite{ATLAS} the centralisers of $x, y$ and $z$ in $A$ have orders  $162=2\cdot 3^4$, $1152=2^7\cdot 3^2$ and $11$, so
\[\frac{|A|^2}{|C(x)|\cdot|C(y)|\cdot|C(z)|}=2^6\cdot 3^2\cdot 5^4 \cdot 7^2\cdot 11. \]
Multiplying this by $2\Sigma_A$ shows that the number of triples in $A$ satisfying~(\ref{2311eqns}) is
 \[t:=2^9\cdot 3^4\cdot 5^3\cdot 7\cdot 11=10\times 11!.\]

From this total we must subtract the number of triples generating proper subgroups of $A$. As we have seen, the only such subgroups are those isomorphic to ${\rm PSL}_2(11)$. Each subgroup $L\cong {\rm PSL}_2(11)$ in $A$ is contained in a unique subgroup $M\cong {\rm M}_{11}$; there are two conjugacy classes of 2520 such subgroups $M$ in $A$, each containing one conjugacy class of $12$ subgroups $L$, so the number of such subgroups $L$ in $A$ is $2\cdot 2520\cdot 12=60480=2^6\cdot 3^3\cdot 5\cdot 7$. As shown in Example~\ref{ex:PSL_2(11)}, each subgroup $L$ is generated by $2^3\cdot 3\cdot 5\cdot 11$ triples satisfying~(\ref{2311eqns}), so the number of triples in $A$ generating subgroups $L\cong {\rm PSL}_2(11)$ is $2^9\cdot 3^4\cdot 5^2\cdot 7\cdot 11=t/5=2\times 11!$.

Since there are $10\times 11!$ triples in $A$, of which $2\times 11!$ generate proper subgroups, there are $8\times11!$ triples generating $A$, forming eight orbits under ${\rm Aut}\,A={\rm S}_{11}$; these correspond to eight normal subgroups of $\Delta$ with quotient $A$, and hence to eight regular dessins with (orientation-preserving) automorphism group $A$. Now $A$ has a single conjugacy class of subgroups of index $11$, namely the natural point-stabilisers isomorphic to 
${\rm A}_{10}$, so these lift to eight more conjugacy classes of subgroups of index $11$ in 
$\Delta$, in addition to the two resulting from Example~\ref{ex:PSL_2(11)}. These give eight more dessins $\mathcal D$ of type $(3,2,11)$ and degree $11$, all with monodromy group $G\cong {\rm A}_{11}$. Thus we have a total of ten dessins of the required degree and type. This confirms the earlier enumeration, and also confirms that for ${\mathcal M}_1$ and $\overline{\mathcal M}_1$ the monodromy group is 
${\rm PSL}_2(11)$, whereas for the other eight dessins it is ${\rm A}_{11}$. Generating triples for these dessins are given in Subsection \ref{sec:kl-11}.
\end{exam}

In each of the examples considered here, it was a fairly easy matter to evaluate the Frobenius formulae~(\ref{Frobformula}) or (\ref{Frobformula2}) by hand, but in other cases, if there are many non-zero summands, or if irrational character values arise, it may be necessary to appeal to GAP for this. For instance, the following GAP commands repeat the calculation in 
Example~\ref{ex:PSL_2(11)}.

\medskip

\begin{quote}
\begin{verbatim}

> G:=PSL(2,11);;
> T:=CharacterTable(G);;
> OrdersClassRepresentatives(T);
  [ 1, 3, 2, 5, 5, 6, 11, 11 ]
> ClassStructureCharTable(T,[2,3,7]);
  660
> ClassStructureCharTable(T,[2,3,8]);
  660

\end{verbatim}
\end{quote}


\noindent Here $2, 3, 7$ and $2, 3, 8$ indicate the positions of the chosen classes $\mathcal X$, $\mathcal Y$ and $\mathcal Z$ in the preceding list. The double semicolons after a command ask GAP to execute the command but not to output the result. 
As in Example~\ref{ex:PSL_2(11)} we obtain $1320$ triples.


\section{The dessins and their regular covers}\label{regcov}

Our enumeration of triples in $L$ shows that there is a single regular dessin $\mathcal R_1$ of type $(3,2,11)$ with automorphism group $L\cong {\rm PSL}_2(11)$; this is the minimal regular cover of each of the two dessins $\mathcal M_1$ and $\overline{\mathcal M}_1$ with monodromy group $L$, and these are quotients of $\mathcal R_1$ by non-conjugate subgroups of $L$ isomorphic to ${\rm A}_5$. By the Riemann--Hurwitz formula a regular dessin of this type with automorphism group $G$ has genus $1+\frac{5}{132}|G|$, in this case equal to $26$. By its uniqueness, or by an observation of Singerman~\cite{Sin74} concerning triples in ${\rm PSL}_2(q)$, $\mathcal R_1$ is regular as a map, corresponding to a subgroup of $\Delta$ which is normal in the extended triangle group
$\Delta[3,2,11]$. It is, in fact, the dual of the unique orientable regular map of genus $26$ and type $\{3,11\}$, denoted by R26.2 in Conder's catalogue of regular maps~\cite{Con}; its full automorphism group as a map, containing also the automorphisms reversing the orientation, is isomorphic to 
${\rm PGL}_2(11)$. For a drawing of R26.2 (combinatorially though not geometrically correct), with $60$ $11$-valent vertices, $330$ edges and $220$ triangular faces, see~\cite[Fig.~6 and Table~1]{ISS}, where it is shown as a $198$-gon with side identifications; in this paper, Ivrissimtzis, Singerman and Strudwick reinterpret the regular maps associated with Klein's papers~\cite{Kle78} and~\cite{Kle79} in terms of Farey fractions.
 
The eight dessins $\mathcal D$ with monodromy group ${\rm A}_{11}$ consist of the three chiral pairs $\mathcal M_i$ and $\overline{\mathcal M}_i$ with $i=2,3,4$, together with $\mathcal M_5$ and $\mathcal M_6$. The minimal regular covers of these dessins are eight regular dessins $\mathcal R$ of genus $756\,001$, and each dessin $\mathcal D$ is the quotient of $\mathcal R$ by a subgroup isomorphic to ${\rm A}_{10}$ in its automorphism group ${\rm A}_{11}$. These eight dessins $\mathcal R$ are mutually non-isomorphic, since ${\rm A}_{11}$ has a unique conjugacy class of such subgroups (the point-stabilisers in the natural representation), so that any isomorphism between the dessins $\mathcal R$ would induce an isomorphism between their quotients $\mathcal D$, which is visibly impossible (see Figure~\ref{Kleinmaps}).

The bilateral symmetry of the dessins $\mathcal M_5$ and $\mathcal M_6$ is easily explained: in each case there is a permutation in $S_{11}$ simultaneously  inverting the generating permutations $x$ and $y$ for the monodromy group $A_{11}$, and this induces an isomorphism of the dessin with its mirror image. (Of course, since $y^2=1$, inverting $y$ means the same as commuting with it.) Equivalently, each of the corresponding map subgroups $M$ in $\Delta(3,2,11)$ is contained with index $2$ in a subgroup of $\Delta[3,2,11]$ which induces this isomorphism. The same applies to the regular covers of $\mathcal M_5$ and $\mathcal M_6$, each corresponding to the core (intersection of conjugates) of $M$ in $\Delta(3,2,11)$; the other six regular covers form three chiral pairs, like their quotients.

\begin{figure}[!ht]

\begin{center}
 \begin{tikzpicture}[scale=0.25, inner sep=0.8mm]
 
 \node (A) at (-10,0) [shape=circle, fill=black] {};
 \node (B) at (0,0) [shape=circle, fill=black] {};
 \node (C) at (10,0) [shape=circle, fill=black] {};
 \node (D) at (-15,-5) [shape=circle, fill=black] {};
 \node (E) at (15,-5) [shape=circle, fill=black] {};

 \draw [thick] (A) to (B) to (C); 
 \draw [thick] (D) to (A) to (-15,5);
 \draw [thick] (15,-5) to (C) to (15,5);
 \draw [thick] (B) to (0,-7);
 
\node at (1,-2) {$11$};
\node at (2,1) {$6$};
\node at (-2,-1) {$1$};
\node at (8,-1) {$10$};
\node at (-8,1) {$5$};
\node at (12,-0.6) {$8$};
\node at (-10.5,-2) {$2$};
\node at (10.8,2) {$7$};
\node at (-12.5,1) {$4$};
\node at (13,-4.2) {$9$};
\node at (-14.5,-3) {$3$};

\end{tikzpicture}

\end{center}
\caption{The map $\mathcal M_5$ corresponding to Klein's plane tree V.} 
\label{mapV}
\end{figure}
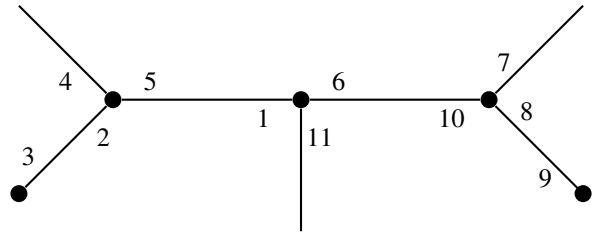

For example, Figure~\ref{mapV} shows the map $\mathcal M_5$ corresponding to $\mathcal D_5$; directed edges are labelled $1, 2, \ldots, 11$, so that
\[x=(1,11,6)(2,5,4)(7,10,8)\quad{\rm and}\quad y=(1,5)(2,3)(6,10)(8,9),\]
and hence
\[z=(1,2,\ldots,11).\]
The permutation $t=(1,6)(2,8)(3,9)(4,7)(5,10)$ inverts $x$ and $y$, inducing the automorphism of $\mathcal M_5$ given by  reflection in the vertical axis. The situation is similar for $\mathcal M_6$, with generators
\[x=(1,11,6)(2,5,3)(7,10,9)\quad{\rm and}\quad y=(1,5)(3,4)(6,10)(7,8)\]
inverted by $t=(1,6)(2,9)(3,7)(4,8)(5,10)$.

\section{Computations and Galois orbits}
\label{sec:CompGal}

Each of the ten dessins $\mathcal D$ discussed here represents a Bely\u\i\/ function, in this case an $11$-sheeted covering $\beta:\Sigma\to\Sigma$. As such, $\beta$ is a rational function of degree $11$, and having a single pole at $\infty$ it must be a polynomial, called a Shabat polynomial.

A computation using MAPLE shows that the Shabat polynomial $P:\Sigma\to\Sigma$ for the dessins $\mathcal M_1$ and $\overline{\mathcal M}_1$ in Figure~\ref{L2(11)dessins} has the form
\[P=\frac{1}{2^{12}3^{14}}p_1^3p_2^3p_3\]
where
\[p_1(x)=2x+(11-3\sqrt{-11}),\]
\[p_2(x)=2x^2-(11-3\sqrt{-11})x-(22+6\sqrt{-11}),\]
\[p_3(x)=x^2+11x+(55+9\sqrt{-11}),\]
so that
\[P-1=-\frac{1}{2^{11}3^{14}}q_1^2q_2^2q_3\]
where
\[q_1(x)=2x+(5+3\sqrt{-11}),\]
\[q_2(x)=2x^3+(15-3\sqrt{-11})x^2-(12-12\sqrt{-11})x+(56+96\sqrt{-11}),\]
\[q_3(x)=2x^3-18x^2+(21+45\sqrt{-11})x-(175+279)\sqrt{-11}).\]
The two possible choices for $\sqrt{-11}$ give the two dessins. They are defined over ${\mathbb Q}(\sqrt{-11})$ (which is not  surprising, since $\sqrt{-11}$ appears in the character table of 
${\rm PSL}_2(11)$), and they are transposed by the Galois group of this field, generated by complex conjugation.


\begin{figure}[!htp]
\begin{center}

\includegraphics[scale=0.3]{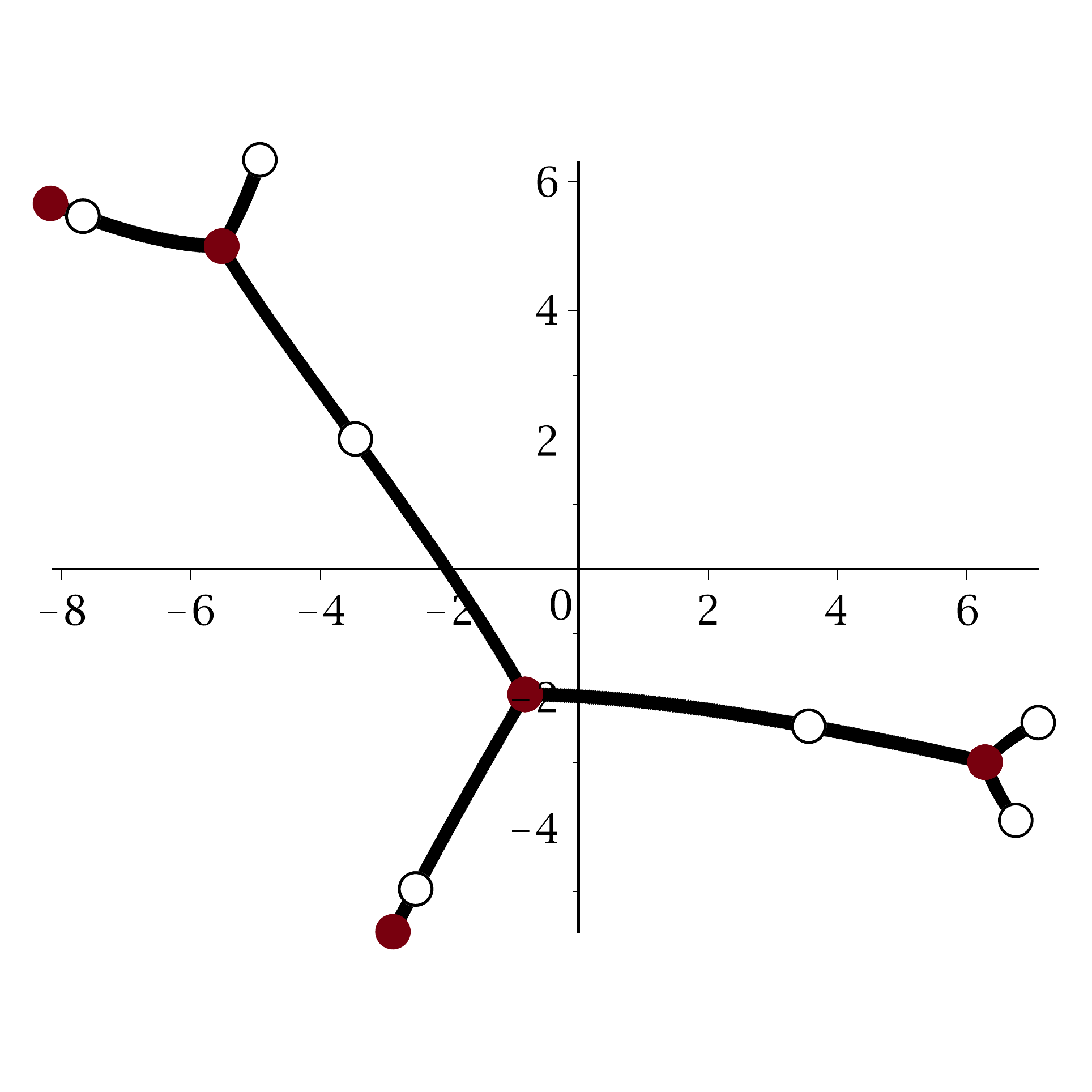}


\caption{A geometrically correct drawing of $\mathcal M_1$.}
\label{Klein-tree-new}
\end{center}
\end{figure}

In each dessin, the black vertices are at the zeros of $P$: the three vertices of valency~$3$ correspond to the zeros of $p_1$ and $p_2$, each appearing with multiplicity $3$ as a zero of~$P$, and the two vertices of valency~$1$ are at the zeros of $p_3$. The white vertices are at the zeros of $P-1$: the four vertices of valency $2$ correspond to the zeros of $q_1$ and $q_2$, each appearing with multiplicity $2$ as a zero of $P-1$, and the three vertices of valency $1$ are the zeros of $q_3$. The unique face-centre is at $\infty$, the pole of $P$ with multiplicity~$11$. For example, Figure~\ref{Klein-tree-new} shows such a geometrically correct version of $\mathcal M_1$ in the plane~$\mathbb C$. 
With a few exceptions (Figures \ref{deg14map}, \ref{loop} and 
\ref{periodicdessin}) all plane trees, dessins and maps shown in this paper are only combinatorially and topologically correct.

The other eight dessins form a single Galois orbit. A computation using Pari/GP produced several defining polynomials for the field of moduli corresponding to this orbit, of which the simplest appears to be
\[Z^8+2Z^6-3Z^5+10Z^4-14Z^3+14Z^2-8Z+1.\]
The Galois group of this field is ${\rm S}_8$, inducing all possible permutations of the eight dessins. The field does not contain ${\mathbb Q}(\sqrt{-11})$, nor, indeed, does it contain ${\mathbb Q}(\sqrt{21})$, even though $\sqrt{21}$ appears in the character table of ${\rm A}_{11}$.


\section{Other related dessins of degree~$11$}
\label{sec:(3,2,r)}

The ten dessins discussed in this paper are not the only dessins of degree~$11$ and type $(3,2,r)$ for some $r$, the least common multiple of the face valencies. It is easy to see that in any such dessin the distribution of valencies of the black vertices must be $3^31^2$, while that of the white vertices must be $2^41^3$ or $2^51^1$. In the former case we obtain the ten dessins based on Klein's diagram, all with $r=11$, but in the latter case we have planar dessins with two faces, rather than one, not considered by Klein.

\begin{figure}[!ht]

\begin{center}
 \begin{tikzpicture}[scale=0.12, inner sep=0.6mm]
  
 \node (a1) at (10,25) [shape=circle, fill=black] {};
 \node (b1) at (10,30) [shape=circle, fill=black] {};
 \node (c1) at (10,35) [shape=circle, fill=black] {};
 \node (d1) at (7,22) [shape=circle, fill=black] {};
 \node (e1) at (13,22) [shape=circle, fill=black] {};

 \draw [thick] (b1) to (a1) to (e1); 
 \draw [thick] (a1) to (d1);
 \draw [thick] (c1) to (10,40);
 \draw [thick] (12.5,32.5) arc (0:360:2.5);
 
 \node at (10,19) {(a)};
 
  
 \node (a2) at (25,25) [shape=circle, fill=black] {};
 \node (b2) at (25,30) [shape=circle, fill=black] {};
 \node (c2) at (25,35) [shape=circle, fill=black] {};
 \node (d2) at (22,22) [shape=circle, fill=black] {};
 \node (e2) at (28,22) [shape=circle, fill=black] {};

 \draw [thick] (b2) to (a2) to (e2); 
 \draw [thick] (a2) to (d2);
 \draw [thick] (c2) to (25,32.5);
 \draw [thick] (27.5,32.5) arc (0:360:2.5);
  
 \node at (25,19) {(b)};
 
 
 \node (a3) at (40,25) [shape=circle, fill=black] {};
 \node (b3) at (50,25) [shape=circle, fill=black] {};
 \node (c3) at (45,35) [shape=circle, fill=black] {};
 \node (d3) at (37,22) [shape=circle, fill=black] {};
 \node (e3) at (53,22) [shape=circle, fill=black] {};
 
\draw [thick]  (a3) to (b3) to (c3) to (a3);
\draw [thick]  (b3) to (e3);
\draw [thick]  (a3) to (d3);
\draw [thick]  (c3) to (45,40);

\node at (45,19) {(c)};

 
 \node (a4) at (65,25) [shape=circle, fill=black] {};
 \node (b4) at (75,25) [shape=circle, fill=black] {};
 \node (c4) at (70,35) [shape=circle, fill=black] {};
 \node (d4) at (62,22) [shape=circle, fill=black] {};
 \node (e4) at (78,22) [shape=circle, fill=black] {};
 
\draw [thick]  (a4) to (b4) to (c4) to (a4);
\draw [thick]  (b4) to (e4);
\draw [thick]  (a4) to (d4);
\draw [thick]  (c4) to (70,29);
 
 \node at (70,19) {(d)};
 
 
 \node (a5) at (0,0) [shape=circle, fill=black] {};
 \node (b5) at (0,5) [shape=circle, fill=black] {};
 \node (c5) at (0,10) [shape=circle, fill=black] {};
 \node (d5) at (5,5) [shape=circle, fill=black] {};
 \node (e5) at (3,-3) [shape=circle, fill=black] {};

 \draw [thick] (c5) to (b5) to (a5) to (e5); 
 \draw [thick] (b5) to (d5);
 \draw [thick] (a5) to (-3,-3);
 \draw [thick] (2,12) arc (0:360:2);
 
 \node at (0,-6) {(e)};
 

 \node (a6) at (15,0) [shape=circle, fill=black] {};
 \node (b6) at (15,5) [shape=circle, fill=black] {};
 \node (c6) at (15,10) [shape=circle, fill=black] {};
 \node (d6) at (20,5) [shape=circle, fill=black] {};
 \node (e6) at (12,-3) [shape=circle, fill=black] {};

 \draw [thick] (c6) to (b6) to (a6) to (e6); 
 \draw [thick] (b6) to (d6);
 \draw [thick] (a6) to (18,-3);
 \draw [thick] (17,12) arc (0:360:2);

 \node at (15,-6) {(f)};
 
 
  \node (a7) at (30,0) [shape=circle, fill=black] {};
 \node (b7) at (30,5) [shape=circle, fill=black] {};
 \node (c7) at (30,10) [shape=circle, fill=black] {};
 \node (d7) at (27,-3) [shape=circle, fill=black] {};
 \node (e7) at (33,-3) [shape=circle, fill=black] {};

 \draw [thick] (c7) to (b7) to (a7) to (e7); 
 \draw [thick] (a7) to (d7);
 \draw [thick] (b7) to (35,5);
 \draw [thick] (32,12) arc (0:360:2);
 
 \node at (30,-6) {(g)};
 
 
 \node (a8) at (45,0) [shape=circle, fill=black] {};
 \node (b8) at (45,5) [shape=circle, fill=black] {};
 \node (c8) at (45,10) [shape=circle, fill=black] {};
 \node (d8) at (45,15) [shape=circle, fill=black] {};
 \node (e8) at (48,-3) [shape=circle, fill=black] {};

 \draw [thick] (b8) to (a8) to (e8); 
 \draw [thick] (c8) to (d8);
 \draw [thick] (a8) to (42,-3);
 \draw [thick] (47.5,7.5) arc (0:360:2.5);
 
 \node at (45,-6) {(h)};
 
 
 \node (a9) at (60,0) [shape=circle, fill=black] {};
 \node (b9) at (60,5) [shape=circle, fill=black] {};
 \node (c9) at (60,10) [shape=circle, fill=black] {};
 \node (d9) at (63,-3) [shape=circle, fill=black] {};
 \node (e9) at (60,7.5) [shape=circle, fill=black] {};

 \draw [thick] (d9) to (a9) to (b9); 
 \draw [thick] (c9) to (e9); 
 \draw [thick] (a9) to (57,-3);
 \draw [thick] (62.5,7.5) arc (0:360:2.5);
 
 \node at (60,-6) {(i)};
 
 
 \node (a10) at (75,0) [shape=circle, fill=black] {};
 \node (b10) at (85,0) [shape=circle, fill=black] {};
 \node (c10) at (80,10) [shape=circle, fill=black] {};
 \node (d10) at (80,4) [shape=circle, fill=black] {};
 \node (e10) at (88,-3) [shape=circle, fill=black] {};
 
\draw [thick]  (a10) to (b10) to (c10) to (a10);
\draw [thick]  (b10) to (e10);
\draw [thick]  (c10) to (d10);
\draw [thick]  (a10) to (72,-3);

 \node at (80,-6) {(j)};

\end{tikzpicture}

\end{center}
\caption{Maps of degree~$11$ and type $(3,2,r)$, $r\ne 11$.} 
\label{othermaps}
\end{figure}
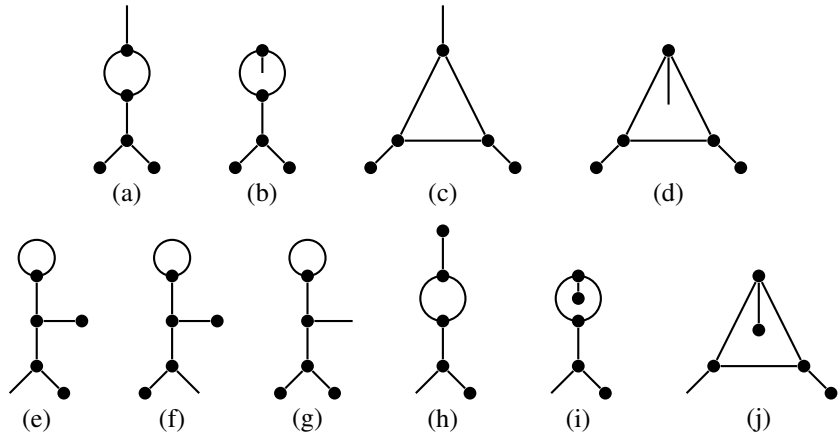

The maps corresponding to these additional dessins can all be formed from the basic map ${\mathcal M}_0$ in Figure~\ref{9dartmap} by adding vertices of valency $1$ to two of its free edges, as in the construction of the maps ${\mathcal M}={\mathcal M}_i$ and $\overline{\mathcal M}_i$, and then joining together two of the remaining three free edges to form a single non-free edge. Equivalently one can perform the latter operation, in three possible ways, on each of the ten maps ${\mathcal M}$, giving $30$ maps of the required form. However, one finds that some of these maps are mutually isomorphic, and that up to isomorphism we obtain just $16$ dessins (confirmed using GAP), consisting of four with bilateral symmetry (in the first row of Figure~\ref{othermaps}) and six chiral pairs (those in the second row and their mirror images). 


There are six where the face valency split is $10+1$, three each where it is $9+2$ or $7+4$, and two each where it is $8+3$ or $6+5$. These dessins correspond to sixteen more conjugacy classes of subgroups of index $11$ in the modular group $\Gamma$, in addition to the ten corresponding to Klein's plane trees. (This is confirmed by~\cite[Theorem~1]{GIR}, which gives a recurrence relation for the number $N_n$  of subgroups of index $n$ in $\Gamma$; here $N_{11}=286$, corresponding to $10+16=26$ conjugacy classes each containing $11$ subgroups.) In all cases except $10+1$, one of the face valencies is a prime $p$, and $z^{11-p}$ is a $p$-cycle with $11-p\ge 3$ fixed points; being transitive of prime degree the monodromy group $G=\langle x, y\rangle$ is primitive, so it follows from Jordan's Theorem and the fact that $y$ is odd that $G={\rm S}_{11}$; in the case $10+1$ we may apply an extension of Jordan's Theorem in~\cite{Jon14} to the cycles of length $10$ and $11$ to show that $G$ acts as ${\rm S}_{11}$ or ${\rm AGL}_1(11)$, with the latter excluded since it has no elements of order $3$. It follows that these dessins correspond to non-congruence subgroups of $\Gamma$. Their regular covers have genus
 \[1+\frac{(r-6)}{12r}11!\]
as given in the following table:

\medskip

\begin{table}[htbp]
\begin{center}
\begin{tabular}{|c||c|c|c|c|c|c|c|c|c|c|c|c|c|c|c|c}
\hline
face valencies & 10, 1 & 9, 2 & 8, 3 & 7, 4 & 6, 5  \\
\hline
$r$ & $10$ & $18$ & $24$ & $28$ & $30$ \\
\hline
genus & $1\,330\,561$ & $2\,217\,601$ & $2\,494\,801$ & $2\,613\,601$ & $2\,661\,121$  \\
\hline
dessins & e, f, g & a, h & b, c & d, i & j \\
\hline
number of dessins & 6 & 3 & 2 & 3 & 2 \\
\hline
\end{tabular}
\caption{Regular coverings of dessins of type $(3,2,r)$.}
\label{tab:regular_3-2-r}
\end{center}
\end{table}


\section{Monodromy groups of dessins of type $(3,2,p)$ and degree $p$}
\label{sec:mono(3,2,p)}

We can generalise our earlier arguments by considering dessins $\mathcal D$ of type $(3,2,p)$ and degree $p$ for {\sl all\/} primes $p$, not just for $p=11$. The monodromy group $G$ of such a dessin $\mathcal D$ must be a  transitive permutation group of degree $p$, and must be perfect if $p>3$, as we will assume, since $3, 2$ and $p$ are mutually coprime. The transitive groups of prime degree $p$ are all known (modulo a very difficult number-theoretic problem which we will discuss in the next section), and are described in~\cite{Cam, DM, Hup}, for example. They are as follows:

\begin{theo}\label{th:primedegree}
The transitive permutation groups $G$ of prime degree are the following:
\begin{itemize}
\item[\rm (a)] 	subgroups $G$ of ${\rm AGL}_1(p)$ containing the translation subgroup;
\item[\rm (b)] 	$G={\rm A}_p$ or ${\rm S}_p$ for a prime $p\ge 5$;
\item[\rm (c)] 	$G={\rm PSL}_2(11)$, ${\rm M}_{11}$ or ${\rm M}_{23}$ with $p=11, 11$ or $23$;
\item[\rm (d)] 	groups $G$ such that ${\rm PSL}_n(q)\le G\le {\rm P\Gamma L}_n(q)$ where the 
				degree $(q^n-1)/(q-1)$ of its natural representation(s) is prime.
\end{itemize}
\end{theo}

By a result of Galois~\cite{Gal}, the solvable groups of prime degree $p$ are those in (a).
Burnside (see~\cite{Bur06} or~\cite[\S251]{Bur}) showed that any nonsolvable group $G$ of prime degree $p$ is doubly transitive; in this case a straightforward argument shows that $S\le G\le{\rm Aut}\,S$ where $S$, the unique minimal normal subgroup of $G$, is a doubly transitive nonabelian simple group. By the classification of finite simple groups, the possibilities for $S$ are
$S={\rm A}_p$ for a prime $p\ge 5$, giving the groups in (b), the three groups $S=G$ in (c), and
$S={\rm PSL}_n(q)$, giving the groups in (d).

In (a) and (b) these are the natural actions of the groups listed. In (c), ${\rm PSL}_2(11)$ has two actions of degree $11$, on the cosets of two conjugacy classes of subgroups isomorphic to ${\rm A}_5$, while the Mathieu groups ${\rm M}_{11}$ and ${\rm M}_{23}$ act on Steiner systems with $11$ and $23$ points. In (d) the groups act on the points of the projective geometry ${\mathbb P}^{n-1}({\mathbb F}_q)$, with a second dual action on its hyperplanes if $n\ge 3$. In particular, since ${\rm Out}\,S$ is solvable in all cases, we have:

\begin{coro}\label{cor:perfect}
The only perfect groups of prime degree are the simple groups ${\rm A}_p$, ${\rm PSL}_2(11)$, ${\rm M}_{11}$, ${\rm M}_{23}$ and\/ ${\rm PSL}_n(q)$ appearing in parts\/ {\rm (b), (c)} and\/ {\rm (d)} of\/ {\rm Theorem~\ref{th:primedegree}}.
\end{coro}

These simple groups $G$ all have elements $x$, $y$ and $z$ of orders $3$, $2$ and $p$ since they have orders divisible by these primes, so there remain the questions of whether there exist such triples with $xyz=1$, whether these can generate $G$, and if so, of determining (or at least enumerating) the corresponding dessins $\mathcal D$. We will address these questions, on a case by case basis, later in this and the following two sections. First we calculate the genera of these maps and of their regular covers.

If $x$ and $y$ have cycle structures $3^a 1^{\alpha}$ and
$2^b 1^{\beta}$, with $\alpha=p-3a$ and $\beta=p-2b$,
then $\mathcal D$ has $p-2a$ black vertices, $p-b$ white vertices, $p$ edges and one face, so it has Euler characteristic
\[\chi=(p-2a)+(p-b)-p+1=p-2a-b+1\,=\frac{4\alpha+3\beta-p}{6}+1\]
and hence genus
\[1-\frac{\chi}{2}=\frac{2a+b+1-p}{2}\;=\frac{p-4\alpha-3\beta+6}{12}.\]
(Note that this implies that $2a+b\ge p-1$, or equivalently $4\alpha+3\beta\le p+6$.) If it exists, such a dessin $\mathcal D$ has a regular cover $\mathcal R$ of type $(3,2,p)$ with automorphism group $G$.
Since $\mathcal{R}$ has $|G|/3$ black vertices of degree 3, $|G|/2$ white
vertices of degree 2, $|G|$ edges, and $|G|/p$ faces of degree $p$, its Euler characteristic is
\[\chi=\left(\frac{1}{3}+\frac{1}{2}+\frac{1}{p}-1\right)\cdot|G|=\frac{6-p}{6p}|G|,\]
so that its genus is
\[g=1-\frac{\chi}{2}=\frac{p-6}{12p}|G|+1.\]


\subsection{Case~(c): sporadic examples}

It is convenient to deal with case~(c) first. We have seen in 
Examples~\ref{ex:PSL_2(11)} and \ref{ex:M_11}, that if $p=11$ then ${\rm M}_{11}$ does not arise as the monodromy group of a dessin of type~$(3,2,p)$ and degree $p=11$, whereas ${\rm PSL}_2(11)$ does, for the chiral pair ${\mathcal M}_1$ and $\overline{\mathcal M}_1$. The Frobenius formula~(\ref{Frobformula}) eliminates ${\rm M}_{23}$: this group contains no triples $(x,y,z)$ of type $(3,2,23)$ such that $xyz=1$.


\subsection{Case~(b): alternating groups} In case~(b) we have $G={\rm A}_p$ for some prime $p$. The following result shows that almost all such groups arise as monodromy groups in the required context.

\begin{theo}\label{degreepmaps}
For each prime $p\ne 2,3,7$ there is a dessin $\mathcal D$ of degree $p$ and type $(3,2,p)$ with monodromy group ${\rm A}_p$.
\end{theo}

\begin{figure}[!ht]

\begin{center}
 \begin{tikzpicture}[scale=0.07, inner sep=0.8mm]

 \node (a) at (-30,0) [shape=circle, fill=black] {};
 \node (b) at (-20,0) [shape=circle, fill=black] {};
 \node (c) at (-10,0) [shape=circle, fill=black] {};
 \node (d) at (10,0) [shape=circle, fill=black] {};
 \node (e) at (20,0) [shape=circle, fill=black] {};
 \node (f) at (30,0) [shape=circle, fill=black] {};
\node (g) at (37,7) [shape=circle, fill=black] {};
 \node (h) at (37,-7) [shape=circle, fill=black] {};

\draw [thick] (a) to (-5,0);
\draw [thick, dashed] (-5,0) to (5,0);
\draw [thick] (5,0) to (f);
\draw [thick] (-37,7) to (a) to (-37,-7);
\draw [thick] (g) to (f) to (h);
\draw [thick] (b) to (-20,10);
\draw [thick] (c) to (-10,10);
\draw [thick] (d) to (10,10);
\draw [thick] (e) to (20,10);

\node at (0,-8) {$\mathcal D$};


 \node (A) at (60,0) [shape=circle, fill=black] {};
 \node (B) at (70,0) [shape=circle, fill=black] {};
 \node (C) at (80,0) [shape=circle, fill=black] {};
 \node (D) at (100,0) [shape=circle, fill=black] {};
 \node (E) at (110,0) [shape=circle, fill=black] {};
 \node (F) at (120,0) [shape=circle, fill=black] {};
\node (G) at (127,-7) [shape=circle, fill=black] {};

\draw [thick] (A) to (85,0);
\draw [thick, dashed] (85,0) to (95,0);
\draw [thick] (95,0) to (F);
\draw [thick] (53,7) to (A) to (53,-7);
\draw [thick] (G) to (F) to (127,7);
\draw [thick] (B) to (70,10);
\draw [thick] (C) to (80,10);
\draw [thick] (D) to (100,10);
\draw [thick] (E) to (110,10);

\node at (90,-8) {${\mathcal D}'$};

\end{tikzpicture}

\end{center}
\caption{Maps $\mathcal D$ and ${\mathcal D}'$ for $p=3k+2$ and $p=3k+1$.}
\label{p=3k+2map}
\end{figure}
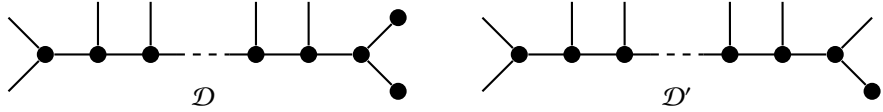

\noindent{\em Proof}. First suppose that $p\equiv 2$ mod~$(3)$, say $p=3k+2$. There is an obvious map for $p=5$, so we may assume that $p\ge 11$, and hence $k\ge 3$. Let $\mathcal D$ be the map on the left in Figure~\ref{p=3k+2map}, where there are $k$ vertices of degree $3$, so that $\mathcal D$ has passport $(3^k1^2, 2^{k+1}1^k, p^1)$. We aim to show that its monodromy group $G$ is ${\rm A}_p$. By Corollary~\ref{cor:perfect} and the discussion of sporadic examples in case~(c) it is sufficient to show that $G$ is not ${\rm PSL}_2(11)$ acting with degree~$11$ or a group 
${\rm PSL}_n(q)$ acting naturally for some $n\ge 2$ and prime power $q$. In fact, when $p=11$, so that $k=3$, $\mathcal D$ is isomorphic to $\overline{\mathcal M}_3$ in Figure~\ref{Kleinmaps}, and we have seen that this has monodromy group ${\rm A}_{11}$. We may therefore suppose, for a contradiction, that $G={\rm PSL}_n(q)$ acting naturally, and that $p\ge 17$, so that $k\ge 5$.

If $\pi$ denotes the permutation character of ${\rm S}_p$, counting fixed points of elements, we see that the involution $y$ corresponding to $\mathcal D$ has fixed point ratio
\[\frac{\pi(y)}{\pi(1)}=\frac{k}{p}=\frac{k}{3k+2}\ge\frac{5}{17}.\]
We will compare this with the corresponding values for involutions in ${\rm PSL}_n(q)$. Any involution in ${\rm PSL}_2(q)$ fixes at most two points, so $n\ne  2$ since $\pi(y)=k\ge 5$. For $n\ge 3$ the involutions $g\in{\rm PSL}_n(q)$ have fixed point ratios $\pi(g)/\pi(1)$ bounded above by
\[b_e=\frac{q^{n-2}+q^{n-3}+\cdots+q+1}{q^{n-1}+q^{n-2}+\cdots+q+1}
\quad\hbox{or}\quad b_o=\frac{(q^{n-3}+\cdots+q+1)+(q+1)}{q^{n-1}+q^{n-2}+\cdots+q+1}\]
as $q$ is even or odd: these bounds are attained when $g$ fixes a hyperplane or two disjoint subspaces of dimensions~$n-3$ and $1$, respectively, corresponding to matrices in ${\rm SL}_n(q)$ with Jordan normal form given by block matrices
\[\Big(\,\begin{array}{cc}
J & O \\
O & I
\end{array}\,\Big)
\quad\hbox{where}\quad 
J=\Big(\,\begin{array}{cc}
 1 & 1 \\ 
 0 & 1
\end{array}\,\Big)
\quad\hbox{or}\quad
\Big(\,\begin{array}{cc}
-1 & 0 \\
 0 &-1
\end{array}\,\Big)
\quad\hbox{and}\quad 
I=I_{n-2}.\]
It is straightforward to check that if $q$ is odd then $b_o<5/17$ and hence $G\ne{\rm PSL}_n(q)$, unless $n=q=3$, when $b_o=5/13$; however, this exceptional case corresponds to $p=13\not\equiv 2$ mod~$(3)$. Thus $q$ is even, so $q=2$ since $b_e<1/q$. This gives $3k+2=p=2^n-1$ and hence $2^n\equiv 0$ mod~$(3)$, which is impossible.

We can apply a similar argument when $p\equiv 1$ mod~$(3)$, say $p=3k+1$. We may assume that $p\ge 13$, that is, $k\ge 4$. By deleting a vertex of valency $1$ from $\mathcal D$ we obtain the map ${\mathcal D}'$ of degree $p$, type $(3,2,p)$ and passport $(3^k1^1, 2^k1^{k+1},p^1)$ shown on the right in Figure~\ref{p=3k+2map}. The corresponding involution $y$ now has fixed point ratio
\[\frac{\pi(y)}{\pi(1)}=\frac{k+1}{p}=\frac{k+1}{3k+1}\ge\frac{5}{13}.\]
As before, if $G={\rm PSL}_n(q)$ then $\pi(y)>2$ implies that $n\ge 3$. We have already seen that if $q$ is odd then $b_o<5/17$, so that $b_o<5/13$ and hence $G\ne{\rm PSL}_n(q)$, unless $n=q=3$; however, in this case $p=13$, and GAP shows that ${\mathcal D}'$ has monodromy group ${\rm A}_{13}$ rather than ${\rm PSL}_3(3)$, so $q$ must be even. Again, $q=2$ since $b_e<1/q$, so if any element of $G={\rm PSL}_n(2)$ has fixed points they must form a projective subspace. Thus $k=\pi(y)$ is odd and hence $p=3k+1$ is even, a contradiction.  \hfill$\square$

\medskip

It is easy to see that there is no corresponding dessin if $p=2$ or $3$. When $p=7$ there is a chiral pair of dessins of the required degree and type, represented by ${\mathcal D}'$ in Figure~\ref{p=3k+2map} with $k=2$, but they have monodromy group ${\rm PSL}_2(7)$, not ${\rm A}_7$ (which is not a Hurwitz group). Simple modifications to the dessins in Figure~\ref{p=3k+2map} show that in the proof of Theorem~\ref{degreepmaps} the number of choices for $\mathcal D$ and $\mathcal D'$ increases exponentially with $p$, even though that construction produces only planar dessins.


\subsection{Case~(d): projective groups} For the rest of this section, we will therefore assume that we are in case~(d), with $G={\rm PSL}_n(q)$ acting on points or hyperplanes of ${\mathbb P}^{n-1}({\mathbb F}_q)$. (This case, together with associated number theoretic estimates and computer searches covered in Section~\ref{sec:BHC}, is investigated in more detail in~\cite{JZ20}.)

\medskip\noindent
{\bf Warning:} here, in spite of an almost universal convention on notation, the unnamed prime of which $q$ is a power is {\sl not\/} the prime $p$; instead, $p$ will denote the degree $(q^n-1)/(q-1)$ of $G$.

\medskip
Before proceeding, we need to make a short digression into number theory. Let us define 
a prime $p$ to be {\em projective\/} if it has the form
\[p=\frac{q^n-1}{q-1}=1+q+q^2+\cdots+q^{n-1}\]
for some prime power $q\ge 2$ and integer $n\ge 2$. These primes include the Fermat primes 
$p=2^{2^k}+1$ with $q=2^{2^k}$ and $n=2$, the Mersenne primes $p=2^n-1$ with  $q=2$ 
and $n$ prime, together with other examples with $q, n\ge 3$, such as $p=13$ with 
$q=n=3$, and $p=31$ with $q=5$ and $n=3$ (notice that $p=31$ is also a Mersenne prime). 
Of course, it is an open problem whether there are infinitely many Fermat or Mersenne primes; at the time of writing, only five Fermat primes (with $k=0,\ldots, 4$) and $51$ Mersenne primes are known to exist. More generally, the existence of infinitely many projective primes seems to be an open problem. In~\cite{JZ20} we give heuristic arguments based on the Bateman--Horn Conjecture, together with associated computational evidence, to support a conjecture that there are infinitely many such primes, even for each fixed prime $n\ge 3$. These are summarised in the next section. However, first we need to prove some necessary conditions for $(q^n-1)/(q-1)$ to be prime.

\begin{lemm}\label{le:irred}
Given integers $n\ge 2$ and $e\ge 1$, the polynomial
\[f(t)=\frac{t^{ne}-1}{t^e-1}=1+t^e+t^{2e}+\cdots+t^{(n-1)e}\]
is irreducible in ${\mathbb Z}[t]$ if and only if $n$ is prime and $e$ is a power 
$n^i\;(i\ge 0)$ of $n$.
\end{lemm}

\noindent{\sl Proof}.
If $k\in\mathbb{N}$ the {\em cyclotomic polynomial}\/ $\Phi_k(x)$ is, by definition,
the polynomial with integer coefficients whose roots are the primitive $k$th roots of unity.
It is irreducible and has degree $\varphi(k)$, where $\varphi$ is the Euler totient function.
For any $n\in\mathbb{N}$ we have $x^n-1=\prod_{d|n}\Phi_d(x)$
(see~\cite[\S 5.2.1]{BS} or \cite[\S 4.3, Problem~26]{NZM}). Putting $x=t^e$ gives
\begin{equation}\label{eq:cyclotomic}
f(t) = \frac{t^{ne}-1}{t^e-1} = \prod_d \Phi_d(t),
\end{equation}
with the product over all $d$ which divide $ne$ but not $e$. 
Thus $f$ is irreducible if and only if there is just one such divisor $d$
(which is $ne$ itself, of course). By considering the prime power decompositions 
of $e$ and $ne$
one can see that this happens if and only if $n$ is prime and $e$ is a power of $n$. 
\hfill$\square$

\begin{lemm}\label{le:necessary}
Let $q=r^e$ where $r$ is prime, and let $n\ge 2$. If $p=(q^n-1)/(q-1)$ is prime then
\begin{itemize}
\item [{\rm(a)}] $n$ is prime and $e$ is a power $n^i\;(i\ge 0)$ of $n$;
\item [{\rm(b)}] $q\not\equiv 1$ {\rm mod}~$(n)$;
\item [{\rm(c)}] ${\rm PSL}_n(q)={\rm SL}_n(q)$.
\end{itemize}
\end{lemm}

\noindent{\sl Proof}.
(a) This follows from Lemma~\ref{le:irred}, using the fact that $p=f(r)=\prod_d\Phi_d(r)$ as in (\ref{eq:cyclotomic}), with each $\Phi_d(r)\ne\pm1$. (To see the latter, note that $\Phi_d(r)=\prod_{\zeta}(r-\zeta)$ where $\zeta$ ranges over the primitive $d$-th roots of unity, so that $|r-\zeta|>1$ for each $\zeta$.)

(b) If $q\equiv 1$ mod~$(n)$ then $p=1+q+\cdots+q^{n-1}\equiv 0$ mod~$(n)$ with $p>n$.

(c) For any $q$ and $n$ we have ${\rm PSL}_n(q)={\rm SL}_n(q)/Z$, where $Z=\{\lambda I\in{\rm GL}_n(q)\mid \lambda^n=1\}\cong{\rm C}_d$ is the centre of ${\rm SL}_n(q)$ and $d=\gcd(q-1,n)$. If $p$ is prime then $d=1$ by (a) and (b). \hfill$\square$


\section{The Bateman--Horn Conjecture}
\label{sec:BHC}

The Bateman--Horn Conjecture~\cite{BH} is a quantified version of Schinzel's Hypothesis H~\cite{SS}. Both of them concern the question of whether a finite set of polynomials $f_1(t),\ldots, f_k(t)\in{\mathbb Z}[t]$ simultaneously take prime values $f_i(t)$ at infinitely many $t\in{\mathbb N}$. (Schinzels's Hypothesis H is in turn a generalisation of the Bunyakovsky Conjecture (1857, see \cite{Bun-1857, Bun-wiki}), which deals with the case $k=1$.) Clearly the following conditions are necessary for this to happen:
\begin{itemize}
\item the leading coefficient of each $f_i$ is positive,
\item each $f_i$ is irreducible,
\item $f:=f_1\ldots f_k$ is not identically zero modulo any prime.
\end{itemize}

In any specific case, the final condition can be verified as follows: $f$ cannot be identically zero modulo any prime $p>\deg(f)$ unless $p$ divides all of its coefficients, and for the remaining primes $p$ one can simply evaluate $f(t)$ at each element $t\in{\mathbb Z}_p$.
Schinzel (following Bunyakovsky in the case $k=1$) conjectured that these conditions are also sufficient. If true, this would imply several other conjectures in Number Theory. For instance, taking $f_1(t)=t^2+1$ would imply the Euler--Landau Conjecture that there are infinitely many primes of the form $t^2+1$, while taking $f_1(t)=t$ and $f_2(t)=t+2$ or $2t+1$ would imply the conjectures that there are infinitely many twin primes and Sophie Germain primes.

The Bateman--Horn Conjecture goes further by proposing a heuristic estimate, based on the Prime Number Theorem, for the number $Q(x)$ of $t\le x$ in $\mathbb N$ such that each $f_i(t)$ is prime. This has the form
\begin{equation}\label{eq:BH-Q}
Q(x)\sim E(x):=\frac{C}{\prod_{i=1}^k\deg f_i}\int_2^x\frac{dt}{(\ln t)^k}
\quad\hbox{as}\quad x\to\infty
\end{equation}
where
\begin{equation}\label{eq:BH-C}
C=C(f_1,\ldots, f_k):=\prod_r\left(1-\frac{1}{r}\right)^{-k}\left(1-\frac{\omega_f(r)}{r}\right)
\end{equation}
with the product over all primes $r$, and $\omega_f(r)$ is the number of $t\in{\mathbb F}_r$ such that $f(t)=0$. If the above conditions on $f_1,\ldots, f_k$ are satisfied then the infinite product in (\ref{eq:BH-C}) converges to a limit $C>0$ (this is far from obvious: see \cite{AFG}, Theorem 5.4.3); since the integral in~(\ref{eq:BH-Q}) diverges as $x\to\infty$, this would imply that $f_1(t),\ldots, f_k(t)$ are simultaneously prime for infinitely many $t\in{\mathbb N}$. So far, these three conjectures have been proved only in the case $k=1$, $\deg f_1=1$: in the case of the Bunyakovsky and Schinzel conjectures, this is Dirichlet's theorem on primes in an arithmetic progression $at+b$; in the case of Bateman--Horn it is equivalent to the stronger form of this theorem, giving an asymptotically equal distribution of primes among the congruence classes of units $b$ mod~$(a)$  (see~\cite[\S5.3.2]{BS} for a proof, and~\cite{Walfisz} for an error estimate).

Nevertheless, the Bateman--Horn Conjecture has provided estimates remarkably close to experimental evidence in cases such as the twin primes conjecture. To give just one example, the number of twin prime pairs up 
to $10^{18}$ is $808\,675\,888\,577\,436$ (entry A007508 of \cite{OEIS}), 
while the Bateman--Horn estimate gives $808\,675\,900\,456\,220$. The relative error 
is $-0.000000147\%$. See also \cite[Figure~10]{AFG}. The constant $C$ is, in this case, 
equal to $1.320323630\ldots$, while $C/2=0.660161815\ldots$ is traditionally called 
the ``twin primes constant''.

\begin{rema}\label{re:Li}\rm
Li~\cite{W.Li} has suggested a modified version of the Bateman--Horn Conjecture
in which $\prod_i\deg f_i$ in the denominator is removed and each factor $1/\ln(t)$ 
within the integral is replaced with $1/\ln(f_i(t))$; the lower limit of integration 
may also be adjusted to avoid a logarithmic singularity at $f_i(t)=1$. 
The result is asymptotically equivalent to the original conjecture, but in cases such 
as Sophie Germain primes, involving a non-monic polynomial $f_i$, the approximation 
is significantly better.
\end{rema}

In our case, if we take $k=2$, with $f_1(t)=t$ and $f_2(t)=1+t^e+t^{2e}+\cdots+t^{(n-1)e}$ for a fixed prime $n\ge 3$ and a fixed exponent $e=n^i$, then $f_2$ is irreducible by Lemma~\ref{le:irred}, so the conjecture gives an estimate $E(x)$ for the number $Q(x)$ of primes $t\le x$ such that $1+q+q^2+\cdots+q^{n-1}$ is prime, where $q=t^e$. Computer searches confirm the obvious guess that $Q(x)$ grows fastest when  $e$ and $n$ take their minimum values $n=3$ and  $e=1$ (so $q$ is prime), so that $f_2(t)=1+t+t^2$. For example,\footnote{The computation was carried out by our colleague Jean B\'etr\'ema \cite{betrema}. He used the system Julia which, on this occasion, turned out to be more efficient than Maple.} all except $301$ of the $1\,974\,311$ projective primes $p\le 10^{18}$ have the form $p=1+q+q^2$ for some prime $q$.

Putting $n=3$ and $e=1$, so that $f=t(1+t+t^2)$, we get
$$
\omega_f(r) \,=\, \left\lbrace 
\begin{array}{ll}
3 & \mbox{if } r\equiv 1 \mbox{ mod } (3), \\
2 & \mbox{for } r=3, \\
1 & \mbox{otherwise}.
\end{array} \right.
$$
We use Maple to calculate $C\approx1.52173006$ (taking the product over all primes 
$r<10^9$) and $\int_2^x(\ln t)^{-2}dt$ (by numerical quadrature) for $x=i\cdot 10^{10}\;(i=1,2,\ldots, 10)$. In Table~\ref{tab:n=3ratios} the resulting estimates $E(x)$ are compared with the true values $Q(x)$ found by computer searches. 
(Here, Li's improvement in Remark~\ref{re:Li} made negligible difference: for instance, using $1/\ln(1+t+t^2)$ instead of $1/2\ln(t)$ with $x=10^{11}$ gave the estimate $129\,297\,407.4$, compared with $129\,297\,407.9$ in Table~\ref{tab:n=3ratios}.) The proportional errors shown are less that $0.00034$, that is, $0.034\%$. Similar comparisons for other values of $n$ and $e$ were reassuring but less convincing, as the numbers of projective primes within our computing range were much smaller: for example, with $n=5$ and $e=1$ we obtained an estimate of $246.718\ldots$\ for the number of projective primes $p=1+q+\cdots+q^4<10^{18}$, whereas the correct number is 252. (Note that in order to estimate the number of primes $p=1+q+q^2\le 10^{18}$ we have to compute the integral in~(\ref{eq:BH-Q}) up to $x=10^{18/2}$; however, for $p=1+q+\cdots+q^4<10^{18}$ the integral is up to $10^{18/4}$, over a much smaller interval.)

\begin{table}[htbp]
\begin{center}
\begin{tabular}{c|c|c|c}
$x$ & $Q(x)$ & $E(x)$ & $E(x)/Q(x)$ \\
\hline
$1 \cdot 10^{10}$ &  15\,801\,827 & $1.579642126 \times 10^7$  & 0.9996579044 \\
$2 \cdot 10^{10}$ &  29\,684\,763 & $2.968054227 \times 10^7$  & 0.9998578150 \\
$3 \cdot 10^{10}$ &  42\,963\,858 & $4.296235691 \times 10^7$  & 0.9999650617 \\
$4 \cdot 10^{10}$ &  55\,877\,571 & $5.587447496 \times 10^7$  & 0.9999445924 \\
$5 \cdot 10^{10}$ &  68\,522\,804 & $6.852175590 \times 10^7$  & 0.9999847043 \\
$6 \cdot 10^{10}$ &  80\,962\,422 & $8.096382889 \times 10^7$  & 1.0000173771 \\
$7 \cdot 10^{10}$ &  93\,236\,613 & $9.323905289 \times 10^7$  & 1.0000261688 \\
$8 \cdot 10^{10}$ & 105\,372\,725 & $1.053741048 \times 10^8$  & 1.0000130940 \\
$9 \cdot 10^{10}$ & 117\,383\,505 & $1.173885689 \times 10^8$  & 1.0000431394 \\
$10^{11}$         & 129\,294\,308 & $1.292974079 \times 10^8$  & 1.0000239757 
\end{tabular}
\end{center}
\vspace{2mm}
\caption{Comparing the BHC estimates $E(x)$ with true values $Q(x)$}
\label{tab:n=3ratios}
\end{table}

On the basis of this evidence, we make the following conjecture:
\begin{conj}
There are infinitely many projective primes of the form $p=1+q+q^2$ where $q$ is prime.
\end{conj}
On the basis of similar evidence for $n>3$ or $e>1$, but necessarily involving smaller numbers of projective primes, we make the following more general conjecture:
\begin{conj}
For each fixed prime $n\ge 3$ and each fixed power $e$ of $n$, there are infinitely many projective primes of the form $p=1+q+\cdots+q^{n-1}$ where $q$ is the $e$-th power of a prime.
\end{conj}

Of course, rather weaker statements than these would also imply the existence of infinitely many projective primes. For example, it would be sufficient to prove that for infinitely many primes $n$ (resp.~primes power $q$) there is at least one prime power $q$ (resp. prime $n$) such that $(q^n-1)/(q-1)$ is prime. This would give an infinite set of such pairs $(n,q)$, with any specific prime represented by at most finitely many of them.


\section{Groups and dessins of projective prime degrees}
\label{sec:dessins(3,2,p)} 

Returning to the groups, the elements~$z$ of order $p=(q^n-1)/(q-1)$ in $G={\rm PSL}_n(q)$ are the Singer cycles (see~\cite[Satz~II.7.3]{Hup}, for example), induced by appropriate non-zero elements of ${\mathbb F}_{q^n}$ acting linearly by multiplication on that field, regarded as an $n$-dimensional vector space over ${\mathbb F}_q$, and permuting its $1$-dimensional subspaces regularly. It follows from results of Berecky~\cite{Ber} that when a Singer cycle $z$ has prime order $p$ the only maximal subgroup of $G$ containing it is the normaliser $N_G(\langle z\rangle)\cong {\rm C}_p\rtimes {\rm C}_n$, where the complement is induced by the Galois group ${\rm Gal}\,{\mathbb F}_{q^n}/{\mathbb F}_q=\langle t\mapsto t^q\rangle\cong {\rm C}_n$. Being solvable, $N_G(\langle z\rangle)$ cannot contain a perfect subgroup $\langle x, y, z\rangle$, so any triple of type $(3,2,p)$, if it exists, must generate $G$. Moreover, there are $(p-1)/n$ conjugacy classes $\mathcal Z$ of Singer cycles $z$ in $G$, one for each orbit of ${\rm Gal}\,{\mathbb F}_{q^n}/{\mathbb F}_q$ on elements of multiplicative order $p$ in ${\mathbb F}_{q^n}$. (Note that since $n$ is prime we have $p-1=q+q^2+\cdots+q^{n-1}\equiv 0$ mod~$(n)$: this is obvious if $q\equiv 0$ or $1$ mod~$(n)$, and otherwise since $n$ is prime the summands represent the non-zero congruence classes mod~$(n)$, so they cancel in pairs.

For the small projective primes $p=5$, $7$, $13$ and $17$ the groups 
${\rm PSL}_2(4)\;(\cong {\rm PSL}_2(5)\cong {\rm A}_5)$, 
${\rm PSL}_3(2)\;(\cong {\rm PSL}_2(7))$, ${\rm PSL}_3(3)$ and ${\rm PSL}_2(16)$ arise as monodromy groups. We will systematically examine various possibilities for monodromy groups ${\rm PSL}_n(q)$, starting with small values of $n$.


\subsection{\bf The case $n=2$: Fermat primes} In the case of the Fermat primes, where $n=2$, we have the following:

\begin{theo}\label{Fermat}
For each Fermat prime $p=2^{2^k}+1\ge 5$ there are $2^{2^k-k-1}$ dessins $\mathcal D$ of type $(3,2,p)$ and degree $p$ with monodromy group $G={\rm PSL}_2(2^{2^k})$, acting naturally. They have genus $(p-5)/12$, while their regular covers $\mathcal R$, which are mutually non-isomorphic, have genus $p(p-4)(p-5)/12$. Dessins $\mathcal D$ and $\mathcal R$ are both mirror symmetric.
\end{theo}

Let $G={\rm PSL}_2(q)$ where $q=2^{2^k}\ge 4$ and $p=q+1$ is prime. The conjugacy classes $\mathcal X$ and $\mathcal Y$ of elements $x$ and $y$ of orders $3$ and $2$ are unique, while there are $q/2$ classes $\mathcal Z$ of elements $z$ of order $p$. Since $3$ and $p$ divide $q-1$ and $q+1$ respectively, the character table of $G$ (see~\cite[\S5.5]{JW}, for example) shows that  only the principal character $\chi=\chi_1$ has $\chi(x)\chi(y)\chi(z)\ne 0$, so for each choice of $\mathcal Z$ we have a character sum $\Sigma=1$. Since $|C(x)|=q-1$, $|C(y)|=q$, $|C(z)|=q+1$ and $|G|=q(q^2-1)$, the number of triples of type $(3,2,p)$ in $G$ is
\[\frac{q}{2}\cdot \frac{|G|^2}{|C(x)|\cdot|C(y)|\cdot|C(z)|}\cdot\Sigma=\frac{q}{2}\cdot|G|.\]

As explained earlier, results of Berecky~\cite{Ber} on Singer cycles imply that each triple generates $G$. We have ${\rm Aut}\,G\cong G\rtimes {\rm Gal}\,{\mathbb F}_q\cong G\rtimes {\rm C}_{2^k}$, acting semi-regularly on these triples, so the number of its orbits on them is
\[\frac{q}{2}\cdot\frac{1}{2^k}=2^{2^k-k-1}.\]
This is therefore the number of normal subgroups of $\Delta=\Delta(3,2,p)$ with quotient group $G$, and thus the number of regular dessins $\mathcal R$ of type $(3,2,p)$ with automorphism group~$G$.

Now $G$ has a unique conjugacy class of subgroups of index $p=q+1$, the stabilisers of points in ${\mathbb P}^1({\mathbb F}_q)$. These lift back to $2^{2^k-k-1}$ conjugacy classes of subgroups of the same index in $\Delta$, so they correspond to that number of dessins $\mathcal D$ of type $(3,2,p)$ and degree $p$ with monodromy group $G$ and regular cover one of the dessins $\mathcal R$. Since $x$, $y$ and $z$ have cycle structures $3^{(q-1)/3}1^2$, $2^{q/2}1^1$ and $p^1$ on ${\mathbb P}^1({\mathbb F}_q)$, the dessins $\mathcal D$ have characteristic
\[\chi=\frac{q-1}{3}+2+\frac{q}{2}+1-(q+1)+1=\frac{16-q}{6}\]
and hence have genus
\[g=1-\frac{\chi}{2}=\frac{q-4}{12}=\frac{p-5}{12}.\]
By the Riemann--Hurwitz formula the regular dessins $\mathcal R$ have genus
\[\frac{p-6}{12p}|G|+1=\frac{p(p-4)(p-5)}{12}.\]

Each of the conjugacy classes $\mathcal Z$ consists of mutually inverse pairs $z, z^q=z^{-1}$, so the maps $\mathcal D$ and $\mathcal R$ are isomorphic to their mirror images. (In fact, as noted by Singerman~\cite{Sin74}, it follows from results of Macbeath~\cite{Macb} that for {\sl any\/} prime power~$q$, each generating pair for ${\rm PSL}_2(q)$ are simultaneously inverted by an automorphism, implying this result more generally.)

\begin{exam}\label{ex:small-proj}\rm
If $k=1$ and $p=5$ we obtain one planar dessin $\mathcal D$, shown on the left in Figure~\ref{deg5,7}, with monodromy group $G={\rm PSL}_2(4)\cong {\rm A}_5$; its regular cover $\mathcal R$ is the dodecahedron, and $\mathcal D$ is the quotient of $\mathcal R$ by a point-stabiliser in $G$, isomorphic to ${\rm A}_4$. If $k=2$ and $p=17$ we obtain two dessins of genus $1$; they are shown in Section~\ref{sec:App}, Figure \ref{fig:psl-2-16}. Their regular covers are of genus $221$. If $k=3$ and $p=257$ we obtain $16$ dessins of genus $21$, with regular covers of genus $1\,365\,441$. If $k=4$ and $p=65\,537$ we obtain $2048$ dessins of genus $5461$, with regular covers of genus $23\,454\,100\,602\,881$.  No Fermat primes are known for $k\ge 5$, and it is conjectured that none exist, but if any do, they are covered by Theorem~\ref{Fermat}.

\begin{figure}[!ht]

\begin{center}
 \begin{tikzpicture}[scale=0.18, inner sep=0.8mm]

 \node (A) at (-5,5) [shape=circle, fill=black] {};
 \node (B) at (0,0) [shape=circle, fill=black] {};
 \node (C) at (5,5) [shape=circle, fill=black] {};

\draw [thick] (A) to (B) to (C);
\draw [thick] (B) to (0,-7);


\node (a) at (15,5) [shape=circle, fill=black] {};
\node (b) at (20,0) [shape=circle, fill=black] {};
\node (c) at (27,0) [shape=circle, fill=black] {};
\draw [thick] (a) to (b) to (c) to (32,-5);
\draw [thick] (b) to (15,-5);
\draw [thick] (c) to (32,5);

\end{tikzpicture}

\end{center}
\caption{Dessins of type $(3,2,p)$ and degree $p$ for $p=5, 7$.}
\label{deg5,7}
\end{figure}
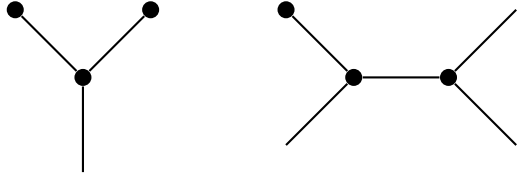

\end{exam}

Having dealt with the case $n=2$, we will now assume that $n\ge 3$. In this case ${\rm PSL}_n(q)$ has two natural representations, on points and hyperplanes of ${\mathbb P}^{n-1}({\mathbb F}_q)$, transposed by the outer automorphism induced by the duality of this geometry. Another difference between this case and the case $n=2$ is that when $n\ge 3$ the conjugacy classes $\mathcal Z$ of Singer cycles are not inverse closed: instead, they form mutually inverse pairs of classes, transposed by duality. This implies that the dessins $\mathcal D$ occur in chiral pairs, each pair having the same regular cover $\mathcal R$, which is regular as a map, that is, isomorphic to its mirror image. 


\subsection{\bf The case $n=3$: projective planes}

The great majority of projective primes, up to a given bound, have $n=3$, so that the group $G={\rm PSL}_3(q)$ acts on the points and lines of the projective plane ${\mathbb P}^2({\mathbb F}_q)$. In this case one can enumerate and describe the associated dessins by using the generic character table for $G$ given by Simpson and Frame in~\cite{SF}, in addition to those in the ATLAS~\cite{ATLAS} for $q\le 9$.

\begin{theo}\label{th:PSL3q}
If $q$ is a prime power such that $p:=1+q+q^2$ is prime, there are $(p-1)/3e$ dessins of type $(3,2,p)$ and degree $p$ with monodromy group ${\rm PSL}_3(q)$ acting naturally, where $q$ is the $e$-th power of a prime. They are all chiral. If $q$ is odd they have passport $(3^{q(q+1)/3}1^1, 2^{(q^2-1)/2}1^{q+2}, p^1)$ and genus $(q-3)(q+1)/12$, whereas if $q=2^e$ their passport is $(3^{q(q+1)/3}1^1, 2^{q^2/2}1^{q+1}, p^1)$ and their genus is $q(q-2)/12$.
\end{theo}

Let $G:={\rm PSL}_3(q)$. By Lemma~\ref{le:necessary}(c) we have $G={\rm SL}_3(q)$. We will apply the Frobenius formula~(\ref{Frobformula}) to $G$, using the generic character table for 
${\rm SL}_3(q)$ in~\cite[Table~1b]{SF}.

First assume that $q$ is coprime to $6$. There are unique conjugacy classes $\mathcal X$ and $\mathcal Y$ of elements of order $3$ and $2$ in $G$, and there are $(p-1)/3=q(q+1)/3$ classes $\mathcal Z$ of elements (Singer cycles) of order $p$; these classes have types ${\mathcal C}_7$, ${\mathcal C}_4$ and ${\mathcal C}_8$ in the notation of~\cite[Table~1a]{SF}, and their elements have centralisers of orders $q^2-1$, $q(q-1)^2(q+1)$ and $p$. The only non-principal irreducible character of $G$ which does not vanish on ${\mathcal X}$, ${\mathcal Y}$ or ${\mathcal Z}$ is the Steinberg character, the unique irreducible character of degree $q^3$, which takes the values $-1$, $q$ and $1$ on them. Since $|G|=q^3(q-1)^2(q+1)(q^2+q+1)$ we find that the number of triples of type $(3,2,p)$ in $G$ is
\[\frac{p-1}{3}\cdot\frac{(q^3(q-1)^2(q+1)(q^2+q+1))^2}{(q^2-1)\cdot q(q-1)^2(q+1)\cdot p}\left(1+\frac{(-1)\cdot q\cdot 1}{q^3}\right)
=\frac{p-1}{3}|G|.\]

By~\cite{Ber} the only maximal subgroup of $G$ containing $z$ is $N_G(\langle z\rangle)$, of order $3p$, so any such triple generates $G$. If $q$ is the $e$-th power of a prime then $|{\rm Out}\,G|=2e$, the two factors corresponding to duality of the projective plane ${\mathbb P}^2({\mathbb F}_q)$ on which $G$ acts and the Galois group of ${\mathbb F}_q$, so $\Delta(3,2,p)$ has $(p-1)/6e$ normal subgroups with quotient $G$. Since $G$ has two conjugacy classes of subgroups of index $p$, the stabilisers of points and lines in ${\mathbb P}^2({\mathbb F}_q)$, we obtain $(p-1)/3e$ dessins of type $(3,2,p)$ and degree $p$ with monodromy group $G$.

In each of the two actions of degree $p$ of $G$, the elements $x$, $y$ and $z$ fix $1$, $q+2$ and $0$ points, so these dessins all have passport $(3^{q(q+1)/3}1^1, 2^{(q^2-1)/2}1^{q+2}, p^1)$ and hence genus $(q-3)(q+1)/12$.

If $q=3^e$ there are two conjugacy classes of elements of order $3$, of types 
${\mathcal C}_2$ and ${\mathcal C}_3$, with centralisers of order $q^3(q-1)$ and $q^2$. 
If $x$ is chosen from the first class, then only the principal character and the constituent 
of the permutation character of degree $p-1=q(q+1)$ contribute to the character 
sum $\Sigma$, and this sum is equal to
\[\Sigma=1+\frac{q\cdot (q+1)\cdot (-1)}{q(q+1)}=0.\]

\noindent
If $x$ is chosen from the second class, then only the principal character contributes to $\Sigma$, so the resulting number of triples in $G$ is
 \[\frac{p-1}{3}\cdot\frac{(q^3(q-1)^2(q+1)(q^2+q+1))^2}{q^2\cdot q(q-1)^2(q+1)\cdot p}=\frac{p-1}{3}|G|\] 
 as before. The number of dessins is therefore again $(p-1)/3e$, and the passport and genus are also as before.

If $q=2^e$ then the main changes are that the unique class $\mathcal Y$ of involutions has type ${\mathcal C}_2$, with centralisers of order $q^3(q-1)$, and that every non-principal irreducible character vanishes on $x$, $y$ or $z$, so $\Sigma=1$. The number of triples of type $(3,2,p)$ in $G$ is therefore
\[\frac{p-1}{3}\cdot\frac{(q^3(q-1)^2(q+1)(q^2+q+1))^2}{(q^2-1)\cdot q^3(q-1)\cdot p}=\frac{p-1}{3}|G|,\]
so the number of dessins is again $(p-1)/3e$. The other change is that $y$ now fixes $q+1$ points, so that the passport is $(3^{q(q+1)/3}1^1, 2^{q^2/2}1^{q+1}, p^1)$ and the genus is $q(q-2)/12$.

In all cases, chirality of the dessins follows from the facts that the two dual representations of $G$ of degree $p$ correspond to a mutually inverse pair of conjugacy classes $\mathcal Z$, and that ${\rm Out}\,G$ acts semiregularly on the set of such pairs.

\begin{exam}\label{ex:p=7}\rm
The smallest example for the above theorem is the chiral pair of dessins of degree 
$p=7$ with monodromy group ${\rm PSL}_3(2)\cong {\rm PSL}_2(7)$; one of this pair 
is shown on the right in Figure~\ref{deg5,7}.
\end{exam}

\begin{exam}\label{ex:p=13}\rm
Let $p=13$ and $G={\rm PSL}_3(3)$. The required conjugacy classes are (in ATLAS~\cite{ATLAS} notation)  ${\mathcal X}=3B$, ${\mathcal Y}=2A$ and ${\mathcal Z}=13A, 13B, 13C$ or $13D$, giving four dessins $\mathcal D$ of type $(3,2,13)$ and degree $13$; they form two chiral pairs, each pair having the same regular cover. Each dessin $\mathcal D$ (see Figure~\ref{PSL_3(3)maps}, or~\cite[Orbit~13.1, Fig.~8.22(A)]{APZ}) has passport $(3^41^1, 2^41^5, 13^1)$ and genus~$0$, while the two regular covers have genus~$253$. (These are the duals of the regular maps R253.1 and R253.2 of type $\{3,13\}_8$ and $\{3,13\}_{12}$ in~\cite{Con}.) 

\begin{figure}[!htp]
\begin{center}

\includegraphics[scale=0.5]{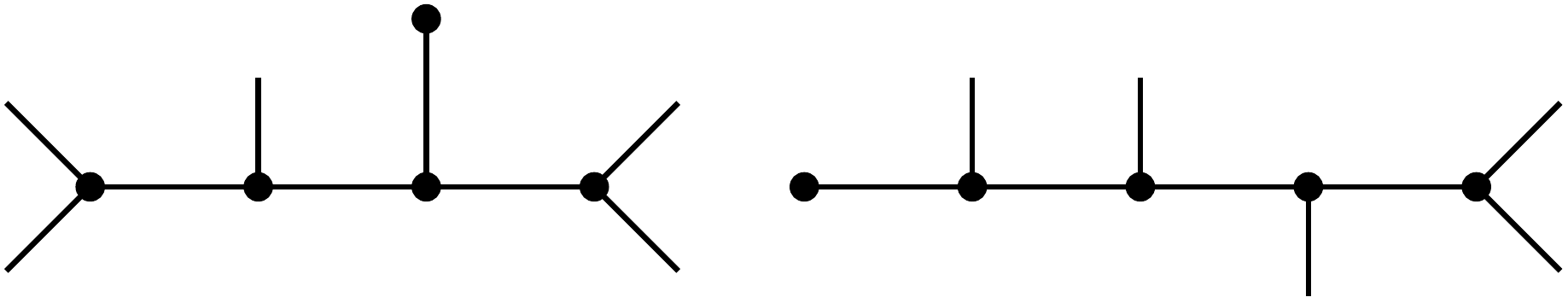}

\caption{Two dessins with monodromy group ${\rm PSL}_3(3)$.}
\label{PSL_3(3)maps}
\end{center}
\end{figure}

\end{exam}

\begin{exam}\label{ex:psl-3-5}\rm
Let $p=31$, with $G={\rm PSL}_3(5)$. In this case ${\mathcal X}=3A$, ${\mathcal Y}=2A$ and there are ten possible classes ${\mathcal Z}=31A, 31B, \ldots$. These give ten dessins $\mathcal D$ of type $(3,2,31)$ and degree $31$, forming five chiral pairs with five regular covers. Each dessin $\mathcal D$ has passport $(3^{10}1^1, 2^{12}1^7, 31^1)$ and genus~$1$, while the regular covers have genus~$25\,001$. Generators of the monodromy groups of these ten dessins $\mathcal D$,
as well as the diagram of one of them, are given in Subsection~\ref{sec:deg-31a}.
\end{exam}

\begin{exam}\label{ex:psl-3-8}\rm
Let $p=73$ and $G={\rm PSL}_3(8)$. Here ${\mathcal X}=3A$ and ${\mathcal Y}=2A$, and there are $24$ classes ${\mathcal Z}=73A,\ldots$ of elements of order $73$. Since $e=3$ we obtain $24/3=8$ dessins $\mathcal D$. They have passport $(3^{24}1^1, 2^{32}1^9, 73^1)$ and genus~$4$; their permutation representations are given in Subsection~\ref{sec:deg-73}. They form four chiral pairs, each pair having a regular cover of genus~$1\,260\,673$.
\end{exam}

Of course, we do not consider cases such as $q=4$, $7$, $9$ and $11$ since the corresponding values of $1+q+q^2$ are not prime. Indeed, squares of primes can be ignored, since for $p$ to be prime we require $e$ to be a power of $n=3$. For any such $e$, the Bateman--Horn conjecture, applied to the polynomials $t$ and $1+t^e+t^{2e}$, suggests that infinitely many projective primes will be obtained. However, even for $e=3$ there are so few examples within our computing range ($10$ primes $p=1+q+q^2<10^{18}$, compared with an estimate of about $12$) that the computational evidence is not convincing.


\subsection{Fixed $n>3$}

In the case of a prime $n>3$ we do not have generic character tables for ${\rm PSL}_n(q)$. Nevertheless, the Bateman--Horn Conjecture still applies, again suggesting that if we fix $n$ and take $e$ to be a fixed power of $n$ we will obtain infinitely many projective primes. The only case in which the computational evidence is persuasive is when $n=5$ and $e=1$, with $252$ projective primes $p<10^{18}$, compared with an estimate of about $246$.


\subsection{The case $q=2$: Mersenne primes}

When $q=2$ we have $G={\rm PSL}_n(2)={\rm GL}_n(2)$, and $p=2^n-1$, currently known to be prime for $51$ values of $n$. No generic character tables for ${\rm GL}_n(2)$ are available, although there is a general theory of characters for ${\rm GL}_n(q)$, largely based on the work of Green~\cite{Gre}; for a concise and readable summary, see~\cite[\S4]{LRS} or~\cite[Ch.~IV]{Macd}. Here are two small examples where character tables in~\cite{ATLAS} can be used.

\begin{exam}\rm
Let $p=2^3-1=7$ and $G={\rm PSL}_3(2)\;(={\rm GL}_3(2)\cong {\rm PSL}_2(7))$. The only conjugacy classes of elements of order $3$ or $2$ in $G$ are $3A$ and $2A$, with cycle structures $3^21^1$ and $2^21^3$ in both of the natural representations of $G$, on points and lines of the Fano plane ${\mathbb P}^2({\mathbb F}_2)$; there are two classes $7A, 7B$ of elements of order $p=7$. The only character which does not vanish in the character sum is $\chi_1$, so $\Sigma=1$.
Since $|C(x)|=3$, $|C(y)|=8=2^3$, $|C(z)|=7$ and $|G|=168=2^3\cdot 3\cdot 7$,
the number of triples of type $(3,2,7)$ in $G$ is
\[2\cdot \frac{2^6\cdot 3^2\cdot 7^2}{3\cdot 2^3\cdot 7}=2^4\cdot 3\cdot 7=2|G|=|{\rm Aut}\,G|.\]
(Note that ${\rm Out}\,G\cong {\rm C}_2$, induced by the duality of the Fano plane, or equivalently by conjugation of ${\rm PSL}_2(7)$ in ${\rm PGL}_2(7)$.) These triples all generate $G$, so $\Delta=\Delta(3,2,7)$ has a unique normal subgroup with quotient $G$, corresponding to a regular dessin $\mathcal R$ of type $(3,2,7)$ with automorphism group $G$. There are two conjugacy classes of subgroups of index $7$ in $G$, the stabilisers of points and lines in the Fano plane; these lift to two conjugacy classes of subgroups of index $7$ in $\Delta$, corresponding to two dessins $\mathcal D$ of type $(3,2,7)$ and degree $7$, each having regular cover $\mathcal R$. Each dessin $\mathcal D$ has passport $(3^21^1, 2^21^3, 7^1)$, so it has genus~$0$, while $\mathcal R$ has genus~$3$.
One of the two dessins $\mathcal D$ is shown on the right in Figure~\ref{deg5,7}, and the other is its mirror image $\overline{\mathcal D}$. These two maps, corresponding to the two choices for the conjugacy class $\mathcal Z$ or to the two choices for the class of subgroups of index~$7$, are not isomorphic since the outer automorphism of $G$ transposes both pairs. The dessin~$\mathcal R$ is Klein's heptagonal map on his quartic curve, introduced in~\cite{Kle78}; it is denoted by $\{7,3\}_8$ in~\cite{CM} and is the dual of R3.1 in~\cite{Con}. As the common regular cover of a chiral pair of dessins (or by its uniqueness) this map is isomorphic to its mirror image.
\end{exam}

\begin{exam}\label{ex:p=31Mersenne}\rm
Let $p=31$ (again), but now in its role as a Mersenne prime $2^5-1$, with 
$G={\rm PSL}_5(2)={\rm GL}_5(2)$. There are two classes each of elements of orders $3$ and $2$, namely $3A$ and $3B$, with cycle structures $3^81^7$ and $3^{10}1^1$, and $2A$ and $2B$, with cycle structures $2^81^{15}$ and $2^{12}1^7$, together with six classes $31A, \ldots$ of elements of order $31$. The inequality $2a+b\ge p-1$ discussed earlier implies that $\mathcal X$ and $\mathcal Y$ must be $3B$ and $2B$. For each of the six choices of $\mathcal Z$, only $\chi_1$ and $\chi_6$ appear in the character sum, so
\[\Sigma=1+\frac{8\cdot(-5)\cdot1}{280}=1-\frac{1}{7}=\frac{6}{7}.\]
We have $|C(x)|=180=2^2\cdot 3^2\cdot 5$, $|C(y)|=1536=2^9\cdot 3$, $|C(z)|=31$ and $|G|=9\,999\,360=2^{10}\cdot 3^2\cdot 5\cdot 7\cdot 31$, so the total number of triples in $G$ is
\[6\cdot 2^9\cdot 3\cdot 5\cdot 7^2\cdot 31\cdot\frac{6}{7}=2^{11}\cdot 3^3\cdot 5\cdot 7\cdot 31=6|G|=3|{\rm Aut}\,G|.\]
Each triple generates $G$, so $\Delta$ has three normal subgroups with quotient $G$. There are two conjugacy classes of subgroups of index $31$ in $G$, the stabilisers of points and hyperplanes in ${\mathbb P}^4({\mathbb F}_2)$, giving six dessins $\mathcal D$ of degree $31$, with three chiral pairs each having the same regular cover $\mathcal R$. Each dessin $\mathcal D$ has passport $(3^{10}1^1,2^{12}1^7,31^1)$, so it has genus~$1$, while the regular covers $\mathcal R$ have genus~$672\,001$. Generators of the monodromy groups of these six dessins $\mathcal D$,
as well as the diagram of one of them, are given in Subsection \ref{sec:deg-31b}.
\end{exam}

\begin{rema}[Goormaghtigh conjecture]\rm
Note that $31$ appears twice as a projective prime, as both $1+5+5^2$ and $1+2+2^2+2^3+2^4$.
In 1917 Goormaghtigh~\cite{Goo} conjectured that this example and $1+2+2^2+\cdots+2^{12}=8191=1+90+90^2$
are the only repetitions among numbers of the form $(q^n-1)/(q-1)$ for integers $q, n\ge 2$.
Although $8191$ is prime (and hence a projective prime for $(q,n)=(2,13)$), 
$90$ is not a prime power, so only the first example is a repetition of projective primes.
A computer search of prime degrees up to $10^{18}$ has revealed no further examples.
The conjecture is still open. For additional information see \cite{Goo-wiki}. 
\end{rema}

Before attempting to find triples $(x,y,z)$ of type $(3,2,p)$ in $G={\rm GL}_n(2)$ for general primes $n\ge 3$, it is useful first to consider the possibilities for elements of orders $3$, $2$ and $p$ in $G$. Any element $x\in G$ of order~$3$ has a rational canonical form 
with $r$ and $(n-r)/2$ blocks of the form
\[\left(\,\begin{array}{c}
1 \\
\end{array}\,\right)
\quad{\rm or}\quad
\left(\,\begin{array}{cc}
0 & 1 \\
1 & 1 \\
\end{array}\right)
\]
for some odd $r<n$, so it fixes a subspace of ${\mathbb P}^{n-1}({\mathbb F}_2)$ of dimension $r-1$; this contains $\alpha=2^r-1$ points, so $x$ has $a=(p-2^r+1)/3$ nontrivial cycles. Similarly, any element $y\in G$ of order~$2$ has a Jordan normal form with $s$ and $(n-s)/2$ blocks of the form
\[\left(\,\begin{array}{c}
1 \\
\end{array}\,\right)
\quad{\rm or}\quad
\left(\,\begin{array}{cc}
1 & 1 \\
0 & 1 \\
\end{array}\right)
\]
for some odd $s<n$; it therefore fixes a subspace of dimension $(n+s-2)/2$ containing $\beta=2^{(n+s)/2}-1$ points, so it has $b=(p-2^{(n+s)/2}+1)/2$ nontrivial cycles. The element $z$ is a $p$-cycle, so provided a triple $(x,y,z)$ of type $(3,2,p)$ exists, it will generate $G$ by~\cite{Ber}; the resulting dessin $\mathcal D$ will have passport $(3^a1^{\alpha}, 2^b1^{\beta}, p^1)$ and genus $(2a+b+1-p)/2$,
while its minimal regular cover $\mathcal R$ will have genus
\[\frac{p-6}{12p}|G|+1=\frac{p-6}{12p}\prod_{i=0}^{n-1}(2^n-2^i)+1=(p-6)2^{(n^2-n-4)/2}\prod_{i=3}^{n-1}(2^i-1)+1,\]
independent of $r$ and $s$. In particular, if we choose the most `active' elements $x$ and $y$ by taking $r=s=1$ (as in the earlier examples with $n=3$ and $5$) we will have $\alpha=1$ and $\beta=2^{(n+1)/2}-1$, so $a=(p-1)/3$ and $b=(p-2^{(n+1)/2}+1)/2$ and hence $\mathcal D$ will have genus
\begin{equation}\label{Mersennegenus}
\frac{1}{3}(2^{n-2}+1)-2^{(n-3)/2}.
\end{equation}

\begin{theo}\label{th:Mersenne}
For each Mersenne prime $p=2^n-1\ge 7$ there are at least $(p-1)/2n$ chiral pairs of dessins of type $(3,2,p)$, degree $p$ and genus given by~{\rm(\ref{Mersennegenus})}, with monodromy group $G={\rm PSL}_n(2)={\rm GL}_n(2)$. 
\end{theo}
 
\noindent
{\em Proof}.
We will use the Frobenius formula~(\ref{Frobformula}) to show that for each conjugacy class $\mathcal Z$ of Singer cycles in $G$ there a triple $(x,y,z)$ of type $(3,2,p)$ where $x$ and $y$ satisfy $r=s=1$ (as above) and $z\in{\mathcal Z}$. The calculation is aided by the fact that most of the non-principal irreducible characters of $G$ vanish on the Singer cycles $z$ (see~\cite[\S4.2]{LRS}, for example), and hence do not contribute to the character sum $\Sigma$. Among the other characters, that of degree $p-1$, the non-principal constituent of the natural permutation character, vanishes on $x$, while the Steinberg character of degree $2^{n(n-1)/2}$ vanishes on $y$, so they can also be ignored. The remaining characters have sufficiently large degrees, compared with their values at $x, y$ and $z$ (see~\cite[\S4.1,4.2]{LRS} for full details), that their total contribution to $\Sigma$ has modulus less than the contribution $1$ of the principal character, so that $\Sigma>0$ and hence the required triples exist. As before, they all generate $G$ by~\cite{Ber}. There are $(p-1)/n$ conjugacy classes $\mathcal Z$ of Singer cycles $z$ in $G$, forming $(p-1)/2n$ mutually inverse pairs; each pair are transposed by the outer automorphism group of $G$, as are the two permutation representations of $G$ of degree $p$; thus we obtain (at least) $(p-1)/2n$ chiral pairs of dessins $\mathcal D$, all with the same passport and genus, as calculated above.
\hfill$\square$
 
On the basis of Theorem~\ref{th:Mersenne}, and of other examples involving small projective primes, we conjecture the following extension of this theorem:
 
\begin{conj}
For each projective prime $p=1+q+\cdots+q^{n-1}>3$ there exists at least one dessin of type $(3,2,p)$ and degree $p$ with monodromy group ${\rm PSL}_n(q)$.
\end{conj}


\subsection{Summary}
We can summarise the results of this rather long section as follows:

\begin{theo}
Suppose that $p$ is a prime such that there is a dessin of type $(3,2,p)$ and degree~$p$ with monodromy group $G$. Then $p\ge 5$, and one of the following holds:
\begin{itemize}
\item[\rm (a)] 	$G={\rm A}_p$, acting naturally, with $p\ne 7$;
\item[\rm (b)] 	$p=11$ and $G={\rm PSL}_2(p)$, acting on the cosets of a subgroup isomorphic 
				to ${\rm A}_5$;
\item[\rm (c)] 	$p$ is a projective prime $(q^n-1)/(q-1)$ for some prime $n$ and prime-power 
				$q$, and $G={\rm PSL}_n(q)={\rm SL}_n(q)$ acting naturally on points or 
				hyperplanes of ${\mathbb P}^{n-1}(q)$.
\end{itemize}
\end{theo}

In case~(a), if $p=5$ there is a unique dessin, shown on the left in Figure~\ref{deg5,7}; if $p>7$ there are exponentially many dessins as $p\to\infty$, even if we restrict to planar dessins. In case~(b) there is a single chiral pair of dessins; these are the dessins ${\mathcal M}_1$ and $\overline{\mathcal M}_1$ studied by Klein in~\cite{Kle79} and shown in Figures~\ref{L2(11)dessins} and ~\ref{Kleinmaps}. In case~(c) we conjecture that there are infinitely many projective primes $p$, even for any fixed prime $n\ge 3$; when $n=2$ (Fermat primes) or $3$ (projective planes) the dessins are all as described in Theorem~\ref{Fermat} or Theorem~\ref{th:PSL3q}. In the case $q=2$ (Mersenne primes) we have an existence result rather than a complete classification.


\section{Modular dessins and curves}
\label{sec:modular}

Just as we generalised our arguments about dessins of type $(3,2,p)$ and degree $p$ from $p=11$ to all primes $p$, we can do the same for dessins of this type and of degree $p+1$. In this section we will consider the dessins ${\mathcal D}_0(p)$ arising naturally from congruence subgroups of level $p$ in the modular group, and in the next section we will consider more general dessins of this type and degree.

\subsection{${\mathcal D}_0(p)$ and its passport}

As in the case $p=11$ considered earlier, for any prime $p$ the natural action of ${\rm PSL}_2(p)$ on the projective line ${\mathbb P}^1({\mathbb F}_p)$ yields a dessin ${\mathcal D}_0(p)$ of type $(3,2,p)$ and degree $p+1$. One can find its passport by using the fact that a non-identity M\"obius transformation (over any field) has at most two fixed points. The cases $p=2$ and $3$ are trivial, so let us assume that $p\ge 5$.

The permutation $x$ representing the black vertices has order $3$, so
\begin{itemize}
\item 	if $p\equiv 1$ mod~$(3)$ then $x$ has two fixed points, and its cycle type 
		is $3^{(p-1)/3}1^2$;
\item	if $p\equiv -1$ mod~$(3)$ then $x$ has no fixed points, and its cycle type 
		is $3^{(p+1)/3}$.
\end{itemize}
Similarly the permutation $y$ representing the white vertices has order $2$, 
and is even, so
\begin{itemize}
\item	if $p\equiv 1$ mod~$(4)$ then $y$ has two fixed points, and its cycle structure 
		is $2^{(p-1)/2}1^2$;
\item	if $p\equiv -1$ mod~$(4)$ then $y$ has no fixed points, and its cycle structure 
		is $2^{(p+1)/2}$.
\end{itemize}
Finally the permutation $z$ representing the faces has order $p$, so its cycle structure 
is~$p^11^1$.

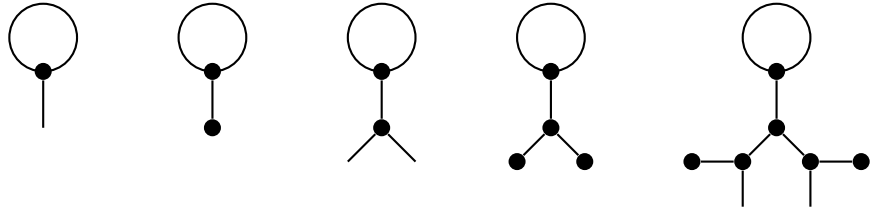
\begin{figure}[!ht]

\vspace*{5mm}

\begin{center}
 \begin{tikzpicture}[scale=0.15, inner sep=0.8mm]
 
 \node (D2) at (-65,5) [shape=circle, fill=black] {};

 \draw [thick] (-65,0) to (D2);
 \draw [thick] (-62,8) arc (0:360:3);
 
 
 \node (B3) at (-50,0) [shape=circle, fill=black] {};
 \node (D3) at (-50,5) [shape=circle, fill=black] {};

 \draw [thick] (B3) to (D3);
 \draw [thick] (-47,8) arc (0:360:3);
 
 
 \node (B5) at (-35,0) [shape=circle, fill=black] {};
 \node (D5) at (-35,5) [shape=circle, fill=black] {};

 \draw [thick] (-38,-3) to (B5) to (-32,-3); 
 \draw [thick] (B5) to (D5);
 \draw [thick] (-32,8) arc (0:360:3);
 

 \node (A7) at (-23,-3) [shape=circle, fill=black] {};
 \node (B7) at (-20,0) [shape=circle, fill=black] {};
 \node (C7) at (-17,-3) [shape=circle, fill=black] {};
 \node (D7) at (-20,5) [shape=circle, fill=black] {};

 \draw [thick] (A7) to (B7) to (C7); 
 \draw [thick] (B7) to (D7);
 \draw [thick] (-17,8) arc (0:360:3);
 

 \node (A13) at (-3,-3) [shape=circle, fill=black] {};
 \node (B13) at (0,0) [shape=circle, fill=black] {};
 \node (C13) at (3,-3) [shape=circle, fill=black] {};
 \node (D13) at (0,5) [shape=circle, fill=black] {};
 \node (E13) at (-7.5,-3) [shape=circle, fill=black] {};
 \node (F13) at (7.5,-3) [shape=circle, fill=black] {};

 \draw [thick] (A13) to (B13) to (C13); 
 \draw [thick] (B13) to (D13);
 \draw [thick] (-3,-7) to (A13) to (E13);
 \draw [thick] (3,-7) to (C13) to (F13);
 \draw [thick] (3,8) arc (0:360:3);

\end{tikzpicture}

\end{center}
\vspace{-5mm}
\caption{${\mathcal D}_0(p)$ for $p=2,3,5,7$ and $13$; white vertices are omitted.}
\label{D0(p)}
\end{figure}

From the passport one can compute the Euler characteristic $\chi$ and hence the genus $g$ of ${\mathcal D}_0(p)$. There are four cases:
\begin{itemize}
\item 	if $p\equiv 1$ mod~$(12)$ the passport is  $(3^{(p-1)/3}1^2, 2^{(p-1)/2}1^2, p^11^1)$, 
		giving the genus $g=(p-13)/12$;
\item 	if $p\equiv 5$ mod~$(12)$ the passport is  $(3^{(p+1)/3}, 2^{(p-1)/2}1^2, p^11^1)$, 
		giving the genus $g=(p-5)/12$;
\item 	if $p\equiv 7$ mod~$(12)$ the passport is  $(3^{(p-1)/3}1^2, 2^{(p+1)/2}, p^11^1)$, 
		giving the genus $g=(p-7)/12$;
\item 	if $p\equiv 11$ mod~$(12)$ the passport is  $(3^{(p+1)/3}, 2^{(p+1)/2}, p^11^1)$, 
		giving the genus $g=(p+1)/12$;
\end{itemize}
Thus the genus of ${\mathcal D}_0(p)$ is $(p-c)/12$ for primes $p\equiv c=-1, 5, 7$ or 
$13$ mod~$(12)$, and~$0$ for $p=2$ or $3$. 
Table \ref{tab:genus-D_0(p)} gives the values of $g$ for the first few primes.
Figure \ref{D0(p)} shows the planar dessins from this table.

\medskip

\begin{table}[htbp]
{\small
\begin{center}
\begin{tabular}{|c||c|c|c|c|c|c|c|c|c|c|c|c|c|c|c|c|} 
\hline
$p$ & 2 & 3 & 5 & 7 & 11 & 13 & 17 & 19 & 23 & 29 & 31 & 37 & 41 & 43 & 47 & 53 \\ 
\hline
$g$ & 0 & 0 & 0 & 0 & 1  &  0 &  1 &  1 &  2 &  2 &  2 &  2 &  3 &  3 &  4 &  4 \\  
\hline
\end{tabular}
\end{center}
}
\caption{The genus of the dessin ${\mathcal D}_0(p)$.}
\label{tab:genus-D_0(p)}
\end{table}

\subsection{The dessin ${\mathcal D}_0(13)$}

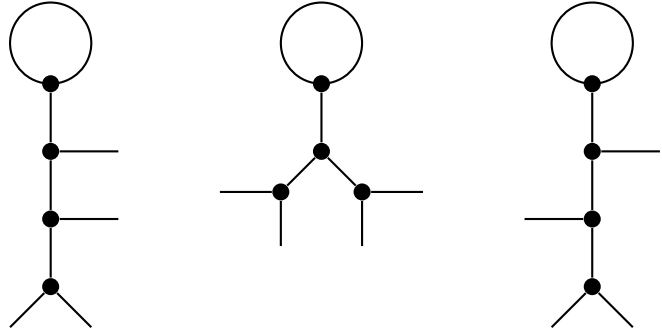
\begin{figure}[ht!]

\begin{center}
 \begin{tikzpicture}[scale=0.18, inner sep=0.8mm]

 \node (A) at (-3,-3) [shape=circle, fill=black] {};
 \node (B) at (0,0) [shape=circle, fill=black] {};
 \node (C) at (3,-3) [shape=circle, fill=black] {};
 \node (D) at (0,5) [shape=circle, fill=black] {};

 \draw [thick] (A) to (B) to (C); 
 \draw [thick] (B) to (D);
 \draw [thick] (-3,-7) to (A) to (-7.5,-3);
 \draw [thick] (3,-7) to (C) to (7.5,-3);
 \draw [thick] (3,8) arc (0:360:3);
 
 \node (b) at (20,0) [shape=circle, fill=black] {};
 \node (d) at (20,5) [shape=circle, fill=black] {};
 \node (e) at (20,-5) [shape=circle, fill=black] {};
 \node (f) at (20,-10) [shape=circle, fill=black] {};
 \draw [thick] (23,8) arc (0:360:3);
 \draw [thick] (d) to (b) to (e) to (f);
  \draw [thick] (b) to (25,0);
    \draw [thick] (e) to (15,-5);
    \draw [thick] (17,-13) to (f) to (23,-13);
    
 \node (b1) at (-20,0) [shape=circle, fill=black] {};
 \node (d1) at (-20,5) [shape=circle, fill=black] {};
 \node (e1) at (-20,-5) [shape=circle, fill=black] {};
 \node (f1) at (-20,-10) [shape=circle, fill=black] {};
 \draw [thick] (-17,8) arc (0:360:3);
 \draw [thick] (d1) to (b1) to (e1) to (f1);
  \draw [thick] (b1) to (-15,0);
    \draw [thick] (e1) to (-15,-5);
    \draw [thick] (-17,-13) to (f1) to (-23,-13);

\end{tikzpicture}

\end{center}
\caption{Three spanning dessins out of which the dessins with the same 
passport as that of ${\mathcal D}_0(13)$, are constructed. White vertices are once again 
omitted.}
\label{3basic}
\end{figure}

The dessin ${\mathcal D}_0(p)$ is planar if and only if $p=2, 3, 5, 7$ or $13$. These five dessins are shown in Figure~\ref{D0(p)}. It is easy to see that the first four are uniquely determined by their passports, and the same applies to the torus dessin ${\mathcal D}_0(11)$ which we will show in Figure~\ref{deg12map}. However, for $p=13$ there are 30 dessins with passport $(3^41^2, 2^61^2, p^11^1)$, consisting of two mirror symmetric dessins, including ${\mathcal D}_0(13)$, and 14 chiral pairs. In fact, a simple argument shows that these dessins can all be formed from the three basic maps in Figure~\ref{3basic}, together with the mirror images of the first and last, by adding black vertices to two of their four free edges. This gives 
$\displaystyle 5\times \left(\!\!\begin{array}{c} 4 \\ 2 \end{array}\!\!\right)=30$ 
dessins.

Whereas  ${\mathcal D}_0(13)$ has monodromy group ${\rm PSL}_2(13)$, the other $29$ dessins all have monodromy group ${\rm A}_{14}$; this follows from the fact that such a group must be a doubly transitive subgroup of ${\rm A}_{14}$, and apart from ${\rm }A_{14}$ itself the only other possibility is ${\rm PSL}_2(13)$, which by the Frobenius formula corresponds to just one dessin, namely ${\mathcal D}_0(13)$.

The fact that ${\mathcal D}_0(13)$ is planar gives us some hope of computing its Bely\u\i\/ function $f:\Sigma\to\Sigma$, which is as follows:
\[f = -\frac{1}{1728}\cdot\frac{(t^4-7t^3+20t^2-19t+1)^3(t^2-5t+13)}{t},\]
\[f-1 = -\frac{1}{1728}\cdot\frac{(t^6-10t^5+46t^4-108t^3+122t^2-38t-1)^2(t^2-6t+13)}{t}.\]
In this case, of course, $f$ is not a polynomial, since there is a pole at $0$ (the face-centre enclosed by the loop) in addition to the pole at $\infty$ (the centre of the outer face).

\begin{figure}[!htp]
\begin{center}

\includegraphics[scale=0.4]{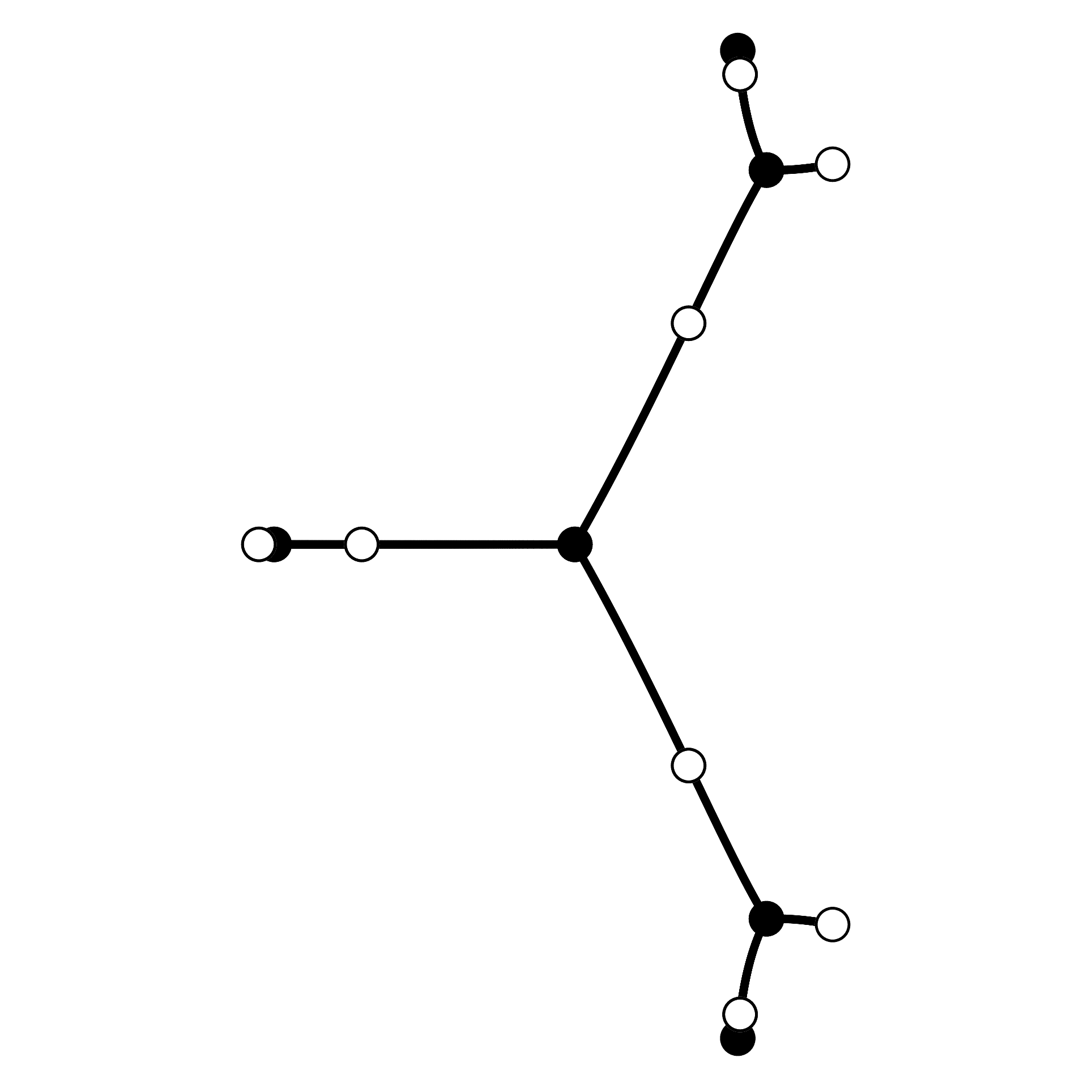}

\caption{A geometrically correct drawing of ${\mathcal D}_0(13)$.}
\label{deg14map}
\end{center}
\end{figure}

\begin{figure}[h!]
\begin{center}

\includegraphics[scale=0.4]{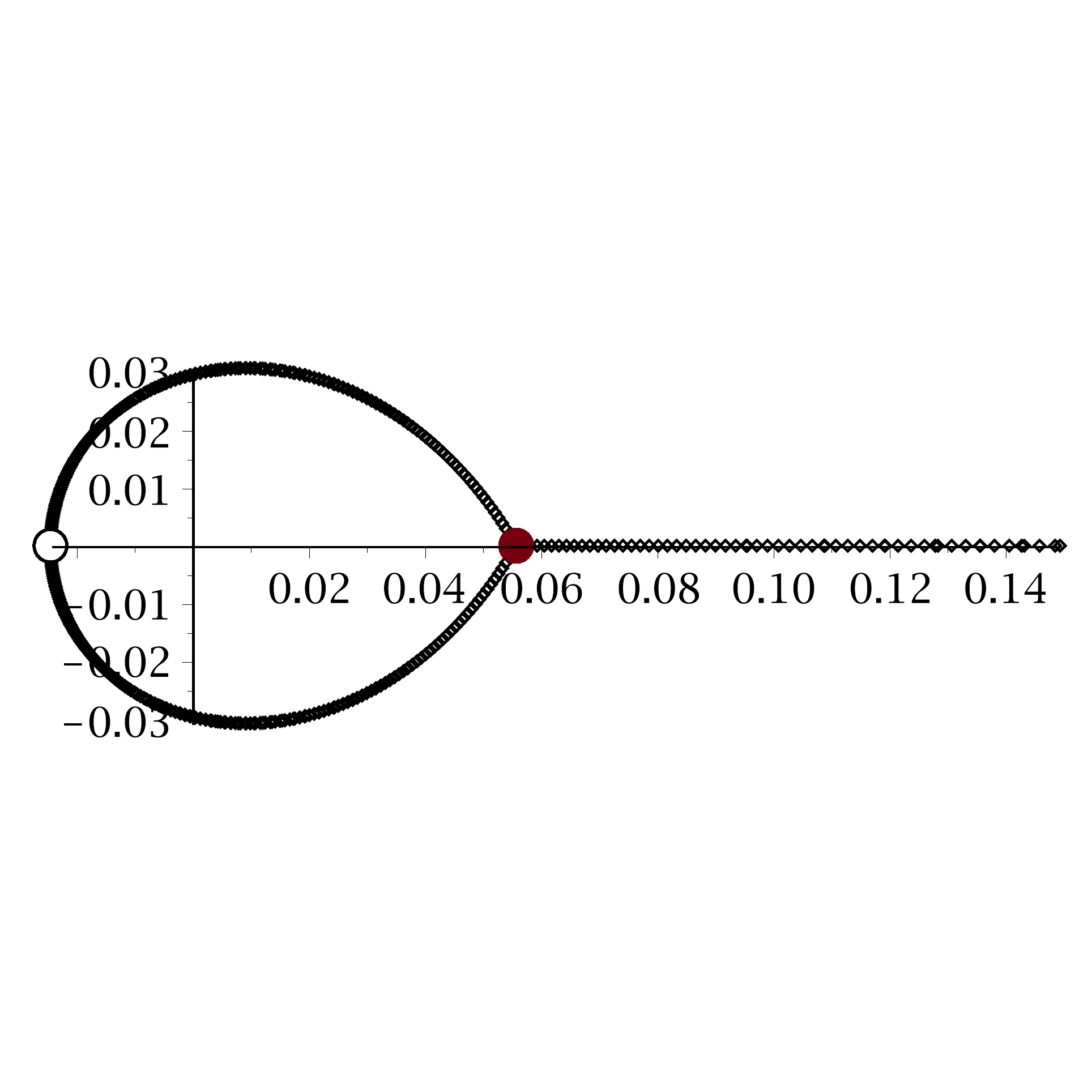}
\vskip -2cm
\caption{The loop in ${\mathcal D}_0(13)$. Its horizontal width is, 
approximately, $0.08$, while the length of the horizontal segment attached to it 
(which goes beyond the limits of the picture) is, approximately, $0.46$.}
\label{loop}
\end{center}
\end{figure}

Figure~\ref{deg14map} shows a geometrically correct version of ${\mathcal D}_0(13)$, with black and white vertices. In this picture, the loop is so small that it is obscured by its incident black and white vertices. (We will meet a similar problem for ${\mathcal D}_0(11)$ in Figure~\ref{periodicdessin}.) Figure~\ref{loop} shows a magnified drawing of the loop, together with nearby parts of the real and imaginary axes.

\subsection{The dessin ${\mathcal D}_0(11)$}

\begin{figure}[htbp]
\begin{center}

\includegraphics[scale=0.5]{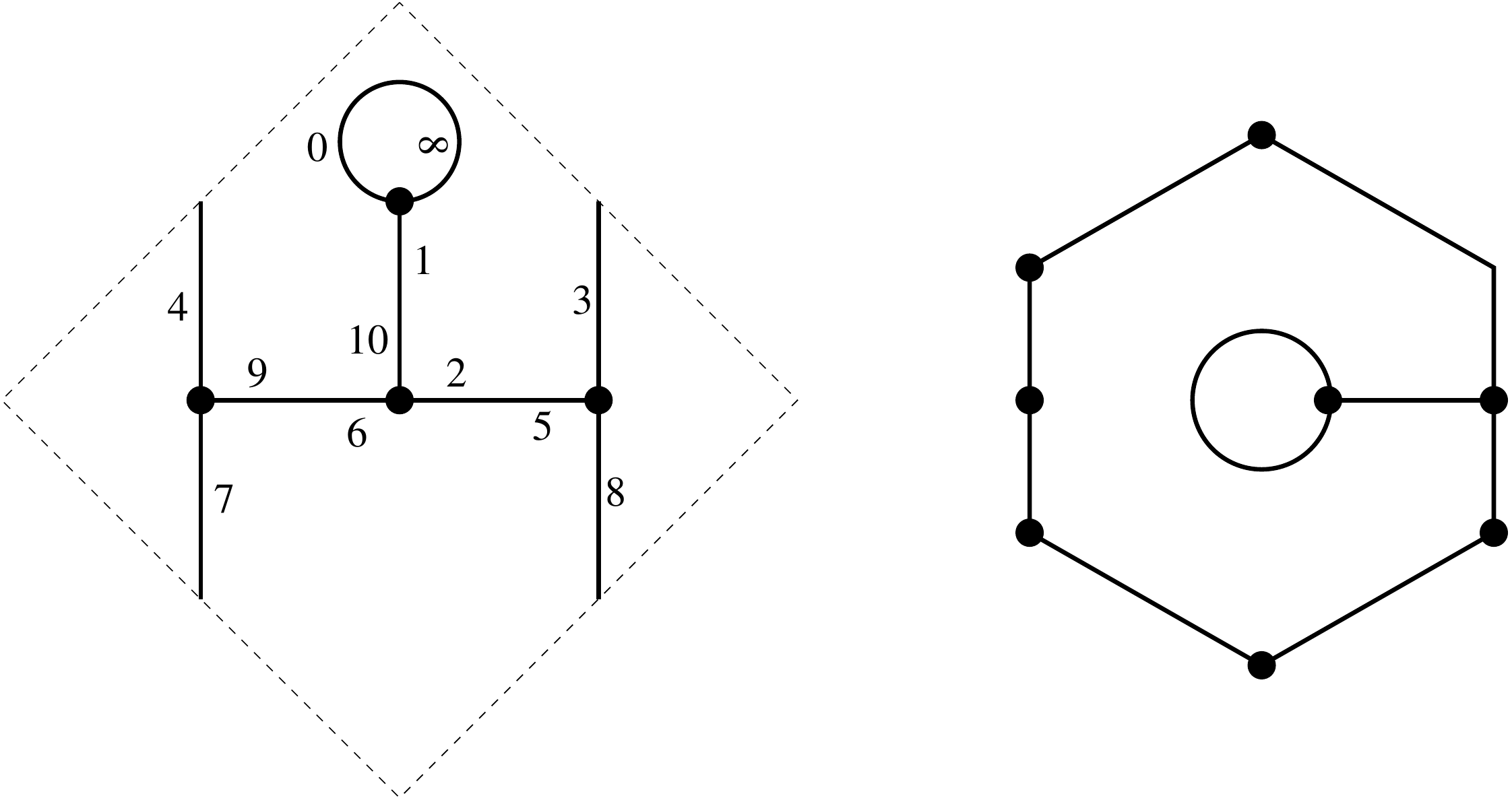}
\end{center}
\caption{Two views of the torus map ${\mathcal D}_0(11)$. The labels of the
half-edges are the elements of $\PP^1(\F_{11})=\F_{11}\cup\{\infty\}$. One may see that
if the label of a half-edge is $t$ than the label of the opposite half-edge of the same
edge is $-1/t$, while the cycle going around the outer face is $(0,1,2,\ldots,10)$.}
\label{deg12map}
\end{figure}

Whereas ${\mathcal D}_0(p)$ is planar for each prime $p\le 13$ except $11$, the dessin ${\mathcal D}_0(11)$ has genus $1$, so we have the additional problem of determining the corresponding elliptic curve. Two views of this dessin are shown in Figure~\ref{deg12map} as maps on a torus, with opposite sides of the outer square or hexagon identified in the usual way. Like the dessins ${\mathcal M}_1$ and $\overline{\mathcal M}_1$, ${\mathcal D}_0(11)$ has the dessin $\mathcal R_1$ as its minimal regular cover. However, in this case it is the quotient of $\mathcal R_1$ by the subgroup of $L\,={\rm PSL}_2(11)$ fixing a point in ${\mathbb P}^1({\mathbb F}_{11})$, of the form ${\rm C}_{11}\rtimes {\rm C}_5$ rather than ${\rm A}_5$. In addition to coming from a more natural representation of $L$, ${\mathcal D}_0(11)$ is also more pleasing because of its greater uniformity: the permutations $x$ and $y$ have no fixed points, so the black vertices all have valency $3$ and there are no free edges. Moreover, the generators $y: t\mapsto -1/t$ and $z: t\mapsto t+1$ of $L$ can both be seen from the edges and faces through the labelling of directed edges with elements of ${\mathbb P}^1({\mathbb F}_{11})$. Since $xyz=1$, the permutation $x$ can now easily be computed as $x:t\mapsto -1/(t-1)$.

Another remarkable feature of this dessin is the fact that it is unique for its passport,
which is
$(3^4,2^6,11^1 1^1)$. One does not have to consider ten different trees and sort out the
corresponding groups. We may wonder why Klein overlooked this example. Maybe the reason
was that the representation of maps via permutations was not yet invented. The first
publications which introduced this construction, by Dyck \cite{Dyck-1888} (1888) and 
Heffter \cite{Heffter-1891} (1891), appeared only a decade after Klein's paper. In fact, 
the pioneer in this subject was Hamilton \cite{Hamilton-1856} (1856), but he explained 
the construction in a private letter. Even at the time close to ours, when the 
permutational model of maps became a part of a common knowledge, the researchers 
studied mainly the automorphism groups of maps while their monodromy groups were 
largely neglected. Only with the advent of the theory of dessins d'enfants the 
monodromy groups rose to the fore since they are Galois invariants.

From the point of view of dessins d'enfants, the uniqueness of the dessin ${\mathcal D}_0(11)$ implies
that it is defined over $\mathbb Q$. However, a difficulty of a different sort appears. 
While in the planar case what we need is a Bely\u\i\/ {\em function}, which is a rational 
function $f:\Sigma\to\Sigma$, in the case of a dessin of genus $g\ge 1$ we need a Bely\u\i\/
{\em pair} $(X,f)$ where $X$ is a Riemann surface (or an algebraic curve) of genus $g$,
and $f$ is a meromorphic function $f:X\to\Sigma$ with no critical values outside
$\{0,1,\infty\}$. There are uncountably many curves of any given genus $g\ge 1$, and the dessin
determines one of them in a unique way. To find such a curve (and a Bely\u\i\/ function
defined on it) is an incredibly difficult task. A relatively simple answer below 
hides numerous obstacles one has to surmount in order to find it. The corresponding 
computation was made for us by John Voight, to whom we are greatly indebted.

In our case $g=1$, so the curve in question is an elliptic curve $E$. One possible
presentation is
\begin{eqnarray}\label{eq:elliptic}
y^2+y=x^3-x^2-10x-20.
\end{eqnarray}
The Bely\u\i\/ function $f=f(x,y)$ is then constructed as follows:
\begin{eqnarray}
a & = & -x^5-23x^4+697x^3-1031x^2-2170x-353; \nonumber \\
b & = & -11x^6+148x^5+643x^4-2704x^3-6780x^2+1781x+3308; \nonumber \\
f & = & \frac{ay+b}{1728(x-16)}. \label{eq:belyi}
\end{eqnarray}
(Alternatively, putting $y=z-\frac{1}{2}$, $x=t+\frac{1}{3}$ gives
$$
z^2 = t^3-\frac{31}{3}t-\frac{2501}{108},
$$
a Weierstrass form for $E$. The $J$-invariant of the curve is
$\displaystyle J(E)=-\frac{2^{12}31^3}{11^5}$.)

Even a verification of the above result is not an easy task. We give an outline
of the main steps of such a verification.

\smallskip

\noindent{\bf Step 0.} Since the map ${\mathcal D}_0(11)$ has 12 directed edges, we expect its Bely\u\i\/ function to have degree~12. However, when one looks at the expression (\ref{eq:belyi}), this property
does not jump to the eyes, to put it mildly. We have an impression that this function
is rather of degree~6. The clue is that it is considered to be defined
on the curve $E$. Let us denote $h(x,y)=y^2+y-x^3+x^2+10x+20$ (see 
equation~(\ref{eq:elliptic})). Expressing $y$ from the equation $f(x,y)=z$ and 
substituting the result in $h$ we indeed get a function of degree~12 in $x$ since 
$h$ contains $y^2$.

\smallskip

\noindent{\bf Step 1.} In order to find critical points on the curve 
we must ensure that the gradients of the functions $f$ and $h$ 
are collinear. Hence the following determinant must be equal to zero:
$$
\left| \begin{array}{cc}
{\displaystyle \frac{\partial f}{\partial x}} & {\displaystyle \frac{\partial f}{\partial y}} \\
& \\
{\displaystyle \frac{\partial h}{\partial x}} & {\displaystyle \frac{\partial h}{\partial y}}
\end{array} \right| \, = \, 0.
$$
This gives us an equation in two variables $x$ and $y$. The second equation is $h(x,y)=0$
since we are looking not for arbitrary points of collineation but for those which belong to $E$.

\smallskip

\noindent{\bf Step 2.} A good idea is to find a Gr\"obner basis for the system thus obtained,
and to do it {\em twice}. Using the `pure lexicographic' monomial order {\tt plex(y,x)} 
we get as the first element of the basis a polynomial in $x$; using {\tt plex(x,y)} we
get a polynomial in $y$. Here they are:
\begin{eqnarray*}
P & = & (x^4-20x^3+62x^2+116x+97)^2(x-16)^2 \times \\
  &   & (x^6-30x^5+243x^4-256x^3-1053x^2+654x+7793), \\
Q & = & (y^4+68y^3+590y^2+556y+7489)^2(y-60)^2 \times \\
  &   & (y^6+102y^5+2619y^4+11272y^3+131211y^2+91074y+253312).
\end{eqnarray*}

\smallskip

\noindent{\bf Step 3.} The roots of these two polynomials are the $x$-coordinates and 
the $y$-coordinates of the critical points of $f$ on $E$. Two factors of degree~4, 
being {\em squared}, indicate that there are four critical points of multiplicity~3. 
We guess that they are the black vertices of our map. But is it true that the value 
of~$f$ at all these points is equal to zero? We must also `couple' $x$- and 
$y$-coordinates in order to show which $y$ correspond to which $x$.

These questions are easy to answer. We solve the system $f=0$, $h=0$ and see that
the values of $x$ are indeed roots of $p_1=x^4-20x^3+62x^2+116x+97$, while 
$$
y = -\frac{1}{23}(x^3-16x^2+90x+62).
$$
It is also easy to verify that these values of $y$ are indeed roots of the polynomial
$q_1=y^4+68y^3+590y^2+556y+7489$. We just substitute these expressions of $y$ into $q_1$
and compute the result modulo $p_1$: it is equal to zero.

Beside these four solutions there is one more: $(x,y)=(16,60)$. What does it mean?
There are two points on $E$ over $x=16$, namely $(16,60)$ and $(16,-61)$. Substituting $y=60$
into the numerator of $f$ we get 
$$
-(x-16)(11x^5+88x^4+2145x^3-4796x^2-8096x-1117).
$$
Thus, $x-16$ in the numerator and in the denominator of $f$ cancel each other:
this is a removable singularity. We conclude that $f$ has a simple pole  
at the point $(16,-61)$ on the curve.

Now, the factors of degree 6 of the polynomials $P$ and $Q$ correspond to simple
critical points, which are mid-points of edges. And, indeed, solving the system 
$f=1$, $h=0$, we obtain, beside the superfluous solution $(16,60)$, the roots of $p_2=x^6-30x^5+243x^4-256x^3-1053x^2+654x+7793$, while
$$
y = -\frac{1}{101617}(104x^5-2811x^4+15943x^3+64714x^2-44258x-217861).
$$
Once again, it is easy to verify that these values of $y$ are roots of the polynomial
\[q_2=y^6+102y^5+2619y^4+11272y^3+131211y^2+91074y+253312.\]

\smallskip

\noindent{\bf Step 4.} Since the function $f$ has only one simple pole on the finite part
of the curve $E$, we may conclude that it also has a pole of multiplicity 11 `over infinity'.
We might stop here. However, we prefer to make the last statement more explicit. 

First of all, let us introduce the projective coordinates $(X\!:\!Y\!:\!Z)$, where $x=X/Z$, 
$y=Y/Z$, and projectivize the equation (\ref{eq:elliptic}) of the curve: 
\begin{eqnarray}\label{eq:projective}
Y^2 Z + YZ^2 = X^3 - X^2 Z - 10XZ^2 - 20Z^3.
\end{eqnarray}
The finite part of the curve corresponds to $Z=1$; the point `at infinity' corresponds to
$Z=0$, which implies also $X=0$; thus, we must take $Y=1$. Near this point, we may introduce 
new affine coordinates $u=X/Y=x/y$, $v=Z/Y=1/y$. Substituting $x=u/v$, $y=1/v$ in the initial 
equation we get the equation of our curve near infinity:
\begin{eqnarray}\label{eq:curve-near-inf}
v + v^2 = u^3 - u^2 v - 10uv^2 - 20v^3.
\end{eqnarray}
Notice that the presence of the term $v$ at the left-hand side of this equation ensures that
the gradient is non-zero and therefore the curve in the vicinity of the point $(u,v)=(0,0)$
is smooth, and $u$ can be taken as a local coordinate. 

Let us represent $v$ as a series $v=c_0 + c_1 u +c_2 u^2 + c_3 u^3 + c_4 u^4+\ldots$ and 
insert it into the difference
$$
(v + v^2) - (u^3 - u^2 v - 10uv^2 - 20v^3).
$$
We would like to get as the result $o(u^3)$. Then, equating to zero the coefficients
in front of the degrees up to the 3rd, we get $c_0=c_1=c_2=0$, $c_3=1$. What remains is to make
the following two operations:
\begin{itemize}
\item[1.]	Substitute $x=u/v$, $y=1/v$ in the expression (\ref{eq:belyi}) for the Bely\u{\i} 
			function.
\item[2.]	Substitute $v=u^3$ in the resulting expression.
\end{itemize}
Then the Bely\u{\i} function, up to smaller terms, becomes
$$
f \, \sim \, -\frac{1}{1728}\frac{1+11u+23u^2-148u^3-\ldots-3308u^{13}}{u^{11}(1-16u^2)} 
\, \sim \,
\frac{\rm Const}{u^{11}} \quad \mbox{when} \quad u\to 0.
$$
Now a pole of degree 11 at the point $u=0$ is apparent.

\medskip

\begin{figure}[!htp]
\begin{center}


\includegraphics[scale=0.4]{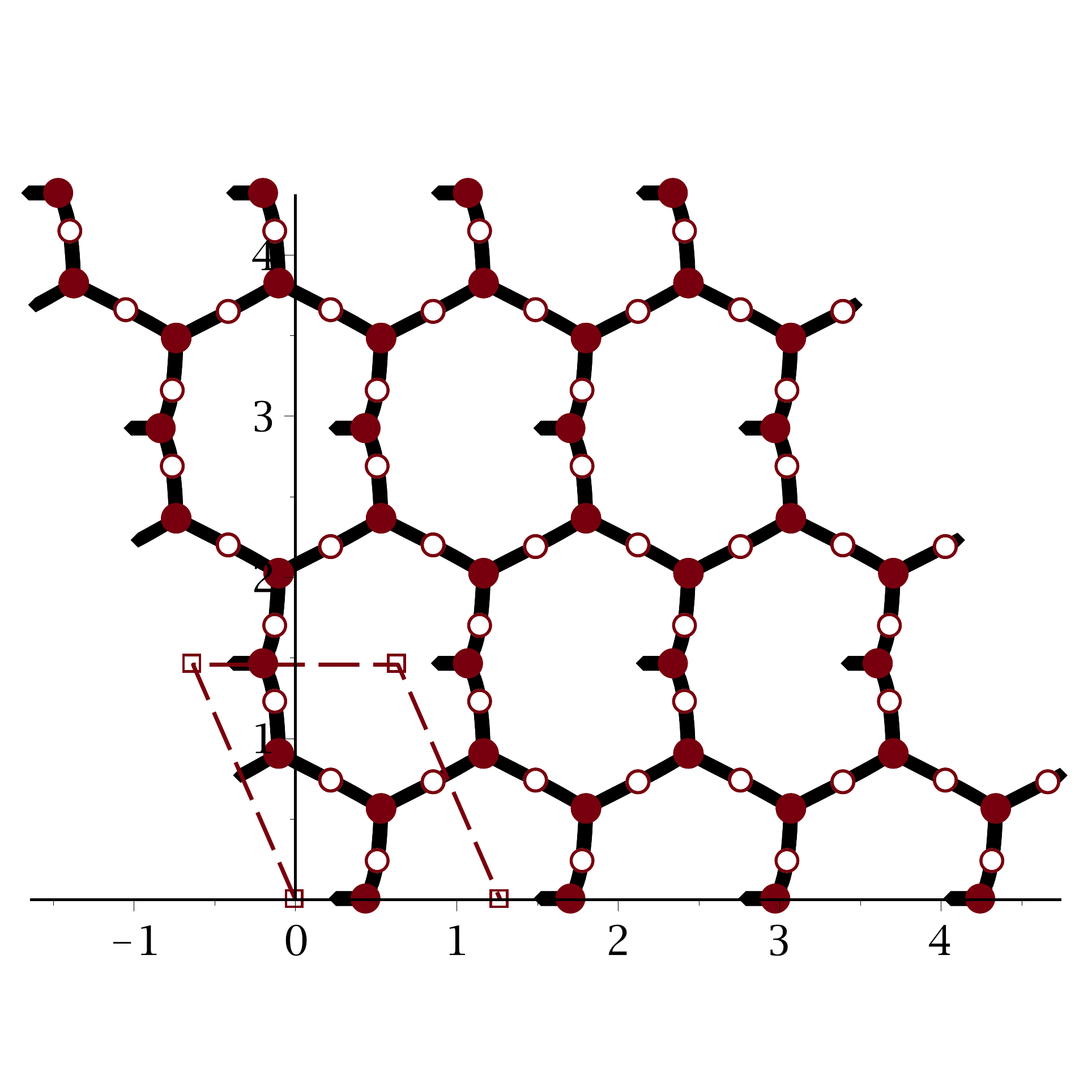}


\caption{A geometrically correct drawing of ${\mathcal D}_0(11)$ 
lifted to the complex plane.}
\label{periodicdessin}
\end{center}
\end{figure}

Figure~\ref{periodicdessin} shows a geometrically correct picture of the dessin ${\mathcal D}_0(11)$, lifted from the torus $E$ to the complex plane to give a doubly periodic map. In addition to the black vertices of degree $3$, the white vertices of degree~$2$ are also shown, as is a fundamental parallelogram for the associated lattice. The repeated copies of ${\mathcal D}_0(11)$ are aligned as in the right-hand map in Figure~\ref{deg12map}. However, the loop (the face of valency~$1$), which is obvious there, is in reality too small to be clearly visible here. Changing accuracy parameters to make this face and its incident vertices visible would create a mess elsewhere, so in this case we have to accept that nature, technology and human eyesight are incompatible. We met a similar problem earlier for ${\mathcal D}_0(13)$, in connection with Figures~\ref{deg14map} and \ref{loop}. 

To obtain this diagram we use the substitution
\[y = \frac{z-1}{2}, \quad x = t+\frac{1}{3}\]
to give the elliptic curve 
\[z^2=4t^3-\frac{124}{3}t-\frac{2501}{27}\]
in Weierstrass normal form.
Then Algorithm~7.4.7 in Cohen's book~\cite{Coh}, using the arithmetic-geometric mean, gives a basis
\[\omega_1 = 1.269209304, \quad \omega_2 = -0.6346046520 + 1.458816617 i\]
for the lattice, as shown in Figure~\ref{periodicdessin}. The associated modulus is
\[\tau = \omega_2/\omega_1 = -0.5000000000+1.149390106 i,\]
on the left-hand vertical border of the standard fundamental region for $\Gamma$, which is consistent with the fact that $J(\tau)<0$. To obtain the embedded graph we lift the unit interval $[0,1]\subset\Sigma$, with black and white vertices at $0$ and $1$, to the curve $E$ via the Bely\u\i\/ function $f$, and then use the standard parametrisation of $E$ by the Weierstrass functions, namely $t=\wp(w)$ and $z=\wp'(w)$ for $w\in\mathbb C$, to lift the resulting dessin to the plane.


\subsection{The dessins ${\mathcal D}(p)$}

The minimal regular cover of ${\mathcal D}_0(p)$ is a regular dessin ${\mathcal D}(p)$ of type $(3,2,p)$ with automorphism group ${\rm PSL}_2(p)$; for $p>2$ it has degree $p(p^2-1)/2$ and genus $(p+2)(p-3)(p-5)/24$, while for $p=2$ it has degree $6$ and genus $0$. The underlying curves of ${\mathcal D}(p)$ and ${\mathcal D}_0(p)$ are the modular curves $X(p)=\overline{{\mathbb H}/\Gamma(p)}$ and $X_0(p)=\overline{{\mathbb H}/\Gamma_0(p)}$ (hence the notation for these two dessins); here $\Gamma(p)$ is the principal congruence subgroup of level $p$ in $\Gamma$, the kernel of the reduction mod~$(p): \Gamma\to {\rm PSL}_2(p)$ (see~\cite[\S 6.9--10]{JS}, for example), while $\Gamma_0(p)$ is the inverse image in $\Gamma$ of the subgroup of ${\rm PSL}_2(p)$ fixing $\infty$.

\begin{theo}\label{th:D(p)}
For each prime $p$ the dessin ${\mathcal D}(p)$ is the only regular dessin of type $(3,2,p)$ 
with automorphism group \linebreak $G={\rm PSL}_2(p)$, and ${\mathcal D}_0(p)$ is the only 
dessin of degree $p+1$ and type $(3,2,p)$ with monodromy group $G$.
\end{theo}

To prove the uniqueness of ${\mathcal D}(p)$ it is sufficient to show that $G$ has, up to automorphisms, only one generating triple $(x,y,z)$ of type $(3,2,p)$. The elements $z\in G$ of order $p$ are all equivalent under automorphisms, while the elements of order $3$ and $2$ are those with traces $\pm 1$ and $0$. We may therefore assume that
\[
y=\pm\left(\,\begin{array}{cc}
a & b \\
c & -a \\
\end{array}
\right)
\quad\hbox{and}\quad
z=\pm\left(\,\begin{array}{cc}
1 & 0 \\
1 & 1 \\
\end{array}
\right),
\quad\hbox{so that}\quad
yz=\pm\left(\,\begin{array}{cc}
a+b & b \\
c-a & -a \\
\end{array}
\right).\]
Since we require $yz=x^{-1}$ to have order $3$ we must have $b=\pm 1$, and without loss of generality we can take $b=1$. Since $\det y=1$ we have $c=-a^2-1$, giving a $1$-parameter set of solutions
\[
y=y_a=\pm\left(\,\begin{array}{cc}
a & 1 \\
-a^2-1 & -a \\
\end{array}
\right)
\quad (a\in{\mathbb F}_p).
\]
Now conjugating the triple $(x,y,z)$ with $z^i$ fixes $z$ and replaces $y=y_a$ with $y_{a+i}$, so all such triples are equivalent under automorphisms, as required. One could also prove this by using the Frobenius triple-counting formula~(\ref{Frobformula}), as in 
Example~\ref{ex:PSL_2(11)} 
for $p=11$ (see~\cite[\S5.5]{JW} for the character table of $G$); this is easy if $p\equiv \pm 5$ mod~$(12)$ since in this case only one nonprincipal irreducible character, of degree $p$, appears in the character sum, but if $p\equiv\pm 1$ mod~$(12)$ then other characters appear and some work with algebraic numbers is required. Since $G$ has a unique conjugacy class of subgroups of index $p+1$, namely the stabilisers of points in ${\mathbb P}^1({\mathbb F}_p)$, it follows that ${\mathcal D}_0(p)$ is the only dessin of degree $p+1$ and type $(3,2,p)$ with monodromy group $G$.


\section{Dessins of type $(3,2,p)$ and degree $p+1$}
\label{sec:dessinsp+1}

Having considered the dessins ${\mathcal D}_0(p)$ in the preceding section, we will now consider arbitrary dessins of their type $(3,2,p)$ and degree $p+1$. The monodromy group $G$ of such a dessin $\mathcal D$ is a transitive permutation group of degree $p+1$, and since it contains a $p$-cycle $z$ it is doubly transitive, and hence primitive. As before, if $p>3$, as we shall assume, $G$ is perfect. The finite primitive groups containing a cycle with one fixed point were determined by M\"uller in~\cite[Theorem~6.2]{Mul}, and as a corollary we have:

\begin{theo}\label{th:degp+1}
The perfect permutation groups of degree $p+1$ {\rm (}$p$ prime\/{\rm )}, containing 
a $p$-cycle, are the following:
\begin{itemize}
\item[\rm (a)] 	${\rm A}_{p+1}$ for primes $p\ge 5$,
\item[\rm (b)] 	${\rm PSL}_2(p)$ for primes $p\ge 5$,
\item[\rm (c)] 	affine groups\/ ${\rm AGL}_n(2)$\, for Mersenne primes $p=2^n-1\ge 5$,
\item[\rm (d)] 	Mathieu groups\/ ${\rm M}_{11}$\, and\/ ${\rm M}_{12}$\, for $p=11$ and\/ 
				${\rm M}_{24}$\, for $p=23$, 
\end{itemize}
all except ${\rm M}_{11}$ {\rm (}see below{\rm )} in their natural representations.
\end{theo}

We will consider these groups in turn, together with their analogues for $p<5$, 
to see which of them can serve as monodromy groups of dessins of type $(3,2,p)$.
The results will be summarized in Theorem~\ref{th:summary}. It is convenient 
to deal with the simplest cases (b) and (d) first, and to treat the alternating groups last.

The groups ${\rm PSL}_2(p)$ in Theorem~\ref{th:degp+1}(b) arise as monodromy groups in this context for all primes~$p$, each of them associated with a unique dessin, the modular dessin ${\mathcal D}_0(p)$ (see Theorem~\ref{th:D(p)}).

In (d), although ${\rm M}_{11}$ has a permutation representation of degree $12$, on the cosets of a subgroup ${\rm PSL}_2(11)$, we have seen in Example~\ref{ex:M_11} that it is not a quotient of $\Delta(3,2,11)$, since all triples of this type in ${\rm M}_{11}$ generate not ${\rm M}_{11}$ but its subgroup ${\rm PSL}_2(11)$; thus it does not arise as a monodromy group in this context.

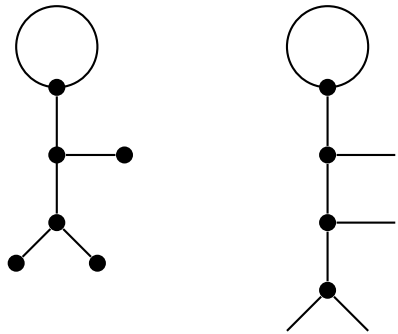
\begin{figure}[!ht]

\begin{center}
 \begin{tikzpicture}[scale=0.18, inner sep=0.8mm]

  \node (b) at (20,0) [shape=circle, fill=black] {};
 \node (d) at (20,5) [shape=circle, fill=black] {};
 \node (e) at (20,-5) [shape=circle, fill=black] {};
 \node (f) at (20,-10) [shape=circle, fill=black] {};
 \draw [thick] (23,8) arc (0:360:3);
 \draw [thick] (d) to (b) to (e) to (f);
 \draw [thick] (b) to (25,0);
 \draw [thick] (e) to (25,-5);
 \draw [thick] (23,-13) to (f) to (17,-13);
 
 
 \node (A) at (-3,-8) [shape=circle, fill=black] {};
 \node (B) at (0,-5) [shape=circle, fill=black] {};
 \node (C) at (3,-8) [shape=circle, fill=black] {};
 \node (D) at (0,5) [shape=circle, fill=black] {};
 \node (E) at (0,0) [shape=circle, fill=black] {};
 \node (F) at (5,0) [shape=circle, fill=black] {};

 \draw [thick] (A) to (B) to (C); 
 \draw [thick] (B) to (D);
\draw [thick] (E) to (F);
 \draw [thick] (3,8) arc (0:360:3);

\end{tikzpicture}

\end{center}
\caption{Two dessins with monodromy group ${\rm M}_{12}$. The field of moduli for
both of them is $\mathbb{Q}(\!\sqrt{-11})$.}
\label{M12dessins}
\end{figure}

The group $G={\rm M}_{12}$ in (d) has two representations of degree $12$, on the cosets of two conjugacy classes of subgroups isomorphic to ${\rm M}_{11}$, transposed by ${\rm Out}\,G$. It is the monodromy group of two chiral pairs of planar dessins of type $(3,2,11)$ and degree $12$. One dessin from each pair is shown in Figure~\ref{M12dessins}: the dessin on the left has passport $(3^31^3, 2^6, 11^11^1)$, with ${\mathcal X}=3A$ and ${\mathcal Y}=2A$, while that on the right has passport $(3^4, 2^41^4, 11^11^1)$, with ${\mathcal X}=3B$ and ${\mathcal Y}=2B$. Their regular covers are a pair of dessins of type~$(3,2,11)$ and genus~$3601$. There are no other dessins of this type and degree with monodromy group $G={\rm M}_{12}$: taking ${\mathcal X}=3B$ and ${\mathcal Y}=2A$, or ${\mathcal X}=3A$ and ${\mathcal Y}=2B$ yields two more orbits of ${\rm Aut}\,G$ on triples, but these generate transitive or intransitive subgroups isomorphic to ${\rm PSL}_2(11)$. The dessin corresponding to these transitive subgroups has already been shown in Figures~\ref{deg12map} and \ref{periodicdessin}, while the intransitive subgroups, with orbits of length $11$ and $1$, are the monodromy groups of the dessins  ${\mathcal M}_1$ and $\overline{\mathcal M}_1$ of degree~$11$ in Figures~\ref{Klein'sdessins}, \ref{L2(11)dessins} and \ref{Kleinmaps}.

\begin{figure}[!ht]

\begin{center}
 \begin{tikzpicture}[scale=0.15, inner sep=0.8mm]

\node (a) at (-10,0) [shape=circle, fill=black] {};
\node (b) at (-5,0) [shape=circle, fill=black] {};
\node (c) at (0,0) [shape=circle, fill=black] {};
\node (d) at (5,0) [shape=circle, fill=black] {};
\node (e) at (10,0) [shape=circle, fill=black] {};
\node (f) at (15,0) [shape=circle, fill=black] {};
\node (B) at (-5,-5) [shape=circle, fill=black] {};
\node (C) at (0,-5) [shape=circle, fill=black] {};
\node (D) at (5,-5) [shape=circle, fill=black] {};
\node (E) at (10,-5) [shape=circle, fill=black] {};
\node (C1) at (-3,-8) [shape=circle, fill=black] {};
\node (C2) at (3,-8) [shape=circle, fill=black] {};

 \draw [thick] (-10,0) arc (0:360:3);
 \draw [thick] (a) to (f);
 \draw [thick] (b) to (B);
 \draw [thick] (c) to (C);
 \draw [thick] (d) to (D);
 \draw [thick] (e) to (E);
 \draw [thick] (C1) to (C) to (C2);
 
 
 \node (a1) at (30,0) [shape=circle, fill=black] {};
\node (b1) at (35,0) [shape=circle, fill=black] {};
\node (c1) at (40,0) [shape=circle, fill=black] {};
\node (d1) at (45,0) [shape=circle, fill=black] {};
\node (e1) at (50,0) [shape=circle, fill=black] {};
\node (f1) at (55,0) [shape=circle, fill=black] {};
\node (g1) at (60,0) [shape=circle, fill=black] {};
\node (E1) at (50,-5) [shape=circle, fill=black] {};

 \draw [thick] (30,0) arc (0:360:3);
 \draw [thick] (a1) to (g1);
 \draw [thick] (b1) to (35,-5);
 \draw [thick] (c1) to (40,-5);
 \draw [thick] (d1) to (45,5);
 \draw [thick] (e1) to (E1);
 \draw [thick] (f1) to (55,5);
 \draw [thick] (47,-8) to (E1) to (53,-8);
 \draw [thick] (63,3) to (g1) to (63,-3);

\end{tikzpicture}

\end{center}
\caption{Two dessins with monodromy group ${\rm M}_{24}$. The field of moduli for
both of them is $\mathbb{Q}(\!\sqrt{-23})$.}
\label{M24dessins}
\end{figure}
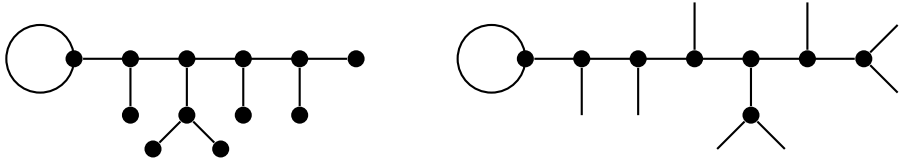

For $p=23$, the group ${\rm M}_{24}$ in (d), found by Mathieu~\cite{Mat73} in 1873, is the monodromy group of two chiral pairs of planar dessins, with passports $(3^61^6, 2^{12}, 23^11^1)$ and $(3^8, 2^81^8, 23^11^1)$: see orbits 24.2 and 24.1 in~\cite[pp.~129--130]{APZ}. For the first pair $\mathcal X=3A$, $\mathcal Y=2B$ and $\mathcal Z=23A$ or $23B$, while for the second pair $\mathcal X=3B$, $\mathcal Y=2A$ and $\mathcal Z=23A$ or $23B$. A member of each pair is shown in Figure~\ref{M24dessins}. Unlike in the case of ${\rm M}_{12}$, there is a single permutation representation of degree $p+1$ (the natural representation, on the cosets of a subgroup ${\rm M}_{23}$), and there is no outer automorphism transposing the two mutually inverse conjugacy classes of elements of order $p$. 

We now consider case (c) of Theorem~\ref{th:degp+1}, the affine groups ${\rm AGL}_n(2)$, where $n$ is prime and $p$ is a Mersenne prime $2^n-1$. The case $p=3$ does not arise since ${\rm AGL}_2(2)\cong{\rm S}_4$, so this group is not a quotient of $\Delta(3,2,2)\cong{\rm S}_3$. For $p=7$ the only dessin of type $(3,2,7)$ and degree~$8$ is ${\mathcal D}_0(7)$, with monodromy group 
${\rm PSL}_2(7)$, so ${\rm AGL}_3(2)$ does not arise. (Note that the Hurwitz group of genus $17$ is a {\sl nonsplit \/}extension of ${\rm V}_8$ by ${\rm GL}_3(2)$, not isomorphic to ${\rm AGL}_3(2)$, and it has no transitive permutation representation of degree~$8$; see~\cite{JZ} for details.) On the other hand, for $p=31$ orbit 32.1 in~\cite[pp.~134--135]{APZ} consists of six planar dessins of type $(3,2,31)$ and degree $32$ with monodromy group ${\rm ASL}_5(2)={\rm AGL}_5(2)$. 
These are shown in Section~\ref{sec:App}, Figure~\ref{fig:asl-5-2}.

\begin{theo}\label{AGLMersenne}
For each Mersenne prime $p=2^n-1>7$ there are dessins of type $(3, 2, p)$ and degree $p+1$ with passport $(3^{(p-1)/3}1^2, 2^{(p+1)/2}, p^11^1)$, genus $(p+5)/12$ and monodromy group 
${\rm AGL}_n(2)$. 
\end{theo}

\noindent{\sl Proof}. We will use the fact that $G:={\rm AGL}_n(2)$ is a semidirect product of its translation subgroup $T\cong({\rm C}_2)^n$ by $G_0:={\rm GL}_n(2)$, the stabiliser of $0$. Let $(x, y, z)$ be a generating triple for ${\rm GL}_n(2)$ of type $(3, 2, p)$, chosen as in the proof of Theorem~\ref{th:Mersenne} by taking $r=s=1$, so that $\alpha=1$ and $\beta=(n+1)/2$. Now let $u=xa$, $v=yb$ and $w=zc$, where $a, b, c\in T$. We will try to choose $a, b$ and $c$ so that $(u,v,w)$ is a generating triple for $G$ of type $(3, 2, p)$.

Since $x^3=1$ we have
\[u^3=(xa)^3=x^3a^{x^2}a^xa=a^{x^2+x+1},\]
so $u^3=1$ if and only if $a\in T_x:=\ker(x^2+x+1)$. (Here we regard $x^2+x+1$ as a linear transformation of the vector space $T$ over ${\mathbb F}_2$.) Similarly $v^2=1$ if and only if $b\in T_y:=\ker(y+1)$, and $w^p=1$ if and only if $c\in T_z:=\ker(z^{p-1}+z^{p-2}+\cdots+z+1)$. Also
\[uvw=xa.yb.zc=xya^yzb^zc=xyza^{yz}b^zc=a^{yz}b^zc,\]
so $uvw=1$ if and only if $a^{yz}b^zc=1$.

Since $z$ is a Singer cycle in $G_0$, $T_z=T$ and hence the condition on $c$ is vacuous, whereas the Jordan forms for $x$ and $y$ used in the proof of Theorem~\ref{th:Mersenne} show that the subspaces $T_x$ and $T_y$ have dimensions $n-1$ and $(n+1)/2$. Within these subspaces one can choose $a$ and $b$ arbitrarily, and then define $c=a^{yz}b^z$, so that $a^{vyz}b^zc=1$ and hence $(u, v, w)$ is a triple of type $(3,2,p)$.

The cycle structures of $u$, $v$ and $w$ can be found as follows. One can regard ${\rm AGL}_n(2)$ as acting naturally on $T$, with elements of $T$ and ${\rm GL}_n(2)$ acting by translations and by linear transformations. The element $u=xa$ acts on $T$ by $t\mapsto t^xa$, so it fixes an element $t\in T$ if and only if $a=t(1-x)$, using additive notation for the vector space $T$. Now the Jordan form for $x$ shows that $T(1-x)=T_x$, so for every choice of $a\in T_x$ there are $|\ker(1-x)|=2$ fixed points $t\in T$ for $u$. Since $u$ has order $3$ it has $(p-1)/3$ cycles of length $3$. In the case of $v$ the Jordan form for $y$ shows that $T(1-y)$ is a subspace of codimension $1$ in $T_y$, so by choosing $b\in T_y\setminus T(1-y)$ we can ensure that $v$ has no fixed points, and therefore consists of $(p+1)/2$ transpositions. Since $w$ has order $p$ it has cycle structure $p^11^1$.

Now let us define $H:=\langle u, v, w\rangle\le G$. We will prove that if $p>7$ then $H=G$. Factoring out $T$ maps $H$ onto $\langle x, y, z\rangle=G_0$, so $HT=G$. If $H$ contains a non-identity translation $t\in T$ then since $H$ contains $w$ it contains all the non-identity translations $t^w$, so $H\ge T$ and hence $H=G$, as required. 

We may therefore assume, for a contradiction, that $H\cap T=1$, so that $H$ is a complement for $T$ in $G$. Since the first cohomology group $H^1({\rm GL}_n(2),T)$ for ${\rm GL}_n(2)$ on its natural module $T$ is trivial for all $n>3$ (see~\cite{Bell}) there is a single conjugacy class of such complements, so $H$ is a point-stabiliser in the natural action of $G$. This contradicts the fact that $v$ has no fixed points, by our choice of the translation $b$. 

Thus $(u, v, w)$ is a generating triple for $G$, so it corresponds to a dessin of type $(3,2,p)$ and degree~$p+1$ with monodromy group $G$. Its genus follows from its passport, determined above. \hfill$\square$

\medskip

The above proof fails in the case $p=7$ since $|H^1({\rm GL}_3(2),T)|=2$ (see~\cite{Bell}), so that there are two conjugacy classes of complements for $T$ in ${\rm AGL}_3(2)$, consisting of the point-stabilisers and of subgroups acting transitively on $T$ as ${\rm PSL}_2(7)$; in this case every triple of type $(3,2,7)$ generates a complement. (The fact that all complements are isomorphic to $G/T$ provides a simple proof that ${\rm GL}_3(2)\cong {\rm PSL}_2(7)$.)

\medskip

\begin{exam}\label{ex:deg32}\rm
Taking $n=5$ we obtain a dessin of type $(3,2,31)$, degree $32$ and genus $3$, with passport $(3^{10}1^2, 2^{16}, 31^11^1)$. This should not be confused with the Galois orbit 32.1 of six planar dessins with the same type, degree and monodromy group in~\cite{APZ}: these have passport $(3^{10}1^2, 2^{12}1^8, 31^11^1)$, and they correspond to taking $b\in T(1-y)$, so that $v$ fixes $|\ker(1-y)|=8$ points (see Subsection~\ref{ex:AGL_5_2_another}).
\end{exam}

Finally, we consider the alternating groups ${\rm A}_{p+1}$, in case~(a) of Theorem~\ref{th:degp+1}, as monodromy groups.

\begin{theo}
The alternating group ${\rm A}_{p+1}$, $p$ prime, is the monodromy group of a dessin of type $(3,2,p)$ and degree $p+1$ if and only if $p=3$ or $p\ge 11$.
\end{theo}

\noindent{\sl Proof.} For $p=2, 5$ or $7$ the triangle group $\Delta(3,2,p)$ does not map onto ${\rm A}_{p+1}$: this is obvious if $p=2$ or $5$, and for $p=7$ it is well known that ${\rm A}_8$ is not a Hurwitz group: the triples of type $(3,2,7)$ in this group all generate proper subgroups isomorphic to ${\rm PSL}_2(7)$. On the other hand, ${\rm A}_4$ is the monodromy group of the dessin of type $(3,2,3)$ and degree $4$ on the left in Figure~\ref{A4,A12,A14dessins}, the quotient of the tetrahedron $\{3,3\}$ by ${\rm C}_3$. (Note that ${\rm A}_4$ is isomorphic to ${\rm PSL}_2(3)$ in case (b).)

\begin{figure}[!ht]

\begin{center}
 \begin{tikzpicture}[scale=0.15, inner sep=0.8mm]
 
\node (a1) at (-30,0) [shape=circle, fill=black] {};
\node (b1) at (-25,0) [shape=circle, fill=black] {};
 \draw [thick] (-30,0) arc (0:360:3);
\draw [thick] (a1) to (b1);


\node (a) at (-10,0) [shape=circle, fill=black] {};
\node (b) at (-5,0) [shape=circle, fill=black] {};
\node (c) at (0,0) [shape=circle, fill=black] {};
\node (d) at (5,0) [shape=circle, fill=black] {};

 \draw [thick] (-10,0) arc (0:360:3);
 \draw [thick] (a) to (d);
 \draw [thick] (b) to (-5,-5);
 \draw [thick] (c) to (0,5);
 \draw [thick] (8,3) to (d) to (8,-3);
 
 
\node (A) at (20,0) [shape=circle, fill=black] {};
\node (B) at (25,0) [shape=circle, fill=black] {};
\node (C) at (30,0) [shape=circle, fill=black] {};
\node (D) at (35,0) [shape=circle, fill=black] {};
\node (E) at (38,3) [shape=circle, fill=black] {};
\node (F) at (38,-3) [shape=circle, fill=black] {};

 \draw [thick] (20,0) arc (0:360:3);
 \draw [thick] (A) to (D);
 \draw [thick] (B) to (25,-5);
 \draw [thick] (C) to (30,5);
 \draw [thick] (E) to (D) to (F);

\end{tikzpicture}

\end{center}
\caption{Dessins with monodromy groups ${\rm A}_4$, ${\rm A}_{12}$ and ${\rm A}_{14}$.}
\label{A4,A12,A14dessins}
\end{figure}
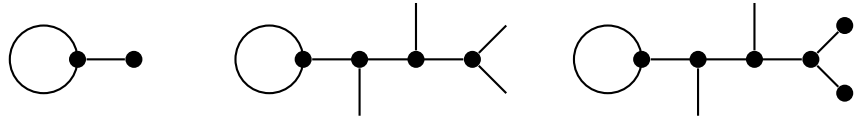

For $p=11$ and $13$ it is sufficient  by Theorem~\ref{th:degp+1} to note that the dessins of degree $12$ and $14$ in Figure~\ref{A4,A12,A14dessins} are not isomorphic to the unique dessins ${\mathcal D}_0(p)$ with monodromy group ${\rm PSL}_2(p)$ shown in Figures~\ref{deg12map} and \ref{D0(p)}.

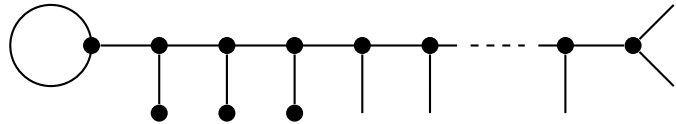
\begin{figure}[!ht]

\begin{center}
 \begin{tikzpicture}[scale=0.18, inner sep=0.8mm]

\node (a) at (-10,0) [shape=circle, fill=black] {};
\node (b) at (-5,0) [shape=circle, fill=black] {};
\node (c) at (0,0) [shape=circle, fill=black] {};
\node (d) at (5,0) [shape=circle, fill=black] {};
\node (e) at (10,0) [shape=circle, fill=black] {};
\node (f) at (15,0) [shape=circle, fill=black] {};
\node (g) at (25,0) [shape=circle, fill=black] {};
\node (h) at (30,0) [shape=circle, fill=black] {};
\node (B) at (-5,-5) [shape=circle, fill=black] {};
\node (C) at (0,-5) [shape=circle, fill=black] {};
\node (D) at (5,-5) [shape=circle, fill=black] {};

 \draw [thick] (-10,0) arc (0:360:3);
 \draw [thick] (a) to (17,0);
 \draw [thick, dashed] (18,0) to (22,0);
  \draw [thick] (23,0) to (h);
 \draw [thick] (b) to (B);
 \draw [thick] (c) to (C);
 \draw [thick] (d) to (D);
 \draw [thick] (e) to (10,-5);
 \draw [thick] (f) to (15,-5);
 \draw [thick] (g) to (25,-5);
 \draw [thick] (33,3) to (h) to (33,-3);

\end{tikzpicture}

\end{center}
\caption{A dessin with monodromy group ${\rm A}_{p+1}$.}
\label{Ap+1dessin}
\end{figure}

Finally, let $p\ge 17$. If $p\equiv 2$ mod~$(3)$, say $p=3k+2$ for some $k\ge 5$, consider the dessin $\mathcal D$ in Figure~\ref{Ap+1dessin}, where there are $k$ vertices of degree $3$, so that the degree of $\mathcal D$ is $p+1$ and its passport is $(3^k1^3, 2^{k+3}1^{k-3}, p^11^1)$. The monodromy group of $\mathcal D$ cannot be ${\rm PSL}_2(p)$ since a non-identity element of that group has at most two fixed points, whereas $x$ has three; similarly it cannot be ${\rm AGL}_n(2)$ for any $n$ since the fixed point set of any element of that group is either empty or an affine subspace of order $2^d$ for some $d$; finally, it cannot be a Mathieu group since $\mathcal D$ is not isomorphic to any of the dessins in Figures~\ref{M12dessins} or~\ref{M24dessins} or their mirror images, so it must be ${\rm A}_{p+1}$.

If $p\equiv 1$ mod~$(3)$, say $p=3k+4$ for some $k\ge5$, one can apply a similar argument to the dessin ${\mathcal D}'$ obtained by adding two vertices of valency $1$ to the rightmost free edges of $\mathcal D$, so that the passport is $(3^k1^5, 2^{k+5}1^{k-5}, p^11^1)$.\hfill$\square$

\medskip

By reflecting an arbitrary subset of the free edges of $\mathcal D$ or ${\mathcal D}'$ one obtains exponentially many dessins of type $(3,2,p)$ and degree $p+1$ with monodromy group ${\rm A}_{p+1}$ as $p\to\infty$.

We summarise the results of the last two sections as follows:

\begin{theo}\label{th:summary}
Suppose that $p$ is a prime such that there is a dessin of type $(3,2,p)$ and degree $p+1$ with monodromy group $G$. Then one of the following holds:
\begin{itemize}
\item[\rm (a)] $G={\rm A}_{p+1}$, with $p=3$ or $p\ge 11$,
\item[\rm (b)] $G={\rm PSL}_2(p)$,
\item[\rm (c)] $G={\rm AGL}_n(2)$ for some Mersenne prime $p=2^n-1 \ge 31$,
\item[\rm (d)] $G={\rm M}_{12}$ with $p=11$ or $G={\rm M}_{24}$ with $p=23$,
\end{itemize}
each in its natural representation.
In case~(a) there are exponentially many dessins as $p\to\infty$, even if we restrict to planar dessins. In case~(b) there is a single dessin ${\mathcal D}_0(p)$ for each prime $p$. In case~(c) there is at least one dessin for each Mersenne prime $p$. In case~(d) there are two chiral pairs for each of the two groups $G$.
\end{theo} 


\section{The road not taken}
\label{sec:3dim}

The starting point for Klein's work in~\cite{Kle79}, and hence for our investigation in this paper, was the embedding of the icosahedral group ${\rm A}_5$ in ${\rm PSL}_2(11)$. This is the third and most complicated of the three instances, known already to Galois, of ${\rm PSL}_2(p)$ having a subgroup of index $p$; the first two, for $p=5$ and $7$, each became the subject of a deservedly famous book (respectively~\cite{Kle84} and \cite{Lev}), and one could easily imagine (though not so easily write) an analogue for $p=11$. In this final section we will briefly sketch a line of research which might form the basis of a chapter in such a text; the road it takes is almost completely disjoint from that we have followed  in the present paper; however, there is an intriguing moment when we catch a brief glimpse of the road starting with Klein's paper~\cite{Kle78} concerning the case $p=7$.

There are a number of mathematical objects and phenomena closely related to the embedding (or more precisely the embeddings) of ${\rm A}_5$ in the group $L={\rm PSL}_2(11)$. A good example is the non-orientable regular abstract polytope ${}_5\{3,5,3\}_5$, the hendecachoron or $11$-cell independently discovered by Gr\"unbaum~\cite{Gru} and Coxeter~\cite{Cox84} (see also~\cite{CE, CW} for this example, and~\cite{McMS} for the general theory of abstract polytopes). 

\medskip

\begin{figure}[!ht]

\begin{center}
 \begin{tikzpicture}[scale=0.15, inner sep=0.8mm]

 \node (A) at (-15,0) [shape=circle, draw] {};
 \node (B) at (-5,0) [shape=circle, draw] {};
 \node (C) at (5,0) [shape=circle, draw] {};
 \node (D) at (15,0) [shape=circle, draw] {};

 \draw [thick] (A) to (B) to (C) to (D); 

\node at (0,-2) {$5$};

\end{tikzpicture}

\vspace{-5mm}

\end{center}
\caption{The Coxeter diagram $\{3,5,3\}$.} 
\label{353}
\end{figure}
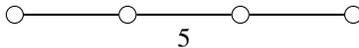

Let $\Gamma$ denote the string Coxeter group with Schl\"afli symbol $\{3,5,3\}$, represented by the Coxeter diagram in Figure~\ref{353}. This group has generators $R_i\;(i=0,\ldots, 3)$ and defining relations
\[R_i^2=(R_0R_1)^3=(R_1R_2)^5=(R_2R_3)^3=1,\quad (R_iR_j)^2=1\;\hbox{whenever}\;|i-j|\ge 2;\]
it acts as a group of isometries of hyperbolic $3$-space ${\mathbb H}^3$, generated by reflections $R_i$ in the sides of a tetrahedron $T\subset{\mathbb H}^3$ with appropriate dihedral angles. As shown in~\cite{JLM} there is a unique normal subgroup $N$ of $\Gamma$ such that $\Gamma/N\cong L$. This is the normal closure of the elements $(R_0R_1R_2)^5$ and $(R_1R_2R_3)^5$, so that putting these equal to $1$ in $\Gamma$ gives a presentation for $L$.

The quotient ${\mathbb H}^3/N$ is a compact non-orientable hyperbolic $3$-orbifold $\mathcal N$, with a tessellation $\mathcal T$ by $|L|=660$ tetrahedra, induced from the tessellation of ${\mathbb H}^3$ by the images of $T$ under $\Gamma$. These tetrahedra can be partitioned into eleven sets of $60$, giving a tessellation $\mathcal I$ of $\mathcal N$ by eleven hemi-icosahedra $\{3,5\}_5$ (antipodal quotients of icosahedra $\{3,5\}$ formed by identifying edges five steps apart along Petrie paths), which has $\mathcal T$ as its barycentric subdivision. Similarly $\mathcal I$ has eleven vertices, each having as its vertex figure a hemi-dodecahedron $\{5,3\}_5$ formed in the same way from a dodecahedron $\{5,3\}$. This tessellation $\mathcal I$ is a realisation of Coxeter's non-orientable regular polytope ${}_5\{3,5,3\}_5$, with the subscripts indicating these identifications, and also the corresponding extra defining relations for $L$. It has automorphism group $L$, with the stabilisers of cells and of vertices forming the two conjugacy classes of subgroups $H\cong {\rm A}_5$ in $L$, represented by the images of the subgroups $\langle R_0, R_1, R_2\rangle\cong\langle R_1, R_2, R_3\rangle\cong\Delta[3,2,5]\cong {\rm A}_5\times {\rm C}_2$ of $\Gamma$. Notice that $N$ is not torsion-free (it contains the central involution of each of these subgroups), so $\mathcal N$ is not a manifold.

\begin{figure}[!ht]

\begin{center}
 \begin{tikzpicture}[scale=0.15, inner sep=0.8mm]

 \node (A) at (-10,20) [shape=circle, fill=black] {};
 \node (B) at (-15,15) [shape=circle, fill=black] {};
 \node (C) at (-5,15) [shape=circle, fill=black] {};
 \node (D) at (-10,10) [shape=circle, fill=black] {};
 \node (E) at (-10,-5) [shape=circle, fill=black] {};
 \node (F) at (-10,-20) [shape=circle, fill=black] {};
 \node (G) at (-15,0) [shape=circle, fill=black] {};

 \draw [thick] (A) to (B) to (D) to (C) to (A); 
 \draw [thick] (D) to (E) to (F); 
 \draw [thick] (B) to (G) to (E); 
 
 \node at (-13,20) {$\Omega$};
 \node at (-18,15) {$\Gamma$};
 \node at (-2,15) {$\Omega^+$};
  \node at (-7,10) {$\Gamma^+$};
 \node at (-18,0) {$N$};
 \node at (-13,-5) {$K$};
 \node at (-13,-20) {$1$};
 
 
  \node (a) at (20,20) [shape=circle, fill=black] {};
 \node (b) at (25,15) [shape=circle, fill=black] {};
 \node (c) at (15,15) [shape=circle, fill=black] {};
 \node (d) at (20,10) [shape=circle, fill=black] {};
 \node (e) at (20,-5) [shape=circle, fill=black] {};
 \node (g) at (15,0) [shape=circle, fill=black] {};

 \draw [thick] (a) to (b) to (d) to (c) to (a); 
 \draw [thick] (d) to (e); 
 \draw [thick] (c) to (g) to (e); 
 
 \node at (31,20) {${\rm PGL}_2(11)\times {\rm C}_2$};
  \node at (10,15) {$L\times {\rm C}_2$};
 \node at (32,15) {${\rm PGL}_2(11)$};
 \node at (23,10) {$L$};
 \node at (12,0) {${\rm C}_2$};
 \node at (23,-5) {$1$};

\end{tikzpicture}

\end{center}
\caption{Subgroups of $\Omega$ and their quotients by $K$.}
\label{Omega}
\end{figure}
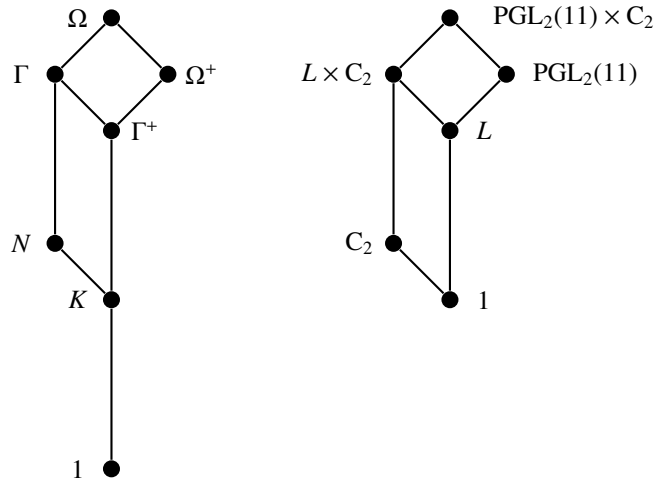

The normaliser $\Omega$ of $\Gamma$ in the isometry group of ${\mathbb H}^3$ is a semidirect product of $\Gamma$ by a group ${\rm C}_2$ induced by the graph automorphism $R_i\mapsto R_{3-i}$ of $\Gamma$; this leaves $N$ invariant, with $\Omega/N\cong {\rm PGL}_2(11)$, and acts as a self-duality of $\mathcal I$, transposing vertices and cells, and edges and faces. The group $\Omega$ plays a similar role in the theory of compact $3$-manifolds to that played by the triangle groups $\Delta(3,2,7)$ and $\Delta[3,2,7]$ for Riemann and Klein surfaces: specifically, as shown by Gehring, Marshall and Martin~\cite{GM,MM}, $\Omega$ has the least covolume among all cocompact discontinuous groups of isometries of ${\mathbb H}^3$. It follows that the torsion-free normal subgroups $K$ of finite index in $\Omega$ uniformise the compact hyperbolic $3$-manifolds which maximise the number of automorphisms, namely $|\Omega/K|$, per unit of volume. As such, these manifolds and automorphism groups are the $3$-dimensional analogues of the Hurwitz curves and groups in dimension~$2$. As shown in~\cite{JLM}, among such subgroups $K$ that of least index is $K=N\cap \Omega^+$, where $\Omega^+$ is the orientation-preserving subgroup of index $2$ in $\Omega$. The corresponding orientable manifold ${\mathcal K}={\mathbb H}^3/K$, a double covering of the orbifold $\mathcal N$, has isometry group $\Omega/K\cong {\rm PGL}_2(11)\times {\rm C}_2$, with orientation-preserving subgroup $\Omega^+/K\cong {\rm PGL}_2(11)$; the subgroup $\Gamma/K\cong L\times {\rm C}_2$ preserves a tessellation of $\mathcal K$ by eleven icosahedra, with the two conjugacy classes of icosahedral subgroups ${\rm A}_5\times {\rm C}_2$ stabilising the cells and vertices, and transposed in $\Omega/K$ by duality. By the minimality of their volumes, $\mathcal K$ and $\mathcal N$ can be regarded as $3$-dimensional analogues of Klein's quartic curve~\cite{Kle78}, the smallest Riemann surface attaining the Hurwitz bound on the number of automorphisms per unit of area (see~\cite{Lev} for connections with ${\rm PSL}_2(7)$ and~\cite{Kle78}). Figure~\ref{Omega} shows the relevant subgroups of $\Omega$, together with their quotients by $K$, that is, the finite groups of isometries they induce on $\mathcal K$.


\input{more-examples}


\bigskip

\noindent{\bf Acknowledgements} We are very grateful to Jean B\'etr\'ema, Yuri Bilu, 
John Voight and Dimitri Zvonkine for valuable assistance and suggestions. 
Alexander Zvonkin was partially supported by the ANR project {\sc Combin\'e}
(ANR-19-CE48-0011).



\bigskip

{\small
\noindent
{\sc School of Mathematical Sciences, University of Southampton, Southampton} SO17 1BJ, UK 

\noindent
{\em E-mail address}: \url{G.A.Jones@maths.soton.ac.uk}

\bigskip

\noindent
{\sc LaBRI, Universit\'e de Bordeaux, 351 cours de la Lib\'eration, F-33405 Talence Cedex, 
France} 

\noindent
{\em E-mail address}: \url{zvonkin@labri.fr}
}

\end{document}

%% file: more-examples.tex

\section{More examples}
\label{sec:App}

In this section we give standard generators for the monodromy groups of some of the smaller and more interesting dessins considered earlier, together with some comments, and diagrams in a few cases of low genus.

\subsection{Klein's dessins, degree 11} 
\label{sec:kl-11}

Generators $x_i, y_i$ and $z_i$ for the monodromy groups of Klein's dessins 
$\mathcal M_i$ ($i=1,\ldots,6)$ (see Section~\ref{sec:klein-trees} and 
Figure~\ref{Kleinmaps}), as subgroups of ${\rm S}_{11}$, are as follows. In all cases
\[z_i = (1,2,\ldots, 11).\]
The other generators $x_i$ and $y_i$ for $i=1,\ldots,6$ (unique up to conjugation 
by powers of~$z_i$) are:

\smallskip

\begin{center}
\begin{tabular}{lcl}
$x_1=(1,4,3)(5,11,9)(6,8,7)$,  && $y_1=(1,2)(4,11)(5,8)(9,10)$; \\
$x_2=(1,4,2)(5,11,9)(6,8,7)$,  && $y_2=(2,3)(4,11)(5,8)(9,10)$; \\
$x_3=(1,3,2)(4,11,10)(5,9,7)$, && $y_3=(3,11)(4,9)(5,6)(7,8)$;  \\
$x_4=(1,4,3)(5,11,10)(6,9,8)$, && $y_4=(1,2)(4,11)(5,9)(6,7)$;  \\
$x_5=(1,11,6)(2,5,4)(7,10,8)$, && $y_5=(1,5)(2,3)(6,10)(8,9)$;  \\
$x_6=(1,11,6)(2,5,3)(7,10,9)$, && $y_6=(1,5)(3,4)(6,10)(7,8)$.
\end{tabular}
\end{center}

\begin{rema}\label{rem:mirror-image}\rm
Here and everywhere, if two dessins $\mathcal{D}$ and $\overline{\mathcal{D}}$ of degree $n$ 
form a chiral pair and $\mathcal{D}=(x,y,z)$ then $\overline{\mathcal{D}}$ can be obtained
as $(x',y',z')=(x^{-1},y^{-1},yx)$. Notice that $z'=yx$ does not necessarily have any 
canonical form. It can, however, be ``standardized'' by a common conjugation of $x',y',z'$
in ${\rm S}_{n}$.
\end{rema}

In what follows, in order to avoid any error we just reproduce the results of the 
GAP sessions.

\subsection{Group ${\rm PSL}_3(3)$, degree 13}

There are, in total, 14 plane trees with the passport $(3^4 1^1,2^4 1^5,13^1)$.
Four of them have the monodromy group ${\rm PSL}_3(3)$. They are given below.
We fix the permutation $z = (1,2,\ldots,13$), and label the only fixed 
point of permutation $x$ by 1. The trees split into two chiral pairs: $(\M_1,\M_3)$ 
and $(\M_2,\M_4)$. Trees $\M_1$ and $\M_2$ are shown in Figure \ref{PSL_3(3)maps}
(see Example~\ref{ex:p=13}).

\medskip

{\footnotesize
\begin{verbatim}
M1 = [ (2,13,6)(3,5,4)(7,12,11)(8,10,9), (1,13)(2,5)(6,12)(7,10),
       (1,2,3,4,5,6,7,8,9,10,11,12,13) ],
       
M2 = [ (2,13,3)(4,12,5)(6,11,10)(7,9,8), (1,13)(3,12)(5,11)(6,9),
       (1,2,3,4,5,6,7,8,9,10,11,12,13) ],
       
M3 = [ (2,13,9)(3,8,4)(5,7,6)(10,12,11), (1,13)(2,8)(4,7)(9,12),
       (1,2,3,4,5,6,7,8,9,10,11,12,13) ],
       
M4 = [ (2,13,12)(3,11,10)(4,9,5)(6,8,7), (1,13)(2,11)(3,9)(5,8),
       (1,2,3,4,5,6,7,8,9,10,11,12,13) ]
\end{verbatim}
}

\medskip

\noindent
The field of moduli of these four trees is the splitting field of the polynomial
$a^4+a^3+2a^2-4a+3$. Its Galois group is the cyclic group ${\rm C}_4$.


\subsection{Group ${\rm PSL}_2(16)$, degree 17}

\begin{figure}[htbp]
\begin{center}
\epsfig{file=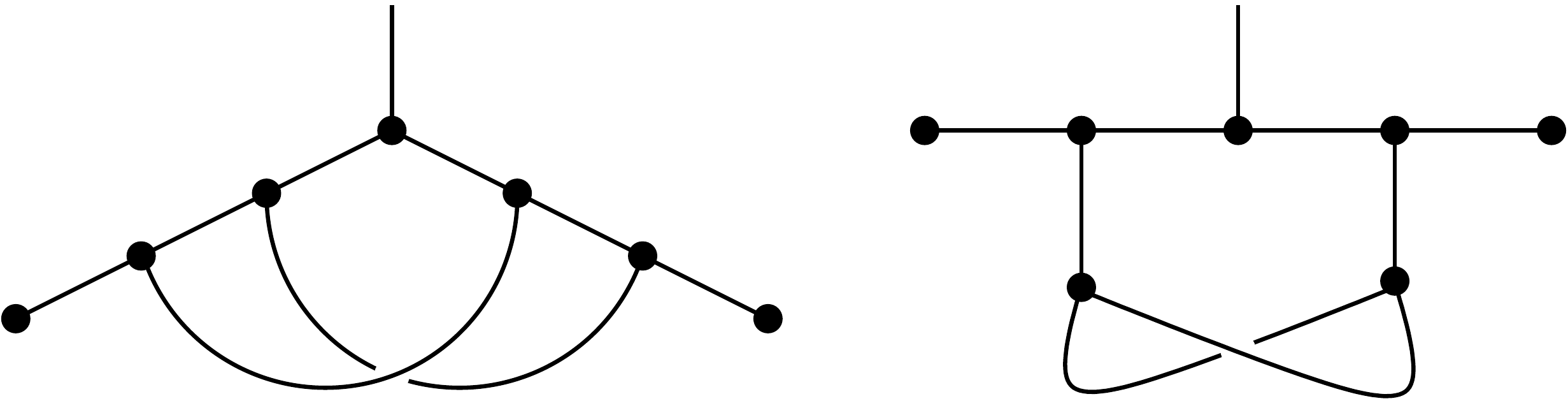,width=10cm}
\end{center}
\caption{Two maps representing the group ${\rm PSL}_2(16)$.}
\label{fig:psl-2-16}
\end{figure}

There are, in this group, eight conjugacy classes of elements of order 17, and 
also four automorphisms of the field $\mathbb{F}_{16}$. Therefore, as predicted by
Theorem~\ref{Fermat}, we obtain two non-isomorphic maps (see Example~\ref{ex:small-proj}).
Their passport is $(3^5 1^2,2^8 1^1,17^1)$, hence they are of genus 1.
We also observe that they are not chiral: both are mirror symmetric. They are shown
in Figure \ref{fig:psl-2-16}.

We fix, as usual, the permutation $z = (1,2,\ldots,17)$, and put the label 1 on the 
only fixed point of the permutation $y$. Here is what we get:

\medskip

{\footnotesize
\begin{verbatim}
M1 = [ (1,10,2)(3,9,13)(4,12,6)(7,11,17)(8,16,14),
       (2,9)(3,12)(4,5)(6,11)(7,16)(8,13)(10,17)(14,15),
       (1,2,3,4,5,6,7,8,9,10,11,12,13,14,15,16,17) ],
       
M2 = [ (1,10,2)(3,9,5)(6,8,13)(7,12,14)(11,17,15),
       (2,9)(3,4)(5,8)(6,12)(7,13)(10,17)(11,14)(15,16),
       (1,2,3,4,5,6,7,8,9,10,11,12,13,14,15,16,17) ]
\end{verbatim}
}

\medskip

The total number of maps with this passport is 70.

\begin{rema}\rm
The next Fermat prime, after 17, is 257. Let us give a cursory glance at the group 
${\rm PSL}_2(256)$. In this group, there are 128 conjugacy classes of cycles of 
order 257, and there are eight automorphisms of the field $\mathbb{F}_{256}$. 
Therefore, there are 16 non-isomorphic maps generating this group. Their passport is
$(3^{85}1^2,2^{128}1^1,257^1)$, and their genus is $g=21$, in accordance with
Theorem~\ref{Fermat}. This theorem also implies that all these maps are mirror symmetric.
\end{rema}

\subsection{Group ${\rm PSL}_3(5)$, degree 31}
\label{sec:deg-31a}

\begin{figure}[!ht]

\begin{center}
 \begin{tikzpicture}[scale=0.16, inner sep=0.8mm]
 
\node (a) at (-15,5) [shape=circle, fill=black] {};
\node (b) at (-10,5) [shape=circle, fill=black] {};
\node (c) at (-5,10) [shape=circle, fill=black] {};
\node (d) at (-5,5) [shape=circle, fill=black] {};
\node (e) at (-5,0) [shape=circle, fill=black] {};
\node (f) at (-5,-8) [shape=circle, fill=black] {};
\node (g) at (0,5) [shape=circle, fill=black] {};
\node (h) at (0,0) [shape=circle, fill=black] {};
\node (i) at (5,0) [shape=circle, fill=black] {};
\node (j) at (10,0) [shape=circle, fill=black] {};
\node (k) at (15,0) [shape=circle, fill=black] {};

\draw [thick] (a) to (d);
\draw [thick] (c) to (f);
\draw [thick] (e) to (k);
\draw [thick] (g) to (h); 

\draw [thick] (a) to (-15,0);
\draw [thick] (b) to (-10,10);
\draw [thick] (j) to (10,-5);

\draw [thick] (-8,13) to (c) to (-2,13);
\draw [thick] (18,3) to (k) to (18,-3);

\draw [thick] (a) to (-19,5);
\draw [thick] (i) to (5,13.5);
\draw [thick] (f) to (17,-8);
\draw [thick] (f) to (-5,-13.5);

\draw [dashed] (-10,20) to (-25,-5) to (10,-20) to (25,5) to (-10,20);
 
\end{tikzpicture}

\end{center}
\caption{Opposite sides of the outer parallelogram are identified to form a torus. 
The torus dessin thus obtained is of degree $31$ and type $(3,2,31)$ with monodromy 
group ${\rm PSL}_3(5)$. It corresponds to the triple $\M_1=(x_1,y_1,z)$ in the above list.}
\label{PSL_3(5)dessin}
\end{figure}
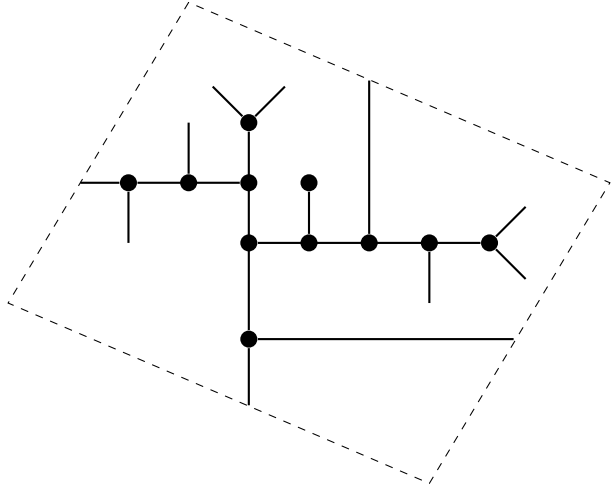

Below we give the ten dessins of type $(3,2,31)$ with the passport 
$(3^{10}1^1,2^{12}1^7,31^1)$, of genus $g=1$, representing the group ${\rm PSL}_3(5)$
(see Example~\ref{ex:psl-3-5}). 
The map $\M_1$ is shown in Figure \ref{PSL_3(5)dessin}.

\begin{itemize}
\item	The permutation $z$ is always the same: $z=(1,2,3,\ldots,31)$. 
\item	The fixed point of $x$ is given the label 1. This convention, together
		with the preceding one, allows us to label all the dessins uniquely
\item	The chiral pairs are $(\mathcal{M}_i,\mathcal{M}_{i+5})$ for $i=1,\ldots,5$.
\end{itemize}

\medskip

{\footnotesize
\begin{verbatim}
M1 = [ (2,31,19)(3,18,11)(4,10,21)(5,20,30)(6,29,25)(7,24,23)(8,22,9)(12,17,16)
       (13,15,14)(26,28,27), (1,31)(2,18)(3,10)(4,20)(5,29)(6,24)(7,22)(9,21)
       (11,17)(12,15)(19,30)(25,28), 
       (1,2,3,4,5,6,7,8,9,10,11,12,13,14,15,16,17,18,19,20,21,22,23,24,25,26,
        27,28,29,30,31) ]
        
M2 = [ (2,31,3)(4,30,16)(5,15,6)(7,14,8)(9,13,21)(10,20,23)(11,22,12)(17,29,18)
       (19,28,24)(25,27,26), (1,31)(3,30)(4,15)(6,14)(8,13)(9,20)(10,22)(12,21)
       (16,29)(18,28)(19,23)(24,27), 
       (1,2,3,4,5,6,7,8,9,10,11,12,13,14,15,16,17,18,19,20,21,22,23,24,25,26,
        27,28,29,30,31) ], 
        
M3 = [ (2,31,3)(4,30,29)(5,28,19)(6,18,13)(7,12,8)(9,11,10)(14,17,21)(15,20,27)
       (16,26,22)(23,25,24), (1,31)(3,30)(4,28)(5,18)(6,12)(8,11)(13,17)(14,20)
       (15,26)(16,21)(19,27)(22,25), 
       (1,2,3,4,5,6,7,8,9,10,11,12,13,14,15,16,17,18,19,20,21,22,23,24,25,26,
        27,28,29,30,31) ], 
        
M4 = [ (2,31,24)(3,23,16)(4,15,14)(5,13,12)(6,11,7)(8,10,9)(17,22,21)(18,20,27)
       (19,26,28)(25,30,29), (1,31)(2,23)(3,15)(4,13)(5,11)(7,10)(16,22)(17,20)
       (18,26)(19,27)(24,30)(25,28), 
       (1,2,3,4,5,6,7,8,9,10,11,12,13,14,15,16,17,18,19,20,21,22,23,24,25,26,
        27,28,29,30,31) ], 
        
M5 = [ (2,31,3)(4,30,8)(5,7,6)(9,29,15)(10,14,19)(11,18,24)(12,23,22)(13,21,20)
       (16,28,27)(17,26,25), (1,31)(3,30)(4,7)(8,29)(9,14)(10,18)(11,23)(12,21)
       (13,19)(15,28)(16,26)(17,24), 
       (1,2,3,4,5,6,7,8,9,10,11,12,13,14,15,16,17,18,19,20,21,22,23,24,25,26,
        27,28,29,30,31) ], 
        
M6 = [ (2,31,14)(3,13,28)(4,27,8)(5,7,6)(9,26,10)(11,25,24)(12,23,29)(15,30,22)
       (16,21,17)(18,20,19), (1,31)(2,13)(3,27)(4,7)(8,26)(10,25)(11,23)(12,28)
       (14,30)(15,21)(17,20)(22,29), 
       (1,2,3,4,5,6,7,8,9,10,11,12,13,14,15,16,17,18,19,20,21,22,23,24,25,26,
        27,28,29,30,31) ], 
        
M7 =  [ (2,31,30)(3,29,17)(4,16,15)(5,14,9)(6,8,7)(10,13,23)(11,22,21)(12,20,24)
        (18,28,27)(19,26,25), (1,31)(2,29)(3,16)(4,14)(5,8)(9,13)(10,22)(11,20)
        (12,23)(17,28)(18,26)(19,24), 
        (1,2,3,4,5,6,7,8,9,10,11,12,13,14,15,16,17,18,19,20,21,22,23,24,25,26,
         27,28,29,30,31) ], 
        
M8 = [ (2,31,30)(3,29,4)(5,28,14)(6,13,18)(7,17,11)(8,10,9)(12,16,19)(15,27,20)
       (21,26,25)(22,24,23), (1,31)(2,29)(4,28)(5,13)(6,17)(7,10)(11,16)(12,18)
       (14,27)(15,19)(20,26)(21,24), 
       (1,2,3,4,5,6,7,8,9,10,11,12,13,14,15,16,17,18,19,20,21,22,23,24,25,26,
        27,28,29,30,31) ], 
        
M9 = [ (2,31,9)(3,8,4)(5,7,14)(6,13,15)(10,30,17)(11,16,12)(18,29,19)(20,28,21)
       (22,27,26)(23,25,24), (1,31)(2,8)(4,7)(5,13)(6,14)(9,30)(10,16)(12,15)
       (17,29)(19,28)(21,27)(22,25), 
       (1,2,3,4,5,6,7,8,9,10,11,12,13,14,15,16,17,18,19,20,21,22,23,24,25,26,
        27,28,29,30,31) ], 
        
M10 = [ (2,31,30)(3,29,25)(4,24,18)(5,17,6)(7,16,8)(9,15,22)(10,21,11)(12,20,13)
        (14,19,23)(26,28,27), (1,31)(2,29)(3,24)(4,17)(6,16)(8,15)(9,21)(11,20)
        (13,19)(14,22)(18,23)(25,28), 
        (1,2,3,4,5,6,7,8,9,10,11,12,13,14,15,16,17,18,19,20,21,22,23,24,25,26,
         27,28,29,30,31) ] 
\end{verbatim}
}

\medskip

\subsection{Group ${\rm PSL}_5(2)$, degree 31}
\label{sec:deg-31b}
 
\begin{figure}[htbp]
\begin{center}
\epsfig{file=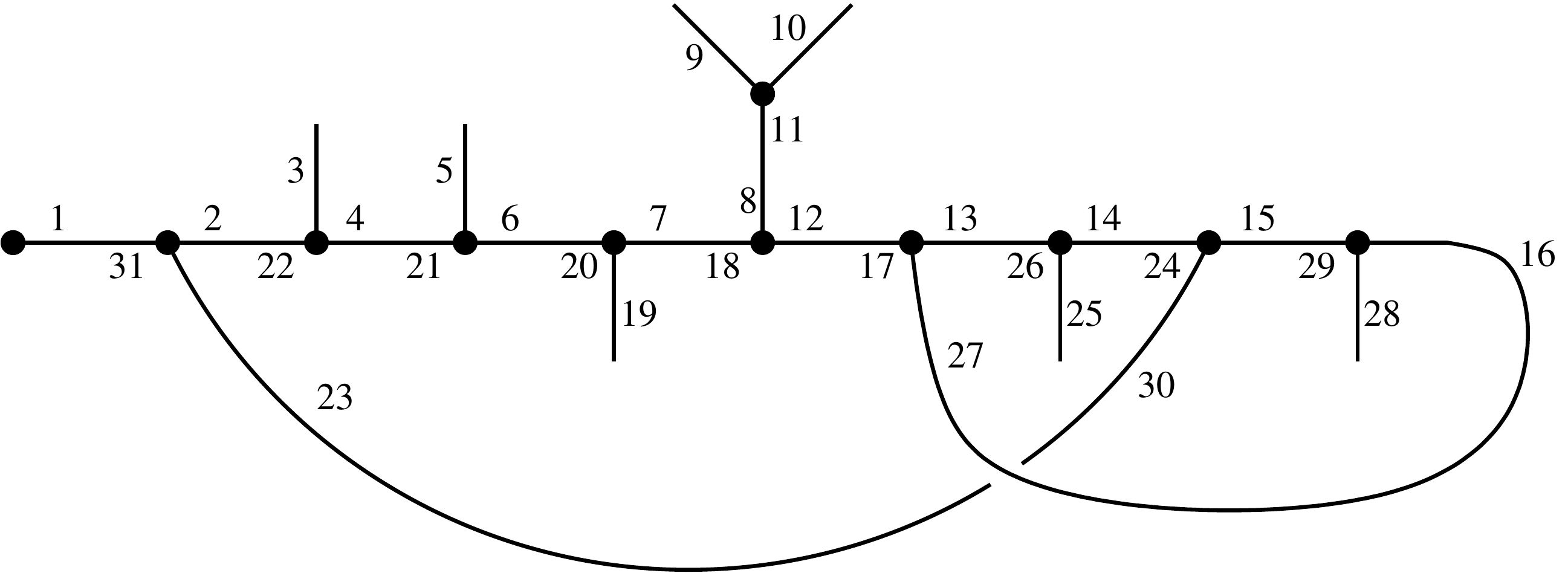,width=10cm}
\end{center}
\caption{The map $\M_1$ of the above list representing the group ${\rm PSL}_5(2)$. 
The map $\M_4$ is mirror symmetric to this one.}
\label{fig:psl-5-2}
\end{figure}

Below we give the six dessins with the same passport as in the previous example, 
that is, $(3^{10}1^1,2^{12}1^7,31^1)$, and therefore also of genus 1, but this time 
representing the group ${\rm PSL}_5(2)$ (see Example~\ref{ex:p=31Mersenne}).
The map $\M_1$ is shown in Figure~\ref{fig:psl-5-2}.

\begin{itemize}
\item	The permutation $z$ is always the same: $z=(1,2,3,\ldots,31)$. 
\item	The fixed point of $x$ is always the point $1$.
\item	The chiral pairs are $(\M_1,\M_4)$, $(\M_2,\M_5)$, and $(\M_3,\M_6)$.
\end{itemize}

\medskip

{\footnotesize
\begin{verbatim}
M1 = [ (2,31,23)(3,22,4)(5,21,6)(7,20,19)(8,18,12)(9,11,10)(13,17,27)(14,26,25)
       (15,24,30)(16,29,28), (1,31)(2,22)(4,21)(6,20)(7,18)(8,11)(12,17)(13,26)
       (14,24)(15,29)(16,27)(23,30),
       (1,2,3,4,5,6,7,8,9,10,11,12,13,14,15,16,17,18,19,20,21,22,23,24,25,26,
        27,28,29,30,31) ], 
        
M2 = [ (2,31,11)(3,10,25)(4,24,20)(5,19,18)(6,17,16)(7,15,27)(8,26,9)(12,30,13)
       (14,29,28)(21,23,22), (1,31)(2,10)(3,24)(4,19)(5,17)(6,15)(7,26)(9,25)
       (11,30)(13,29)(14,27)(20,23), 
       (1,2,3,4,5,6,7,8,9,10,11,12,13,14,15,16,17,18,19,20,21,22,23,24,25,26,
        27,28,29,30,31) ],
         
M3 = [ (2,31,24)(3,23,22)(4,21,8)(5,7,6)(9,20,10)(11,19,28)(12,27,13)(14,26,15)
       (16,25,30)(17,29,18), (1,31)(2,23)(3,21)(4,7)(8,20)(10,19)(11,27)(13,26)
       (15,25)(16,29)(18,28)(24,30), 
       (1,2,3,4,5,6,7,8,9,10,11,12,13,14,15,16,17,18,19,20,21,22,23,24,25,26,
        27,28,29,30,31) ], 
        
M4 = [ (2,31,10)(3,9,18)(4,17,5)(6,16,20)(7,19,8)(11,30,29)(12,28,27)(13,26,14)
       (15,25,21)(22,24,23), (1,31)(2,9)(3,17)(5,16)(6,19)(8,18)(10,30)(11,28)
       (12,26)(14,25)(15,20)(21,24), 
       (1,2,3,4,5,6,7,8,9,10,11,12,13,14,15,16,17,18,19,20,21,22,23,24,25,26,
        27,28,29,30,31) ], 
        
M5 = [ (2,31,22)(3,21,20)(4,19,5)(6,18,26)(7,25,24)(8,23,30)(9,29,13)(10,12,11)
       (14,28,15)(16,27,17), (1,31)(2,21)(3,19)(5,18)(6,25)(7,23)(8,29)(9,12)
       (13,28)(15,27)(17,26)(22,30), 
       (1,2,3,4,5,6,7,8,9,10,11,12,13,14,15,16,17,18,19,20,21,22,23,24,25,26,
        27,28,29,30,31) ], 
        
M6 = [ (2,31,9)(3,8,17)(4,16,15)(5,14,22)(6,21,20)(7,19,18)(10,30,11)(12,29,25)
       (13,24,23)(26,28,27), (1,31)(2,8)(3,16)(4,14)(5,21)(6,19)(7,17)(9,30)
       (11,29)(12,24)(13,22)(25,28), 
       (1,2,3,4,5,6,7,8,9,10,11,12,13,14,15,16,17,18,19,20,21,22,23,24,25,26,
        27,28,29,30,31) ]
\end{verbatim}
}

\begin{rema}\rm
We see that there are different ways to draw maps on surfaces of genus $g>0$: compare
Figures \ref{deg12map}, \ref{fig:psl-2-16}, \ref{PSL_3(5)dessin}, \ref{fig:psl-5-2}.
One may choose a representation according to his or her taste and convenience.

Also, this time, in Figure \ref{fig:psl-5-2}, we decided once again to put the labels 
explicitly. We recommend to the reader to compare them with the permutations given 
in the above triple {\tt M1} and, specifically, to ``go'' along the cycle $(1,2,\ldots,31)$.
\end{rema}


\medskip

\subsection{Group ${\rm AGL}_5(2)$, degree 32}
\label{ex:AGL_5_2_another}

Below we give six dessins representing the group ${\rm AGL}_5(2)$
(see Example~\ref{ex:deg32}).
Their passport is $(3^{10}1^2,2^{12}1^8,31^1 1^1)$, so these are alternatives to 
the dessins of degree~$32$ and passport $(3^{10}1^2, 2^{16}, 31^11^1)$
given by the construction used to prove Theorem~\ref{AGLMersenne}.

\begin{figure}[htbp]
\begin{center}
\epsfig{file=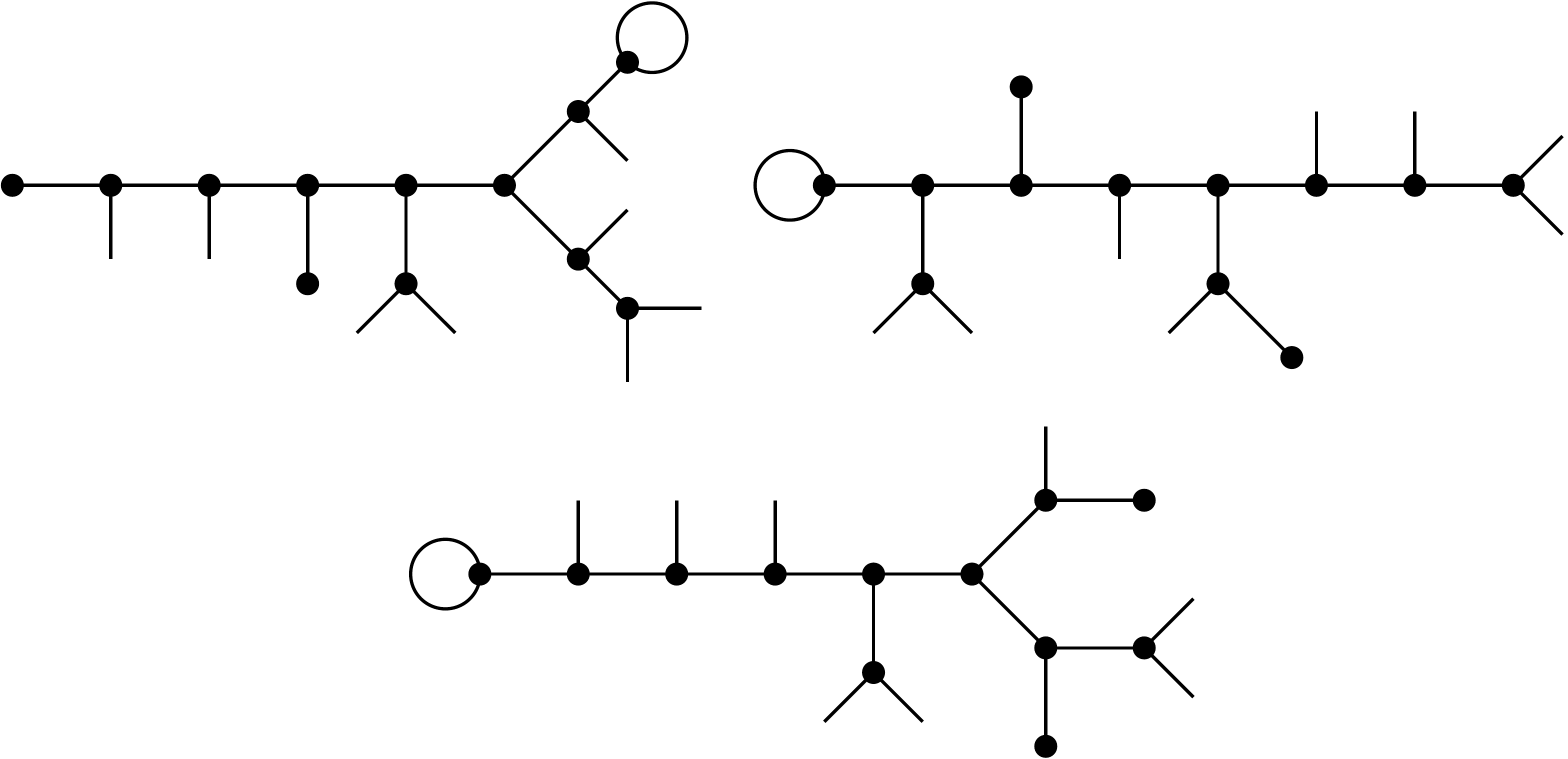,width=10.5cm}
\end{center}
\caption{The maps $\M_1,\M_2,\M_3$ with the passport $(3^{10}1^2,2^{12}1^8,31^1 1^1)$ 
representing the group ${\rm AGL}_5(2)$. Three other maps are their mirror images.}
\label{fig:asl-5-2}
\end{figure}

\begin{itemize}
\item	The permutation $z$ is always the same: $z=(1,2,3,\ldots,31)$. 
\item	The label 1 is attached to the outer half-edge of the only loop; the label of the
		face of degree~1 is 32.
\item	The chiral pairs are $(\M_1,\M_4)$, $(\M_2,\M_5)$, and $(\M_3,\M_6)$.
\end{itemize}
The first three maps are shown in Figure \ref{fig:asl-5-2}.

\bigskip

{\footnotesize
\begin{verbatim}
M1 = [ (1,2,32)(3,31,4)(5,30,12)(6,11,7)(8,10,9)(13,29,17)(14,16,15)(18,28,20)
       (21,27,22)(23,26,24), (1,32)(2,31)(4,30)(5,11)(7,10)(12,29)(13,16)(17,28)
       (18,19)(20,27)(22,26)(24,25), 
       (1,2,3,4,5,6,7,8,9,10,11,12,13,14,15,16,17,18,19,20,21,22,23,24,25,26,
        27,28,29,30,31) ], 
        
M2 = [ (1,2,32)(3,31,27)(4,26,6)(7,25,24)(8,23,18)(9,17,10)(11,16,12)(13,15,14)
       (19,22,21)(28,30,29), (1,32)(2,31)(3,26)(4,5)(6,25)(7,23)(8,17)(10,16)
       (12,15)(18,22)(19,20)(27,30), 
       (1,2,3,4,5,6,7,8,9,10,11,12,13,14,15,16,17,18,19,20,21,22,23,24,25,26,
        27,28,29,30,31) ], 
        
M3 = [ (1,2,32)(3,31,4)(5,30,6)(7,29,8)(9,28,24)(10,23,15)(11,14,12)(16,22,20)
       (17,19,18)(25,27,26), (1,32)(2,31)(4,30)(6,29)(8,28)(9,23)(10,14)(12,13)
       (15,22)(16,19)(20,21)(24,27), 
       (1,2,3,4,5,6,7,8,9,10,11,12,13,14,15,16,17,18,19,20,21,22,23,24,25,26,
        27,28,29,30,31) ], 
        
M4 = [ (1,32,31)(2,30,29)(3,28,21)(4,20,16)(5,15,13)(6,12,11)(7,10,9)(17,19,18)
       (22,27,26)(23,25,24), (1,30)(2,28)(3,20)(4,15)(5,12)(6,10)(7,8)(13,14)
       (16,19)(21,27)(22,25)(31,32), 
       (1,2,3,4,5,6,7,8,9,10,11,12,13,14,15,16,17,18,19,20,21,22,23,24,25,26,
        27,28,29,30,31) ], 
        
M5 = [ (1,32,31)(2,30,6)(3,5,4)(7,29,27)(8,26,9)(10,25,15)(11,14,12)(16,24,23)
       (17,22,21)(18,20,19), (1,30)(2,5)(6,29)(7,26)(9,25)(10,14)(12,13)(15,24)
       (16,22)(17,20)(27,28)(31,32), 
       (1,2,3,4,5,6,7,8,9,10,11,12,13,14,15,16,17,18,19,20,21,22,23,24,25,26,
        27,28,29,30,31) ], 
        
M6 = [ (1,32,31)(2,30,29)(3,28,27)(4,26,25)(5,24,9)(6,8,7)(10,23,18)(11,17,13)
       (14,16,15)(19,22,21), (1,30)(2,28)(3,26)(4,24)(5,8)(9,23)(10,17)(11,12)
       (13,16)(18,22)(19,20)(31,32), 
       (1,2,3,4,5,6,7,8,9,10,11,12,13,14,15,16,17,18,19,20,21,22,23,24,25,26,
        27,28,29,30,31) ]
\end{verbatim}
}

\subsection{Group ${\rm PSL}_3(8)$, degree 73}
\label{sec:deg-73}

Below we give the eight dessins of genus 4 with the passport $(3^{24}1^1,2^{32}1^9,73^1)$
representing the group ${\rm PSL}_3(8)$ of degree 73 (see Example~\ref{ex:psl-3-8}). 

\begin{itemize}
\item	The permutation $z$ is always the same: $z=(1,2,3,\ldots,73)$. 
\item	The fixed point of $x$ is always labelled by $1$.
\item	The chiral pairs are $(\M_1,\M_5)$, $(\M_2,\M_6)$, $(\M_3,\M_7)$ and $(\M_4,\M_8)$. 
\end{itemize}

\medskip

{\footnotesize
\begin{verbatim}
M1 = [ (2,73,18)(3,17,36)(4,35,47)(5,46,42)(6,41,58)(7,57,50)(8,49,32)(9,31,
        10)(11,30,71)(12,70,69)(13,68,21)(14,20,28)(15,27,39)(16,38,37)(19,72,
        29)(22,67,66)(23,65,64)(24,63,55)(25,54,60)(26,59,40)(33,48,34)(43,45,
        44)(51,56,62)(52,61,53), (1,73)(2,17)(3,35)(4,46)(5,41)(6,57)(7,49)(8,
        31)(10,30)(11,70)(12,68)(13,20)(14,27)(15,38)(16,36)(18,72)(19,28)(21,
        67)(22,65)(23,63)(24,54)(25,59)(26,39)(29,71)(32,48)(34,47)(40,58)(42,
        45)(50,56)(51,61)(53,60)(55,62), (1,2,3,4,5,6,7,8,9,10,11,12,13,14,15,
        16,17,18,19,20,21,22,23,24,25,26,27,28,29,30,31,32,33,34,35,36,37,38,
        39,40,41,42,43,44,45,46,47,48,49,50,51,52,53,54,55,56,57,58,59,60,61,
        62,63,64,65,66,67,68,69,70,71,72,73) ], 
        
M2 = [ (2,73,72)(3,71,28)(4,27,48)(5,47,41)(6,40,7)(8,39,22)(9,21,13)(10,12,
        11)(14,20,51)(15,50,24)(16,23,38)(17,37,59)(18,58,31)(19,30,52)(25,49,
        26)(29,70,53)(32,57,62)(33,61,34)(35,60,36)(42,46,68)(43,67,66)(44,65,
        55)(45,54,69)(56,64,63), (1,73)(2,71)(3,27)(4,47)(5,40)(7,39)(8,21)(9,
        12)(13,20)(14,50)(15,23)(16,37)(17,58)(18,30)(19,51)(22,38)(24,49)(26,
        48)(28,70)(29,52)(31,57)(32,61)(34,60)(36,59)(41,46)(42,67)(43,65)(44,
        54)(45,68)(53,69)(55,64)(56,62), (1,2,3,4,5,6,7,8,9,10,11,12,13,14,15,
        16,17,18,19,20,21,22,23,24,25,26,27,28,29,30,31,32,33,34,35,36,37,38,
        39,40,41,42,43,44,45,46,47,48,49,50,51,52,53,54,55,56,57,58,59,60,61,
        62,63,64,65,66,67,68,69,70,71,72,73) ],
         
M3 = [ (2,73,9)(3,8,51)(4,50,66)(5,65,53)(6,52,7)(10,72,71)(11,70,47)(12,46,
        13)(14,45,37)(15,36,57)(16,56,17)(18,55,28)(19,27,20)(21,26,22)(23,25,
        32)(24,31,33)(29,54,64)(30,63,34)(35,62,58)(38,44,43)(39,42,68)(40,67,
        49)(41,48,69)(59,61,60), (1,73)(2,8)(3,50)(4,65)(5,52)(7,51)(9,72)(10,
        70)(11,46)(13,45)(14,36)(15,56)(17,55)(18,27)(20,26)(22,25)(23,31)(24,
        32)(28,54)(29,63)(30,33)(34,62)(35,57)(37,44)(38,42)(39,67)(40,48)(41,
        68)(47,69)(49,66)(53,64)(58,61), (1,2,3,4,5,6,7,8,9,10,11,12,13,14,15,
        16,17,18,19,20,21,22,23,24,25,26,27,28,29,30,31,32,33,34,35,36,37,38,
        39,40,41,42,43,44,45,46,47,48,49,50,51,52,53,54,55,56,57,58,59,60,61,
        62,63,64,65,66,67,68,69,70,71,72,73) ],
         
M4 = [ (2,73,54)(3,53,10)(4,9,33)(5,32,6)(7,31,30)(8,29,34)(11,52,44)(12,43,
        65)(13,64,14)(15,63,62)(16,61,23)(17,22,18)(19,21,20)(24,60,38)(25,37,
        50)(26,49,70)(27,69,47)(28,46,35)(36,45,51)(39,59,40)(41,58,57)(42,56,
        66)(48,68,71)(55,72,67), (1,73)(2,53)(3,9)(4,32)(6,31)(7,29)(8,33)(10,
        52)(11,43)(12,64)(14,63)(15,61)(16,22)(18,21)(23,60)(24,37)(25,49)(26,
        69)(27,46)(28,34)(35,45)(36,50)(38,59)(40,58)(41,56)(42,65)(44,51)(47,
        68)(48,70)(54,72)(55,66)(67,71), (1,2,3,4,5,6,7,8,9,10,11,12,13,14,15,
        16,17,18,19,20,21,22,23,24,25,26,27,28,29,30,31,32,33,34,35,36,37,38,
        39,40,41,42,43,44,45,46,47,48,49,50,51,52,53,54,55,56,57,58,59,60,61,
        62,63,64,65,66,67,68,69,70,71,72,73) ],
         
M5 = [ (2,73,57)(3,56,46)(4,45,64)(5,63,6)(7,62,54)(8,53,9)(10,52,11)(12,51,
        20)(13,19,24)(14,23,22)(15,21,50)(16,49,35)(17,34,69)(18,68,25)(26,67,
        43)(27,42,41)(28,40,71)(29,70,33)(30,32,31)(36,48,60)(37,59,38)(39,58,
        72)(44,66,65)(47,55,61), (1,73)(2,56)(3,45)(4,63)(6,62)(7,53)(9,
        52)(11,51)(12,19)(13,23)(14,21)(15,49)(16,34)(17,68)(18,24)(20,50)(25,
        67)(26,42)(27,40)(28,70)(29,32)(33,69)(35,48)(36,59)(38,58)(39,71)(43,
        66)(44,64)(46,55)(47,60)(54,61)(57,72), 
        (1,2,3,4,5,6,7,8,9,10,11,12,13,14,15,16,17,18,19,20,21,22,23,24,25,26,
        27,28,29,30,31,32,33,34,35,36,37,38,39,40,41,42,43,44,45,46,47,48,49,50,
        51,52,53,54,55,56,57,58,59,60,61,62,63,64,65,66,67,68,69,70,71,72,73) ], 
        
M6 = [ (2,73,3)(4,72,47)(5,46,22)(6,21,30)(7,29,33)(8,32,9)(10,31,20)(11,19,
        12)(13,18,43)(14,42,41)(15,40,39)(16,38,58)(17,57,44)(23,45,56)(24,55,
        61)(25,60,51)(26,50,49)(27,48,71)(28,70,34)(35,69,68)(36,67,53)(37,52,
        59)(54,66,62)(63,65,64), (1,73)(3,72)(4,46)(5,21)(6,29)(7,32)(9,
        31)(10,19)(12,18)(13,42)(14,40)(15,38)(16,57)(17,43)(20,30)(22,45)(23,
        55)(24,60)(25,50)(26,48)(27,70)(28,33)(34,69)(35,67)(36,52)(37,58)(44,
        56)(47,71)(51,59)(53,66)(54,61)(62,65), 
        (1,2,3,4,5,6,7,8,9,10,11,12,13,14,15,16,17,18,19,20,21,22,23,24,25,26,
        27,28,29,30,31,32,33,34,35,36,37,38,39,40,41,42,43,44,45,46,47,48,49,50,
        51,52,53,54,55,56,57,58,59,60,61,62,63,64,65,66,67,68,69,70,71,72,73) ], 
        
M7 = [ (2,73,66)(3,65,4)(5,64,28)(6,27,34)(7,33,36)(8,35,26)(9,25,71)(10,70,
        22)(11,21,46)(12,45,41)(13,40,17)(14,16,15)(18,39,60)(19,59,58)(20,57,
        47)(23,69,68)(24,67,72)(29,63,62)(30,61,38)(31,37,32)(42,44,51)(43,50,
        52)(48,56,55)(49,54,53), (1,73)(2,65)(4,64)(5,27)(6,33)(7,35)(8,25)(9,
        70)(10,21)(11,45)(12,40)(13,16)(17,39)(18,59)(19,57)(20,46)(22,69)(23,
        67)(24,71)(26,34)(28,63)(29,61)(30,37)(32,36)(38,60)(41,44)(42,50)(43,
        51)(47,56)(48,54)(49,52)(66,72), (1,2,3,4,5,6,7,8,9,10,11,12,13,14,15,
        16,17,18,19,20,21,22,23,24,25,26,27,28,29,30,31,32,33,34,35,36,37,38,
        39,40,41,42,43,44,45,46,47,48,49,50,51,52,53,54,55,56,57,58,59,60,61,
        62,63,64,65,66,67,68,69,70,71,72,73) ], 
        
M8 = [ (2,73,21)(3,20,8)(4,7,27)(5,26,49)(6,48,28)(9,19,33)(10,32,63)(11,62,
        61)(12,60,13)(14,59,52)(15,51,37)(16,36,35)(17,34,18)(22,72,65)(23,64,
        31)(24,30,39)(25,38,50)(29,47,40)(41,46,67)(42,66,71)(43,70,69)(44,68,
        45)(53,58,57)(54,56,55), (1,73)(2,20)(3,7)(4,26)(5,48)(6,27)(8,19)(9,
        32)(10,62)(11,60)(13,59)(14,51)(15,36)(16,34)(18,33)(21,72)(22,64)(23,
        30)(24,38)(25,49)(28,47)(29,39)(31,63)(37,50)(40,46)(41,66)(42,70)(43,
        68)(45,67)(52,58)(53,56)(65,71), (1,2,3,4,5,6,7,8,9,10,11,12,13,14,15,
        16,17,18,19,20,21,22,23,24,25,26,27,28,29,30,31,32,33,34,35,36,37,38,
        39,40,41,42,43,44,45,46,47,48,49,50,51,52,53,54,55,56,57,58,59,60,61,
        62,63,64,65,66,67,68,69,70,71,72,73) ]
\end{verbatim}
}

\bigskip

The next candidate for our series of examples would be the group ${\rm PSL}_7(2)$ 
of degree 127. Since we don't see any particular interest in writing explicitly 
the generating permutations of that large degree, we stop here.

%% file: Klein-11-Second-Version.bbl
\medskip

\begin{thebibliography}{99}


\bibitem{APZ} N.~M.~Adrianov, F.~Pakovich and A.~K.~Zvonkin, {\em Davenport--Zannier 
	Polynomials and Dessins d'Enfants}, AMS Mathematical Surveys and Monographs 249,
	Providence, RI (2020). %

\bibitem{AFG} S.~L.~Aletheia-Zomlefer, L.~Fukshansky and S.~R.~Garcia, ``The Bateman--Horn 
	conjecture: heuristics, history, and applications'', {\em Expo. Math.}, {\bf 38}, 
	430--479 (2020). Also available at \url{arXiv-math[NT]:1807.08899v4}. %

\bibitem{BH} P.~T.~Bateman and R.~A.~Horn, ``A heuristic asymptotic formula concerning the 
	distribution of prime numbers'', {\em Math.~Comp.}, {\bf 16}, 220--228 (1962). %

\bibitem{Bell} G.~W.~Bell, ``On the cohomology of the finite special linear groups I, II'', 
	{\em J.~Algebra}, {\bf 54}, 216--238, 239--259 (1978). %

\bibitem{betrema} J. B\'etr\'ema, Private communication (2020).

\bibitem{Bel} G.~V.~Bely\u\i, ``On Galois extensions of a maximal cyclotomic field'', 
	{\em Izv. Akad. Nauk SSSR Ser. Mat.}, {\bf 43}, 267--276, 479 (1979) (Russian). 
	English translation: {\em Math. USSR Izvestya}, {\bf 14}, 247--256 (1980). %

\bibitem{Ber} \'A.~Berecky, ``Maximal overgroups of Singer elements in classical groups'', 
	{\em J.~Algebra}, {\bf 234}, 187--206 (2000). %

\bibitem{BS} Z.~I.~Borevich and I.~R.~Shafarevich, {\em Number Theory}, Academic Press, 
	New York and London (1966). %

\bibitem{Bun-1857} V. Bouniakowsky, ``Sur les diviseurs num\'eriques invariables des 
	fonctions rationnelles enti\`eres'', {\em M\'em. Acad. Sci. St. P\'eteresbourg}, 
	$6^{\rm e}$ s\'erie, vol. {\bf VI}, 305--329 (1857).\footnote{Numerous publications
	give the following wrong title for Bunyakovsky's paper: ``Nouveaux th\'eor\`emes 
	relatifs \`a la distinction des nombres premiers et \`a la d\'ecomposition des 
	entiers en facteurs''. According to the French Wikipedia (see~\cite{Bun-wiki}), 
	an article with this title does indeed exist, but it was published in 1840 and 
	not in 1857, and it does not discuss the conjecture in question. The reader 
	may also consult the original paper reproduced in the Google archive.} Available at
	the Google archive: 
	\url{https://books.google.fr/books?hl=fr&id=wXIhAQAAMAAJ&pg=PA305#v=onepage&q&f=false}. %
	
\bibitem{Bun-wiki} Bunyakovsky conjecture in Wikipedia: \\ 
	English: {\small \url{https://en.wikipedia.org/wiki/Bunyakovsky_conjecture}}; \\
	French: {\small \url{https://fr.wikipedia.org/wiki/Conjecture_de_Bouniakovski}}. %

\bibitem{Bur06} W.~Burnside, ``On simply transitive groups of prime degree'', 
	{\em Quart. J. Math. (Oxford)}, {\bf 37}, 215--221 (1906). %

\bibitem{Bur} W.~Burnside, {\em Theory of Groups of Finite Order}\/ (2nd ed.), 
	Cambridge University Press, Cambridge (1911). Reprinted by Dover, NewYork (1955). %
	
\bibitem{Cam} P.~J.~Cameron, {\em Permutation Groups}, London Mathematical Society 
	Student Texts 45, Cambridge University Press, Cambridge (1999). %

\bibitem{Coh} H.~Cohen, {\em A Course in Computational Algebraic Number Theory}, Springer, 
	Berlin -- Heidelberg (1993). %

\bibitem{Con91} M.~D.~E.~Conder, ``The symmetric genus of the Mathieu groups'', 
	{\em Bull.~London                                                                                                                                                                                                                                                                                                                                                                                                      Math.~Soc.}, {\bf 23}, 445--453 (1991). %

\bibitem{Con} M.~D.~E.~Conder, ``Regular maps and hypermaps of Euler characteristic 
	$-1$ to $-200$'', {\em J. Combin. Theory, Ser.~B}, {\bf 99}, 455--459 (2009) 
	(with the associated lists of computational data available at 
	\url{http://www.math.auckland.ac.nz/~conder/hypermaps.html}). %

\bibitem{ATLAS} J.~H.~Conway, R.~T.~Curtis, S.~P.~Norton, R.~A.~Parker and R.~A.~Wilson, 
	{\em ATLAS of Finite Groups\/}, Clarendon Press, Oxford (1985). 
	Reprinted with corrections in~2005. %

\bibitem{Cox84} H.~S.~M.~Coxeter, ``A symmetrical arrangement of eleven hemi-icosahedra'', 
	in: {\em Convexity and Graph Theory (Jerusalem, 1981)}, North-Holland Math.~Stud.~87, 
	Ann. Discrete Math.~20, North-Holland, Amsterdam (1984), pp. 103--114. 

\bibitem{CE} H.~S.~M.~Coxeter and W.~L.~Edge, ``The simple groups ${\rm PSL}(2,7)$ 
	and ${\rm PSL}(2,11)$'', {\em C.~R. Math. Rep. Acad. Sci. Canada}, {\bf 5}, 
	201--206 (1983). %

\bibitem{CM} H.~S.~M.~Coxeter and W.~O.~J.~Moser, {\em Generators and Relations for 
	Discrete Groups}, 4th ed., Springer, Berlin -- Heidelberg -- New York (1980). %

\bibitem{CW} H.~S.~M.~Coxeter and A.~I.~Weiss, ``Twisted honeycombs $\{3,5,3\}_t$ and their 
	groups'', {\em Geometriae Dedicata}, {\bf 17}, 169--179 (1984). 

\bibitem{Dic} L.~E.~Dickson, {\em Linear Groups}, Dover, New York (1958). %

\bibitem{DM} J.~D.~Dixon and B.~Mortimer, {\em Permutation Groups}, Graduate Texts in 
	Mathematics 163, Springer, New York (1996). %

\bibitem{Dyck-1888} W. Dyck, ``\"Uber das Problem der Nachbargebiete'', {\em Math. Ann.}, 
	{\bf 32}, 457--512 (1888). %

\bibitem{Fro} F.~G.~Frobenius, ``\"Uber Gruppencharaktere'', {\em Sitzber. K\"oniglich 
	Preuss. Akad. Wiss.~Berlin}, 985--1021 (1896).%

\bibitem{Gal} \'E.~Galois, ``M\'emoire sur les conditions de r\'esolubilit\'e des 
	\'equations par radicaux'', {\em J. Math. Pures et Appl.}, 
	vol. {\bf XI}, 417--433 (1846). Available at: 
\url{http://www.bibnum.education.fr/sites/default/files/galois_memoire_sur_la_resolubiblite.pdf}. 
\\
	See also: {\em Ecrits et M\'emoires Math\'ematiques d'\'Evariste Galois}, R.~Bourgne and 
	J.-P.~Azra, editors, Gauthier-Villars, Paris (1962). \\
	English translation: ``Memoir on the conditions for solvability of equations by
	radicals'', in: H.~M.~Edwards, {\em Galois Theory}, Springer (Graduate Texts
	in Mathematics, vol. {\bf 101}), pp. 101--113 (1984). \\
	See also: P.~M.~Neumann, {\em The Mathematical Writings of \'Evariste Galois}, 
	European Math.~Soc., Zurich, Chapter~IV (2011). %
	
\bibitem{last-letter} \'E.~Galois, ``Lettre de Galois \`a M.~August Chevalier'',
	{\em J. Math. Pures Appl.}, vol. {\bf XI}, 408--415 (1846).
	Available at: \url{http://visualiseur.bnf.fr/CadresFenetre?O=NUMM-16390}). %

\bibitem{GAP}  The GAP Group, GAP -- Groups, Algorithms, and Programming, 
	\url{http://www.gap-system.org}. %

\bibitem{GM} F.~W.~Gehring and G.~J.~Martin, ``Minimal co-volume hyperbolic lattices. 
	I.~The spherical points of a Kleinian group'', {\em Ann.~of Math.} (2), {\bf 170},
	123--161 (2009). 

\bibitem{GG} E.~Girondo and G.~Gonz\'alez-Diez, {\em Introduction to Compact Riemann 
	Surfaces and Dessins d'Enfants}, London Math. Soc. Student Texts 79, Cambridge 
	Univ. Press, Cambridge (2012). %

\bibitem{GIR} C.~Godsil, W.~Imrich and R.~Razen, ``On the number of subgroups of given 
	index in the modular group'', {\em Monatsh.~Math.}, {\bf 867}, 273--280 (1979). %

\bibitem{Goo} R. Goormaghtigh, {\em L'Interm\'ediaire des Math\'ematiciens}, {\bf 24}, 
	88 (1917). %

\bibitem{Goo-wiki} Goormaghtigh conjecture in Wikipedia: 
	\url{https://en.wikipedia.org/wiki/Goormaghtigh_conjecture}. %
	
\bibitem{Gre} J.~A.~Green, ``The characters of the finite general linear groups'', 
	{\em Trans. Amer. Math. Soc.}, {\bf 80}, 402--447 (1955). %

\bibitem{Gro} A.~Grothendieck, ``Esquisse d'un Programme'', in: {\em Geometric Galois 
	Actions~$1$. Around Grothendieck's Esquisse d'un Programme}, P.~Lochak and L.~Schneps, 
	editors, London Math. Soc. Lecture Note Ser., vol. {\bf 242}, Cambridge Univ. Press, 
	Cambridge, pp.~5--48 (1997). 
	English translation: ``Sketch of a programme'', the same volume, pp.~243--283.
	(Original publication appeared as a preprint in 1984.) %

\bibitem{Gru} B.~Gr\"unbaum, ``Regularity of graphs, complexes and designs'', in: 
	{\em Probl\`emes Combinatoires et Th\'eorie des Graphes}, Coll. Int. CNRS, {\bf 260}, 
	191--197 (1977). %

\bibitem{Hamilton-1856} W. R. Hamilton, ``On the Icosian'', Letter to John T. Graves
	(17th October 1856). Reprinted in: W. R. Hamilton, {\em Mathematical Papers, vol. {\rm III}: 
	Algebra}, H. Halberstam and R. E. Ingram, editors, Cambridge Univ. Press, Cambridge, 
	\linebreak pp.~612--625 (1967). %
		
\bibitem{Heffter-1891} L. Heffter, ``\"Uber das Problem der Nachbargebiete'',
	{\em Math. Ann.}, {\bf 38}, 477--508 (1891). %
		
\bibitem{Hup} B.~Huppert, {\em Endliche Gruppen}\/~I, Springer, Berlin -- Heidelberg -- 
	New York (1979). %

\bibitem{ISS} I.~Ivrissimtzis, D.~Singerman and J.~Strudwick, ``From Farey fractions to the 
	Klein quartic and beyond'', {\em Ars Math.~Contemp.}, {\bf 20}, 37--50 (2021). 
	Available at \url{https://arXiv.math:1909.08568} and 
	\url{https://doi.org/10.26493/1855-3974.2046.cb6}

\bibitem{Jon14} G.~A.~Jones, ``Primitive permutation groups containing a cycle'', 
	{\em Bull. Aust. Math. Soc.} {\bf 89}, 159--165 (2014). %

\bibitem{JS} G.~A.~Jones and D.~Singerman, {\em Complex Functions: an Algebraic and 
	Geometric Viewpoint}, Cambridge Univ. Press, Cambridge (1987). %

\bibitem{JLM} G.~A.~Jones, C.~D.~Long and A.~D.~Mednykh, ``Hyperbolic manifolds and 
	tessellations of type $\{3,5,3\}$ associated with ${\rm L}_2(q)$'', 
	\url{arXiv.math: 1106.0867}. %

\bibitem{JW} G.~A.~Jones and J.~Wolfart, {\em Dessins d'Enfants on Riemann Surfaces}, 
	Springer, Cham (2016). %

\bibitem{JZ} G.~A.~Jones and A.~K.~Zvonkin, ``Hurwitz groups as monodromy groups of dessins: 
	several examples'', \url{https://arxiv.org/pdf/2012.07107.pdf} (2020). To appear in:
	{\em Teichm\"uller Theory and its Impact}, Higher Education Press 
	and International Press, Beijing. %

\bibitem{JZ20} G.~A.~Jones and A.~K.~Zvonkin, ``Primes in geometric series and finite 
	permutation groups'', \url{https://arxiv.org/pdf/2010.08023.pdf} (2020). %

\bibitem{Jor71} C.~Jordan, ``Th\'eor\`emes sur les groupes primitifs'', 
	{\em J. Math. Pures Appl.}\/ (2), {\bf 16}, 383--408 (1871). %

\bibitem{Jor73} C.~Jordan, ``Sur la limite de transitivit\'e des groups non altern\'es'', 
	{\em Bull. Soc. Math. France}, {\bf 1}, 40--71 (1873).. %

\bibitem{Kle} F.~Klein, ``\"Uber die Erniedrigung der Modulgleichungen'', {\em Math.~Ann.}, 
	{\bf 14}, No.~3, 417--427 (1878). %

\bibitem{Kle78} F.~Klein, ``\"Uber die Transformationen siebenter Ordnung der
	elliptischen Funktionen'', {\em Math.~Ann.}, {\bf 14}, No.~3, 428--471 (1878). %

\bibitem{Kle79} F.~Klein, ``\"Uber die Transformationen elfter Ordnung der
	ellipti\-schen Funktionen'', {\em Math.~Ann.}, {\bf 15}, No.~3--4, 533--555, (1879). %
		
\bibitem{Kle84} F.~Klein, {\em Vorlesungen \"uber das Ikosaeder und die Aufl\"osung der 
	Gleichungen vom f\"unften Grade}, Teubner, Leipzig (1884). 
	English translation: {\em Lectures on the Icosahedron and the Solution
	of Equations of the Fifth Degree}, Dover Publications (2003). 

\bibitem{KleGMA} F.~Klein, {\em Gesammelte Mathematische Abhandlungen}\/ (3 vols.), Springer, 
	Berlin (1923--1973). %

\bibitem{LZ} S.~K.~Lando and A.~K.~Zvonkin, {\em Graphs on Surfaces and their Applications}, 
	Encyclopedia of Math.~Sciences, vol. {\bf 141}, Springer, Berlin (2004). %

\bibitem{Lev} S.~Levy (ed.), {\em The Eightfold Way. The Beauty of Klein's Quartic Curve}, 
	Cambridge Univ.~Press, Cambridge (1999). %

\bibitem{LRS} J.~B.~Lewis, V.~Reiner and D.~Stanton, ``Reflection factorizations of Singer 
	cycles'', {\em J. Algebr. Combin.}, {\bf 40}, 663--691 (2014).%

\bibitem{W.Li} W.~Li, ``A note on the Bateman--Horn conjecture'',
	{\em J.~Number Theory}, {\bf 208}, 390--399 (2020).
	Also available at \url{https://arxiv.org/pdf/1906.03370.pdf}. %

\bibitem{Macb} A.~M.~Macbeath, ``Generators of the linear fractional groups'', in:
	{\em Number Theory (Houston\/ $1967$)}, W.~J.~Leveque and E.~G.~Straus, editors, 
	{\em Proc. Sympos. Pure Math.}, {\bf 12}, Amer. Math. Soc., Providence, RI, 14--32 (1969). %

\bibitem{Macd} I.~G.~Macdonald, {\em Symmetric Functions and Hall Polynomials},  
	Oxford Math. Monographs, Clarendon Press, Oxford (1979). %

\bibitem{McMS} P.~McMullen and E.~Schulte, {\em Abstract Regular Polytopes}, Encyclopedia 
	of Mathematics and its Applications, vol.~{\bf  92}, Cambridge University Press, 
	Cambridge (2002). 

\bibitem{MM} T.~H.~Marshall and G.~J.~Martin, ``Minimal co-volume hyperbolic lattices. {\rm II}. 
	Simple torsion in a Kleinian group'', {\em Ann.~of Math.} (2), {\bf 176}, 261--301 (2012). %

\bibitem{Mat61} \'E.~Mathieu, ``M\'emoire sur l'\'etude des fonctions de plusieurs quantit\'es, 
	sur la mani\`ere de les former et sur les substitutions qui les laissent invariables'', 
	{\em J. Math. Pures Appl.}, {\bf 6}, 241--323 (1861). %

\bibitem{Mat73} \'E.~Mathieu, ``Sur la fonction cinq fois transitive de $24$ quantit\'es'', 
	{\em J.~Math.~Pures Appl.}, {\bf 18}, 25--46 (1873). %

\bibitem{McQ} D.~L.~McQuillan, ``Classification of normal congruence subgroups of the modular 
	group'', {\em Amer.~J.~Math.}, {\bf 87}, 285--296 (1965). %

\bibitem{Mul} P.~M\"uller, ``Reducibility behavior of polynomials with varying coefficients'', 
	{\em Israel J.~Math.}, {\bf 94}, 59--91 (1996). %

\bibitem{NZM} I.~Niven, H.~S.~Zuckerman and H.~L.~Montgomery, {\em An Introduction to 
	the Theory of Numbers} (5th ed.), Wiley, New York (1991). %
	
\bibitem{OEIS} Online Encyclopedia of Integer Sequences, \url{https://oeis.org/}. %

\bibitem{SS} A.~Schinzel and W.~Sierpi\'nski, ``Sur certaines hypoth\`eses concernant les 
	nombres premiers'', {\em Acta Arith.}, {\bf 4}, 185--298 (1958).
	Erratum: {\bf 5}, 259 (1958). %

\bibitem{SF} W.~A.~Simpson and J.~S.~Frame, ``The character tables for ${\rm SL}(3,q)$, 
	${\rm SU}(3,q^2)$, ${\rm PSL}(3,q)$, ${\rm PSU}(3,q^2)$'' {\em Canadian J.~Math.},
	{\bf 25}, 486--494, (1973). %

\bibitem{Sin74}  D.~Singerman, ``Symmetries of Riemann surfaces with large automorphism 
	group'', {\em Math.~Ann.}, {\bf 210}, 17--32 (1974). %

\bibitem{Walfisz} A. Walfisz, ``Zur additiven Zahlentheorie. {\rm II}'', {\em Math.~Z.},
	{\bf 40}, 592--607 (1936). %

\bibitem{Wie} H.~Wielandt, {\em Finite Permutation Groups}, Academic Press, New York (1964). %
 
\bibitem{Wol90} A.~J.~Woldar, ``Representing ${\rm M}_{11}$, ${\rm M}_{12}$, ${\rm M}_{22}$ 
	and ${\rm M}_{23}$  on surfaces of least genus'', {\em Comm. Algebra}, {\bf 18}, 15--86
	(1990). Corrigendum: p.~605. %

\end{thebibliography}
